\documentclass{article}
\usepackage{amsmath}
\usepackage{amsfonts,amssymb,xspace,amsthm}

\newcommand{\lspan}{\operatorname{span}}
\newcommand{\Xtil}{\tilde{X}}
\newcommand{\Btil}{\tilde{B}}

\newcommand{\fg}{\mathfrak{g}}
\newcommand{\fh}{\mathfrak{h}}

\newcommand{\Mhh}{\hat{M}
\hspace{-1.05ex}\hat{\rule{0ex}{2.0ex}}\hspace{1.05ex}}
\newcommand{\dehh}{\hat{\de}
\hspace{-.95ex}\hat{\rule{0ex}{2.05ex}}\hspace{.95ex}}

\newcommand{\cC}{\mathcal{C}}
\newcommand{\acts}{\triangleright}
\newcommand{\beh}{\hat{\be}}
\newcommand{\hchi}{\raisebox{.45ex}[0pt][0pt]{$\chi$}}
\newcommand{\gtwo}{\mathfrak{g}_2}
\newcommand{\gone}{\mathfrak{g}_1}
\newcommand{\cA}{\mathcal{A}}
\newcommand{\Wtil}{\tilde{W}}
\newcommand{\cros}{M \, \mbox{$_{\al,\cU}$}\hspace{-.2ex}\mbox{$\ltimes$} \, N}
\newcommand{\kruisje}[1]{\, \mbox{$_{#1}$}\hspace{-.2ex}\mbox{$\ltimes$}
\,}
\newcommand{\strong}{\mbox{$\si$-strong$^*$}\xspace}
\newcommand{\Vtil}{\tilde{V}}
\newcommand{\alh}{\hat{\alpha}}
\newcommand{\Mh}{\hat{M}}
\newcommand{\Si}{\Sigma}
\newcommand{\cO}{\mathcal{O}}
\newcommand{\dpr}{^{\prime\prime}}
\newcommand{\pite}{\pi_\theta}
\newcommand{\late}{\Lambda_\theta}
\newcommand{\te}{\theta}
\newcommand{\tetil}{\tilde{\theta}}
\newcommand{\latil}{\tilde{\Lambda}}
\newcommand{\Nte}{\cN_\theta}
\newcommand{\Tr}{{\operatorname{Tr}}}
\newcommand{\nsf}{n.s.f.\xspace\!}
\newcommand{\Jtil}{\tilde{J}}
\newcommand{\Jte}{J_\te}
\newcommand{\alhh}{\hat{\alpha}\hspace{-.7ex}\hat{\rule{0ex}{1.45ex}}\hspace{.7ex}}
\newcommand{\alhhklein}{\hat{\alpha}\hspace{-.55ex}\hat{\rule[-.4ex]{0ex}{1.45ex}}\hspace{.55ex}}
\newcommand{\Mo}{M_1}
\newcommand{\Mt}{M_2}
\newcommand{\deo}{\de_1}
\renewcommand{\det}{\de_2}
\newcommand{\deth}{\hat{\de}_2}
\newcommand{\Mth}{\hat{M}_2}
\newcommand{\deoh}{\hat{\de}_1}
\newcommand{\Moh}{\hat{M}_1}

\newcommand{\Wh}{\hat{W}}
\newcommand{\Wht}{\hat{W}_2}
\newcommand{\tautil}{\tilde{\tau}}
\newcommand{\cVtil}{\tilde{\cV}}
\newcommand{\cUtil}{\tilde{\cU}}
\newcommand{\Nfih}{\cN_{\vfih}}
\newcommand{\sluit}{{- \; \strong}}
\newcommand{\lrecht}{\longrightarrow}
\newcommand{\mult}{\operatorname{m}}
\newcommand{\eh}{\otimes_{\mbox{\scriptsize eh}}}
\newcommand{\Nh}{\hat{N}}
\newcommand{\cWh}{{\hat{\cW}}}
\newcommand{\Jn}{J_N}
\newcommand{\Jhn}{\hat{J}_N}
\newcommand{\Ytil}{\tilde{Y}}
\newcommand{\Mtil}{\tilde{M}}
\newcommand{\Phitil}{\tilde{\Phi}}
\newcommand{\betil}{\tilde{\be}}
\newcommand{\detil}{\tilde{\de}}
\newcommand{\Wc}{\check{W}}
\newcommand{\vfitil}{\tilde{\vfi}}
\newcommand{\cR}{\mathcal{R}}
\newcommand{\sdeo}{\delta_1}
\newcommand{\sdet}{\delta_2}
\newcommand{\Xh}{\hat{X}}
\newcommand{\itil}{\tilde{i}}
\newcommand{\jtil}{\tilde{j}}
\newcommand{\lahh}{\hat{\la}
\hspace{-.8ex}\hat{\rule{0ex}{2.0ex}}\hspace{.8ex}}
\newcommand{\vfihh}{\hat{\vfi}
\hspace{-.58ex}\hat{\rule{0ex}{1.45ex}}\hspace{.58ex}}
\newcommand{\klvfihh}{\hat{\vfi}\hspace{-.52ex}\hat{\rule{0ex}{1ex}}\hspace{.52ex}}
\newcommand{\sihh}{\hat{\si}
\hspace{-.65ex}\hat{\rule{0ex}{1.45ex}}\hspace{.65ex}}
\newcommand{\Jhh}{\hat{J}
\hspace{-.4ex}\hat{\rule{0ex}{2.0ex}}\hspace{.4ex}}
\newcommand{\Rh}{\hat{R}}
\newcommand{\Ma}{M_a}
\newcommand{\Mb}{M_b}
\newcommand{\dea}{\de_a}
\newcommand{\deb}{\de_b}
\newcommand{\ala}{\al_a}
\newcommand{\alb}{\al_b}
\newcommand{\bea}{\be_a}
\newcommand{\beb}{\be_b}
\newcommand{\Mah}{\Mh_a}
\newcommand{\Mbh}{\Mh_b}
\newcommand{\deah}{\deh_a}
\newcommand{\debh}{\deh_b}
\newcommand{\taua}{\tau_a}
\newcommand{\taub}{\tau_b}
\newcommand{\cUa}{\cU_a}
\newcommand{\cUb}{\cU_b}
\newcommand{\cVa}{\cV_a}
\newcommand{\cVb}{\cV_b}
\newcommand{\Wa}{W_a}
\newcommand{\Wb}{W_b}
\newcommand{\alha}{\hat{\alpha}_a}
\newcommand{\alhb}{\hat{\alpha}_b}
\newcommand{\Wha}{\Wh_a}
\newcommand{\Whas}[1]{\Wh_{a, #1}}
\newcommand{\Whb}{\Wh_b}
\newcommand{\Whbs}[1]{\Wh_{b, #1}}
\newcommand{\Wtila}{\Wtil_a}
\newcommand{\Wtilas}[1]{\Wtil_{a, #1}}
\newcommand{\Wtilb}{\Wtil_b}
\newcommand{\Wtilbs}[1]{\Wtil_{b, #1}}
\newcommand{\vfia}{\vfi_a}
\newcommand{\vfib}{\vfi_b}
\newcommand{\laa}{\la_a}
\newcommand{\lab}{\la_b}

\newcommand{\R}{\mathbb{R}}
\newcommand{\Z}{\mathbb{Z}}
\newcommand{\C}{\mathbb{C}}

\newcommand{\Ga}{\Gamma}
\newcommand{\deh}{\hat{\Delta}}

\newcommand{\vfih}{\hat{\vfi}}
\newcommand{\tauh}{\hat{\tau}}

\newcommand{\lah}{\hat{\la}}
\newcommand{\Gah}{\hat{\Ga}}
\newcommand{\pih}{\hat{\pi}}
\newcommand{\sih}{\hat{\si}}
\newcommand{\sdeh}{\hat{\delta}}
\newcommand{\cU}{{\cal U}}

\newcommand{\cJ}{{\cal J}}
\newcommand{\nab}{\nabla}
\newcommand{\Jh}{\hat{J}}
\newcommand{\nabh}{\hat{\nab}}
\newcommand{\cI}{{\cal I}}
\newcommand{\cL}{{\cal L}}
\newcommand{\cF}{{\cal F}}

\newcommand{\cH}{{\cal H}}
\newcommand{\cN}{{\cal N}}
\newcommand{\cM}{{\cal M}}

\newcommand{\cD}{{\cal D}}

\newcommand{\ot}{\otimes}

\newcommand{\la}{\Lambda}
\newcommand{\Om}{\Omega}
\newcommand{\om}{\omega}
\newcommand{\io}{\iota}
\newcommand{\vfi}{\varphi}
\newcommand{\eps}{\varepsilon}

\newcommand{\al}{\alpha}
\newcommand{\be}{\beta}
\newcommand{\ga}{\gamma}
\newcommand{\sde}{\delta}
\newcommand{\de}{\Delta}

\newcommand{\si}{\sigma}
\newcommand{\Mfi}{{\cal M}_{\vfi}}
\newcommand{\Nfi}{{\cal N}_{\vfi}}
\newcommand{\Mpsi}{{\cal M}_{\psi}}
\newcommand{\Npsi}{{\cal N}_{\psi}}

\newcommand{\cV}{{\cal V}}
\newcommand{\cW}{{\cal W}}

\newcommand{\cT}{{\cal T}}

\newcommand{\Th}{\Theta}

\newcommand{\cst}{\text{C}$\hspace{0.1mm}^*$}

\newcommand{\deop}{\de \hspace{-.3ex}\raisebox{0.9ex}[0pt][0pt]{\scriptsize\fontshape{n}\selectfont op}}
\newcommand{\deopo}{\de_1 \hspace{-1ex}\raisebox{1ex}[0pt][0pt]{\scriptsize\fontshape{n}\selectfont op}}
\newcommand{\deopt}{\de_2 \hspace{-1ex}\raisebox{1ex}[0pt][0pt]{\scriptsize\fontshape{n}\selectfont op}}

\newcommand{\dehop}{\deh \hspace{-.3ex}\raisebox{0.9ex}[0pt][0pt]{\scriptsize\fontshape{n}\selectfont op}}
\newcommand{\dehopo}{\deh_1 \hspace{-1ex}\raisebox{1ex}[0pt][0pt]{\scriptsize\fontshape{n}\selectfont op}}
\newcommand{\dehopt}{\deh_2 \hspace{-1ex}\raisebox{1ex}[0pt][0pt]{\scriptsize\fontshape{n}\selectfont op}}
\newcommand{\Gaop}{\Ga\Op}

\newcommand{\Op}{\raisebox{0.9ex}[0pt][0pt]{\scriptsize\fontshape{n}\selectfont op}\,}

\newcommand{\recht}{\rightarrow}
\newcommand{\tekst}[1]{\;\;\text{#1}\;\;}

{\theoremstyle{definition}\newtheorem{definition}{Definition}[section]
\newtheorem{notation}[definition]{Notation}
\newtheorem{terminology}[definition]{Terminology}
\newtheorem{remark}[definition]{Remark}}
\newtheorem{proposition}[definition]{Proposition}
\newtheorem{lemma}[definition]{Lemma}
\newtheorem{theorem}[definition]{Theorem}
\newtheorem{corollary}[definition]{Corollary}

\numberwithin{equation}{section}

\begin{document}

\begin{center}
\Large\bf Extensions of locally compact quantum groups and the bicrossed product construction
\end{center}

\bigskip

\begin{center}
Stefaan Vaes\footnote{Research Assistant of the
Fund for Scientific Research -- Flanders (Belgium)\ (F.W.O.)} \\
Department of Mathematics \\
K.U.Leuven \\
Celestijnenlaan 200 B \\
B-3001 Leuven \\
Belgium

\medskip

e-mail : Stefaan.Vaes@wis.kuleuven.ac.be

\bigskip

Leonid Vainerman \\
Universit\'e Louis Pasteur Strasbourg \\
D\'epartement de Math\'ematiques \\
7, rue Ren\'e Descartes \\
F-67084 Strasbourg Cedex \\
France

\medskip

e-mail : vaynerma@math.u-strasbg.fr
\end{center}

\begin{abstract}
\noindent In the framework of locally compact quantum groups, we study cocycle actions. We develop the
cocycle bicrossed product construction, starting from a matched pair of locally compact quantum
groups. We define exact sequences and establish a one-to-one correspondence between cocycle
bicrossed products and cleft extensions. In this way, we obtain new examples of locally compact
quantum groups.
\end{abstract}

\section*{Introduction}

The major motivation for our paper is the fundamental work \cite{Kac} of G.I.~Kac on
extensions of finite groups which are, in modern terms, finite dimensional Kac algebras (or
Hopf $*$-algebras). Being invented in the early sixties (see \cite{Kac1}) in order to explain
in a symmetric way the duality for locally compact (l.c.) groups, Kac algebras gave
historically the first wide class of quantum groups that included besides usual groups and
their duals also non-trivial (i.e., non-commutative and non-cocommutative) objects \cite{KP1},
\cite {KP2}. In fact, one of the goals of \cite{Kac} was to give a systematic approach to the
construction of such objects. The general theory of Kac algebras was completed independently
on the one hand by G.I.~Kac and the second author \cite{KVai} and on the other hand by
M.~Enock and J.-M.~Schwartz (for a survey see \cite{E-S}).

After the appearence of quantum groups in the eighties a lot of efforts were spent in order to
construct their general theory in operator algebraic framework that would cover known examples
and that would be as elegant and symmetric as the one of Kac algebras. Important steps in this
direction were made by S.~Baaj and G.~Skandalis \cite{B-S1}, S.L.~Woronowicz \cite{Wor1,Wor2}, T.~Masuda and Y. Nakagami \cite{Mas-Nak} and A.~Van Daele \cite{VanDaele1}. The general theory of l.c.\ quantum groups was
proposed by J.~Kustermans and the first author \cite{KV1,KV2} (see \cite{KV3} for an
overview).

It is natural to expect that the main ideas of \cite{Kac} still work in this much more general
framework. Indeed, generalizing the classical group extension theory, G.I. Kac explained there
that, in order to construct an extension of finite groups $G_1$ and $G_2$, it is necessary and
sufficient: 1) To define a pair of compatible actions of $G_1$ and $G_2$ on each other (as on
sets) or, equivalently, $G_1$ and $G_2$ must be subgroups of a certain group $G$ such that
$G_1\cap G_2=\{e\}$ and any $g\in G$ can be written as $g=g_1g_2\ (g_1\in G_1, g_2\in G_2)$;
so $G_1$ and $G_2$ must form a matched pair (this term was introduced later by W.M.~Singer
\cite{Sin} and M.~Takeuchi \cite{Tak}). 2) To define a pair of compatible 2-cocycles for these
actions, so $G_1$ and $G_2$ must form a cocycle matched pair. Then an arbitrary extension has
the structure of their cocycle bicrossed product.

This last construction was studied intensively by S.~Majid both in algebraic and in analytic
aspects \cite{Maj1,Maj2,Maj3,Majbook,Maj-So}. In particular, in \cite{Maj1} he defines a
matched pair of Hopf algebras and studies their bicrossed product, later in \cite{Maj3,Maj-So}
and \cite{Majbook}, 6.3 he considers cocycle bicrossed products of Hopf algebras. In
\cite{Maj2} S.~Majid defines a matched pair of l.c.\ groups with continuous mutual actions and
constructs the corresponding bicrossed product Hopf-von Neumann algebra, which is a Kac
algebra under some additional assumption (modularity). Later T.~Yamanouchi proved \cite{Yam}
that in this situation one always gets a quasi Woronowicz algebra. Our analysis, involving
measurable and almost everywhere defined actions,  shows that in general this bicrossed
product is a l.c.\ quantum group. S.~Majid also gave very interesting concrete examples of
bicrossed products (without cocycles).

S.~Baaj and G.~Skandalis defined in \cite{B-S1} a matched pair of Kac systems and studied
their bicrossed product. In particular they considered matched pairs of locally compact groups
$G_1$ and $G_2$, and to keep themselves in the framework of regular multiplicative unitaries,
they assumed that $G_1$ and $G_2$ are closed subgroups of a locally compact group $G$ such
that the map $(g_1,g_2) \mapsto g_1g_2$ is a homeomorphism of $G_1\times G_2$ onto $G$. In
\cite{B-S2} they extended this notion requiring the above map to be a homeomorphism of
$G_1\times G_2$ only onto an open subset of $G$ with complement of measure zero. We will use
the same setting in Subsection~4.2 and Section~5. In \cite{B-S1} and \cite{Skand} S.~Baaj and
G.~Skandalis gave important concrete examples of bicrossed products.

We also want to mention the algebraic papers \cite{A,A-D,Bes-Dra1,Bes-Dra2,Hof,Mas1,Schn,Sin,Tak} on extensions of Hopf algebras and especially \cite{Mas2},
where A.~Masuoka established interesting connections with extensions of Lie bialgebras.
All this, as well as the recent paper \cite{Pal}, served as a motivation for our work aimed at
the construction of new examples of l.c.\ quantum groups. We explain the structure of our
work.

In Preliminaries we establish notations and briefly summarize the main definitions and results
of the theory of l.c.\ quantum groups in the von Neumann algebraic setting.

Section 1 is in a sense preparatory, although we believe its results are of independent
interest. Here we define and study cocycle actions of l.c.\ quantum groups on von Neumann
algebras generalizing twisted group actions \cite{Nak-Sut,Bus-Smi,Pac-Rae} and
ordinary (without cocycles) l.c.\ quantum group actions \cite{SV}. We construct the von
Neumann cocycle crossed product algebra which is an operator algebra version of the Hopf
algebraic cocycle crossed product \cite{Mon}, Chapter 7, we perform the dual weight construction
and obtain a canonical unitary implementation for these cocycle actions. Finally, we explain
relations between cocycle crossed products and cleft extensions of von Neumann algebras, which
is also motivated by the Hopf algebraic results \cite{Mon}, Chapter 7. Some facts concerning
operator spaces \cite{blecher,blecher2,blecher-paulsen,Efr-Rua} are needed for this
discussion.

Section 2 starts with the most general (and inevitably quite complicated) definition of a
cocycle matched pair of l.c.\ quantum groups. Then we construct the corresponding cocycle
bicrossed product (which is an operator algebra version of the Hopf algebraic construction
described in  \cite{Majbook}, 6.3), we show that it is a l.c.\ quantum group and compute its
invariant weights, fundamental unitary and the dual l.c.\ quantum group. Finally, we
investigate when the above cocycle bicrossed product is compact and discrete.

In Section 3 the notions of extensions and cleft extensions of l.c.\ quantum groups are
introduced and a one-to-one correspondence between cleft extensions and cocycle bicrossed
products is established. We also define isomorphic extensions and give necessary and
sufficient conditions for two cleft extensions to be isomorphic. These results are similar to
Hopf algebraic results obtained in \cite{A-D}.

Section~4 is devoted to the study of the special, but the most important, case when the
cocycle matched pair is formed by usual locally compact groups. When both groups are discrete,
we extend the result obtained by G.I.~Kac \cite{Kac} for finite groups and show that every
extension is cleft and so has a cocycle bicrossed product structure. Then for a matched pair
of locally compact groups we describe the corresponding bicrossed product and compute all the
ingredients of this locally compact quantum group. The same formulas were given by S.~Baaj and
G.~Skandalis \cite{B-S2} in the absence of cocycles. We also characterize the Kac algebra case
and here we extend \cite{Maj2}, Theorem 2.12. Finally, we introduce the group of extensions
related to a matched pair of locally compact groups and discuss relations with cocycles in the
sense of S.~Baaj and G.~Skandalis.

Discussing examples of matched pairs of locally compact groups and corresponding extensions in
Section~5, we start with a brief survey of known examples distinguishing the cases when all
extensions are or are not Kac algebras. Then we discuss the relation between locally compact
quantum groups obtained from matched pairs of Lie groups and the corresponding infinitesimal
objects: Hopf algebras and Lie bialgebras. Next, in Subsection~5.3, we discuss an example
which was already presented by S.~Baaj and G.~Skandalis in a talk in Oberwolfach
\cite{Skand}, but we include it here because it seems worthwhile to compute explicitly the
associated infinitesimal Hopf algebra and to show how this example is a deformation of the
$ax+b$-group. Finally, we construct a new example of a matched pair of Lie groups and compute
the corresponding extension getting this way a new example of a locally compact quantum group.
The description of this example was announced in \cite{SV2}. In Subsection~5.5 we present
the first, up to our knowledge, series of extensions with non-trivial cocycles which are not
Kac algebras. For this we compute explicitly a family of corresponding 2-cocycles.

{\bf Acknowledgements.} Both authors are grateful to Professor A.~Van Daele for stimulating discussions.
The first author wants to thank the research group in Leuven for the nice working atmosphere.
The second author is grateful to the Katholieke Universiteit Leuven and to the
Max-Planck-Institut f\"ur Mathematik in Bonn for kind hospitality during his work on this
paper and to Professor D.~Gurevich for helpful discussions. Both authors want to thank
Professors S.~Baaj and G.~Skandalis for their interesting comments.

\section*{Preliminaries}

We denote by $\ot$ the tensor product of Hilbert spaces (resp., von Neumann algebras) and by
$\Si$ (resp., $\si$) the flip map on it. We also use the leg-numbering notation. For example,
if $H, K$ and $L$ are Hilbert spaces and $X \in B(H \ot L)$, we denote by $X_{13}$ (resp.,
$X_{12},\ X_{23}$) the operator $(1 \ot \Si^*)(X \ot 1)(1 \ot \Si)$ (resp., $X\ot 1,\ 1\ot X$)
defined on $H \ot K \ot L$. If now $H = H_1 \ot H_2$ is itself a tensor product of two Hilbert
spaces, then we sometimes switch from the leg-numbering notation with respect to $H \ot K \ot
L$ to the one with respect to the finer tensor product $H_1 \ot H_2\ot K \ot L$, for example,
from $X_{13}$ to $X_{124}$. There is no confusion here, because the number of legs changes.

We denote the \strong closure of a subset $A$ of a von Neumann algebra $N$ by $A^\sluit$ and
we use \cite{Stra} as a general reference to the modular theory of {\it normal semifinite
faithful} (n.s.f.) weights on von Neumann algebras. If $\te$ is a weight on a von Neumann
algebra $N$, we use the notations $$ \cM_\te^+ = \{ x \in N^+ \mid \te(x) < + \infty \},
\qquad \cN_\te = \{ x \in N \mid x^*x \in \cM_\te^+ \} \quad\text{and}\quad \cM_\te =
\operatorname{span} \cM_\te^+ \; . $$ If $N_0$ is a von Neumann subalgebra of $N$ and if $T$
is an operator valued weight from $N$ to $N_0$ \cite{Stra}, Section~11.5, we denote by $\cN_T$
the set of elements $x \in N$ such that $T(x^*x)$ is a bounded operator. Recall that if $T$ is
a \nsf operator valued weight and $\te_0$ is a \nsf weight on $N_0$, then $\te:=\te_0 \, T$
defines a \nsf weight on $N$. The modular automorphism group of $\te_0$ is the restriction of
the modular automorphism group of $\te$ to the von Neumann algebra $N_0$, \cite{Stra},
Section~11.9.

If $\te$ is a \nsf weight on a von Neumann algebra $N$ and if $(H_\te,\pi_\te,\la_\te)$ is a
GNS-construction for $\te$, we recall that the GNS-map $\la_\te$ is \strong -- norm closed.
This expression means that the map $\la_\te$ is closed for the \strong topology on $N$ and the
norm topology on $H_\te$. We denote by $\cT_\te$ the Tomita $^*$-algebra, consisting of
elements $x \in N$ analytic w.r.t.\ the modular group $\si^\te$ of $\te$ and such that
$\si^\te_z(x) \in \cN_\te \cap \cN_\te^*$ for all $z \in \C$.

In this paper we use systematically l.c.\ quantum groups in von Neumann algebraic setting
\cite{KV2} including special types of morphisms in the definition of extensions in Section 3.
However it should be mentioned that a comprehensive definition of morphisms and then of a
category of l.c.\ quantum groups can be done only in the universal $C^*$-algebraic setting
\cite{Kus}; this was known already for Kac algebras \cite{E-S}, Chapter 5.

A pair $(M,\de)$ is called a (von Neumann algebraic) l.c.\ quantum group \cite{KV2} when
\begin{itemize}
\item $M$ is a von Neumann algebra and $\de : M \recht M \ot M$ is
a normal and unital $*$-homomorphism satisfying the coassociativity relation : $(\de \ot \io)\de = (\io \ot
\de)\de$.
\item There exist \nsf weights $\vfi$ and $\psi$ on $M$ such that
\begin{itemize}
\item $\vfi$ is left invariant in the sense that $\vfi \bigl( (\om \ot
\io)\de(x) \bigr) = \vfi(x) \om(1)$ for all $x \in \Mfi^+$ and $\om \in M_*^+$,
\item $\psi$ is right invariant in the sense that $\psi \bigl( (\io \ot
\om)\de(x) \bigr) = \psi(x) \om(1)$ for all $x \in \Mpsi^+$ and $\om \in M_*^+$.
\end{itemize}
\end{itemize}
From \cite{KV1}, Theorem~7.14 we know that left invariant weights on $(M,\de)$ are unique up
to a positive scalar and the same holds for right invariant weights.

Let us fix a left invariant \nsf weight $\vfi$ on $(M,\de)$ and represent $M$ on the GNS-space
of $\vfi$ such that $(H,\io,\la)$ is a GNS-construction for $\vfi$. Then we can define a
unitary $W$ on $H \ot H$ by
$$W^* (\la(a) \ot \la(b)) = (\la \ot \la)(\de(b)(a \ot 1)) \quad\text{for all}\; a,b \in \Nfi
\; .$$
Here $\la \ot \la$ denotes the canonical GNS-map for the tensor product weight $\vfi \ot
\vfi$.
One proves that $W$ satisfies the pentagonal equation: $W_{12} W_{13} W_{23} = W_{23} W_{12}$.
We say that $W$ is a multiplicative unitary.
The comultiplication can be given in terms of $W$ by the formula $\de(x) = W^* (1 \ot x) W$
for all $x \in M$. Also the von Neumann algebra $M$ can be written in terms of $W$ as
$$M = \{ (\io \ot \om)(W) \mid \om \in B(H)_* \}^\sluit \; .$$

Next the l.c.\ quantum group $(M,\de)$ has an antipode $S$, which is the unique \strong closed
linear map from $M$ to $M$ satisfying $(\io \ot \om)(W) \in \cD(S)$ for all $\om \in B(H)_*$,
$S(\io \ot \om)(W) = (\io \ot \om)(W^*)$ and such that the elements $(\io \ot \om)(W)$ form a
\strong core for $S$. $S$ has a polar decomposition $S = R \tau_{-i/2}$ where $R$ is an
anti-automorphism of $M$ and $(\tau_t)$ is a strongly continuous one-parameter group of
automorphisms of $M$. We call $R$ the unitary antipode and $(\tau_t)$ the scaling group of
$(M,\de)$. From \cite{KV1}, Proposition~5.26 we know that $\si (R \ot R) \de = \de R$. So
$\vfi R$ is a right invariant weight on $(M,\de)$ and we take $\psi:= \vfi R$. Let us denote
by $(\si_t)$ the modular automorphism group of $\vfi$. From \cite{KV1}, Proposition~6.8 we get
the existence of a number $\nu > 0$, called the scaling constant, such that $\psi \, \si_t =
\nu^{-t} \, \psi$ for all $t \in \R$. Hence we get the existence of a unique positive,
self-adjoint operator $\sde_M$ affiliated to $M$, such that $\si_t(\sde_M) = \nu^t \, \sde_M$
for all $t \in \R$ and $\psi = \vfi_{\sde_M}$, see \cite{KV1}, Definition~7.1. Formally this
means that $\psi(x) = \vfi(\sde_M^{1/2} x \sde_M^{1/2})$, and for a precise definition of
$\vfi_{\sde_M}$ we refer to \cite{SV3}. The operator $\sde_M$ is called the modular element of
$(M,\de)$. If $\sde_M=1$ we call $(M,\de)$ unimodular. Let us also mention that the scaling
constant can be characterized as well by the relative invariance $\vfi \, \tau_t = \nu^{-t} \,
\vfi$.

We use the notation $\deop$ to denote the opposite comultiplication defined by $\deop:= \si \de$.

From \cite{KV1}, Notation~7.2 we get the canonical GNS-construction $(H,\io,\Gamma)$ for the
\nsf weight $\psi$, formally given by $\Gamma(x) = \la(x \sde_M^{1/2})$. Then we can define another multiplicative unitary $V$ by
$$V (\Ga(a) \ot \Ga(b)) = (\Ga \ot \Ga)(\de(a)(1 \ot b)) \quad\text{for all}\; a,b \in \Npsi
\; .$$ The comultiplication can be expressed by $\de(x) = V(x \ot 1) V^*$ for all $x \in M$.

The dual l.c.\ quantum group $(\Mh,\deh)$ is defined in \cite{KV1}, Section~8. Its von Neumann
algebra $\Mh$ is $$\Mh = \{(\om \ot \io)(W) \mid \om \in B(H)_* \}^\sluit$$ and the
comultiplication is given by $\deh(x) = \Si W (x \ot 1) W^* \Si$ for all $x \in \Mh$. If we
turn the predual $M_*$ into a Banach algebra with product $\om \, \mu = (\om \ot \mu)\de$ and
define $$\lambda: M_* \recht \Mh : \lambda(\om) = (\om \ot \io)(W),$$ then $\lambda$ is a
homomorphism and $\lambda(M_*)$ is a \strong dense subalgebra of $\Mh$. To construct
explicitly a left invariant \nsf weight $\vfih$ with GNS-construction $(H,\io,\lah)$, first
introduce the space $$\cI = \{ \om \in M_* \mid \; \text{there exists a vector} \; \xi(\om)
\in H \tekst{such that} \om(x^*) = \langle \xi(\om), \la(x) \rangle \tekst{for all} x \in \Nfi
\} \; .$$ If $\om \in \cI$ then such a vector $\xi(\om)$ clearly is uniquely determined. Next
one proves that there exists a unique \nsf weight $\vfih$ on $\Mh$ with GNS-construction
$(H,\io,\lah)$ such that $\lambda(\cI)$ is a \strong -- norm core for $\lah$ and
$$\lah(\lambda(\om)) = \xi(\om) \tekst{for all} \om \in \cI \; .$$ One proves that the weight
$\vfih$ is left invariant, and the associated multiplicative unitary is denoted by $\hat{W}$.
From \cite{KV1}, Proposition~8.16 it follows that $\hat{W} = \Si W^* \Si$.

Since $(\Mh,\deh)$ is again a l.c.\ quantum group, we can introduce the antipode $\hat{S}$,
the unitary antipode $\hat{R}$ and the scaling group $(\hat{\tau}_t)$ exactly as we did it for
$(M,\de)$. Also, we can again construct the dual of $(\Mh,\deh)$, starting from the left
invariant weight $\vfih$ with GNS-construction $(H,\io,\lah)$. From \cite{KV1}, Theorem~8.29
we get that the bidual l.c.\ quantum group $(\Mhh,\dehh)$ is  isomorphic to $(M,\de)$. Also,
defining $\hat{\cI}$ and $\hat{\xi}$ similarly to the definition of $\cI$ and $\xi$, it
follows from \cite{KV1}, Proposition~8.30 that, for all $\om \in \hat{\cI}$, $$\la((\io \ot
\om)(W^*)) = \hat{\xi}(\om) \; .$$

We denote by $(\sih_t)$ the modular automorphism groups of the weight $\vfih$. The modular
conjugations of the weights $\vfi$ and $\vfih$ will be denoted by $J$ and $\Jh$ respectively.
Then it is worthwhile to mention that $$R(x) = \Jh x^* \Jh \quad\text{for all} \; x \in M
\qquad\text{and}\qquad \Rh(y) = J y^* J \quad\text{for all}\; y \in \Mh \; .$$ From
\cite{KV2}, Proposition~2.15 we know that $$V = (\Jh \ot \Jh) \Si W^* \Si (\Jh \ot \Jh)$$ and
in particular $V \in {\Mh}^\prime \ot M$.

Let us mention important special cases of l.c.\ quantum groups.

a) If $({\tau}_t)$ is trivial and the modular element $\delta_M$ is affiliated to the center
of $M$, $(M,\de)$ becomes a {\it Kac algebra} \cite{E-S}. Let us explain this in more detail.
From \cite{E-S} we know that $(M,\de)$ is a Kac algebra if and only if $(\tau_t)$ is trivial
and $\si_t \, R = R \, \si_{-t}$ for all $t \in \R$. Now denote by $(\si'_t)$ the modular
automorphism group of $\psi$. Because $\psi = \vfi R$ we get that $\si'_t \, R = R \,
\si_{-t}$ for all $t \in \R$. Hence $(M,\de)$ is a Kac algebra if and only if $(\tau_t)$ is
trivial and $\si'=\si$. From \cite{SV3} we know that $\si'_t(x) = \sde_M^{it} \si_t(x)
\sde_M^{-it}$ for all $x \in M$ and $t \in \R$. Hence $\si'=\si$ if and only if $\sde_M$ is
affiliated to the center of $M$.

In particular, $(M,\de)$ is a Kac algebra if $M$ is commutative. Then $(M,\de)$ is generated
by a usual l.c.\ group $G:\ M=L^{\infty}(G),\ (\de f)(g,h) = f(gh),\ (Sf)(g) = f(g^{-1}),\
\vfi(f)=\int f(g)\; dg$, where $f\in L^{\infty}(G),\ g,h\in G$ and we integrate with respect
to the left Haar measure $dg$ on $G$. The right invariant weight $\psi$ is given by $\psi(f) =
\int f(g^{-1}) \; dg$. The modular element $\sde_M$ is given by the strictly positive function
$g \mapsto \sde_G(g)^{-1}$, where $\sde_G$ is the modular function of the l.c.\ group $G$.

The von Neumann algebra $M$ acts on $H=L^2(G)$ by multiplication and
$$(W\xi)(g,h)=\xi(g,g^{-1}h)$$ for all $\xi\in H\ot H=L^2(G\times G)$. Then $\Mh=\cL(G)$ is
the group von Neumann algebra generated by the operators $(\lambda_g)_{g\in G}$ of the left
regular representation of $G$ and $\deh(\lambda_g)=\lambda_g\ot\lambda_g$. Clearly,
$\dehop:=\si\deh=\deh$; so, $\deh$ is cocommutative.

b) A l.c.\ quantum group is called compact if its Haar measure is finite: $\vfi(1)<+\infty$,
which is equivalent to the fact that the norm closure of $\{(\io\ot\om)(W)\vert\om\in
B(H)_*\}$ is a unital $C^*$-algebra. A l.c.\ quantum group $(M,\de)$ is called discrete if
$(\Mh,\deh)$ is compact.

Sometimes we refer to Hopf algebra theory, and then we use the Sweedler notation for $\de$:
$$\de(a) = \sum a_{(1)} \ot a_{(2)} \; .$$

\section{Cocycle crossed products and the dual weight construction}
\subsection{Cocycle actions, crossed products and the dual action}
\begin{definition} \label{11}
We call a pair $(\al,\cU)$ a cocycle action of a l.c.\ quantum group $(M,\de)$ on a von
Neumann algebra $N$ if $$\al:N \recht M \ot N$$ is a normal, injective and unital
$*$-homomorphism, $$\cU \in M \ot M \ot N$$ is a unitary, and if $\al$ and $\cU$ satisfy
\begin{align*}
(\io \ot \al)\al(x) &= \cU \, (\de \ot \io)\al(x) \, \cU^* \tekst{for all} x \in N \\
(\io \ot \io \ot \al)(\cU) \; (\de \ot \io \ot \io)(\cU) &= (1 \ot \cU) \; (\io \ot \de \ot \io)(\cU)
\; .
\end{align*}
\end{definition}
Consider some special cases. First, let $(\al,u)$ be a twisted action of a (separable) l.c.\
group $G$ on a ($\si$-finite) von Neumann algebra $N$. This means that $G \recht
\operatorname{Aut} N : s \mapsto \al_s$ and $u:G \times G \recht N$ are Borel maps, such that
$u$ takes values in the set of unitaries of $N$ and
\begin{align*}
\al_s \; \al_t &= \operatorname{Ad} u(s,t) \; \al_{st} \\
\al_r(u(s,t)) \; u(r,st) &= u(r,s) \; u(rs,t)
\end{align*}
nearly everywhere. Putting $M=L^\infty(G)$, one can identify $M \ot N$ with $L^\infty(G,N)$
and hence one can define $\al: N \recht L^\infty(G) \ot N$ by $(\al(x))(s) = \al_{s^{-1}}(x)$
for $x \in N$ and $s \in G$. We also identify $M \ot M \ot N$ with $L^\infty(G \times G,N)$
and define $\cU \in L^\infty(G) \ot L^\infty(G) \ot N$ by $\cU(s,t)=u(t^{-1},s^{-1})$. Then
$(\al,\cU)$ is a cocycle action of the commutative l.c.\ quantum group $(L^\infty(G),\de_G)$
on $N$.

Second, if $\cU=1$, we obtain the definition of an
ordinary action $\al$ of $(M,\de)$ on $N$ \cite{SV}, Definition 1.1.

Third, let us explain the link with the Hopf algebra theory \cite{Mon}, Definition 4.1.1 and
Lemma 7.1.2. Let $(H,\de,S,\eps)$ be a Hopf algebra and let $A$ be a unital algebra equipped
with a linear map $H \ot A \recht A:h \ot a\recht  h \cdot a$ such that $$h \cdot (ab) = \sum
(h_{(1)} \cdot a) \; (h_{(2)} \cdot b) \tekst{and} h \cdot 1 = \eps(h)1 \; $$ and equipped
with a convolution invertible linear map $\si:H \ot H\recht A$ such that
\begin{align*}
h \cdot (k \cdot a) &= \sum \si(h_{(1)},k_{(1)}) \; \bigl( (h_{(2)}k_{(2)}) \cdot a \bigr) \;
\si^{-1}(h_{(3)},k_{(3)}) \\
\sum \bigl( h_{(1)} \cdot \si(k_{(1)},m_{(1)}) \bigr) \; \si(h_{(2)},k_{(2)} m_{(2)}) &= \sum
\si(h_{(1)},k_{(1)}) \; \si(h_{(2)} k_{(2)}, m) \; .
\end{align*}
Here the convolution invertibility of $\si$ means that there exists a linear map $\si^{-1}$
from $H \ot H$ to $A$ such that $$\sum \si(h_{(1)},k_{(1)}) \; \si^{-1}(h_{(2)},k_{(2)}) =
\eps(h)\eps(k)1 = \sum \si^{-1}(h_{(1)},k_{(1)}) \; \si(h_{(2)},k_{(2)}) \; . $$ If now $H$ is
finite dimensional, the linear dual $\hat{H}$ is a Hopf algebra by $(\om \mu)(h) = \sum
\om(h_{(1)}) \mu(h_{(2)})$ and $\sum \om_{(1)}(h) \om_{(2)}(k) = \om(kh)$ for all $\om,\mu \in
\hat{H}$ and $h,k \in H$. Observe that we use the opposite comultiplication on $\hat{H}$.
Identifying $\hat{H} \ot A$ with the space of linear maps from $H$ to $A$ one defines $\al : A
\recht \hat{H} \ot A$ by $\al(a)(h) = h \cdot a$. Then $\al$ is a unital homomorphism. Further
identifying $\hat{H} \ot \hat{H} \ot A$ with the space of linear maps from $H \ot H$ to $A$ we
define $\cU \in \hat{H} \ot \hat{H} \ot A$ by $\cU(h,k)=\si(k,h)$. Then $\cU$ is invertible
and
\begin{align*}
(\io \ot \al)\al(a) &= \cU \, (\de \ot \io)\al(a) \, \cU^{-1} \\
(\io \ot \io \ot \al)(\cU) \; (\de \ot \io \ot \io)(\cU) &= (1 \ot \cU) \; (\io \ot \de \ot \io)(\cU)
\; .
\end{align*}
In this setting without involution it is natural to replace $\cU^*$ by $\cU^{-1}$, so the Hopf algebraic definition of a cocycle action agrees with our definition.

Let us now start with the operator algebraic theory of cocycle actions and their crossed products.
\begin{notation} \label{12}
If $(\al,\cU)$ is a cocycle action of $(M,\de)$ on a von Neumann algebra $N$ we introduce the
notation $$\Wtil = (W \ot 1)\cU^* \; ,$$ and then $\Wtil$ is a unitary in $M \ot B(H) \ot N$.
\end{notation}
\begin{definition} \label{13}
Given a cocycle action $(\al,\cU)$ of a l.c.\ quantum group $(M,\de)$ on a von Neumann algebra
$N$, the cocycle crossed product $\cros$ is the von Neumann subalgebra of $B(H) \ot N$
generated by $$\al(N) \quad\text{and}\quad \{ (\om \ot \io \ot \io)(\Wtil) \mid \om \in M_* \}
\; .$$
\end{definition}
As in the case of ordinary actions \cite{SV}, Proposition 2.2 we can define a dual action $\alh$ of
$(\Mh,\dehop)$ on $\cros$ with trivial cocycle.
\begin{proposition} \label{14}
There exists a unique action $\alh$ of $(\Mh,\dehop)$ on $\cros$ such that
\begin{align*}
\alh(\al(x)) & = 1 \ot \al(x) \tekst{for all} x \in N \\
(\io \ot \alh)(\Wtil) &= W_{12} \Wtil_{134} \; .
\end{align*}
Moreover denoting by $\Vtil$ the unitary $(J \ot J)\Si W \Si (J \ot J)$ we have
$$\alh(z) = (\Vtil \ot 1)(1 \ot z)(\Vtil^* \ot 1) \tekst{for all} z \in \cros.$$
\end{proposition}
\begin{proof}
Observe that $\Vtil \in \Mh \ot M'$ and $\Vtil (1 \ot x) \Vtil^* = \dehop(x)$ for all $x
\in \Mh$. So one gets for every $x \in N$
$$(\Vtil \ot 1) (1 \ot \al(x)) (\Vtil^* \ot 1) = 1 \ot \al(x) \; .$$
Next we have
\begin{align*}
\Vtil_{23} \Wtil_{134} \Vtil_{23}^* &= \Vtil_{23} W_{13} \cU^*_{134} \Vtil_{23}^*
=\Vtil_{23} W_{13} \Vtil_{23}^*  \; \cU^*_{134} \\
&=(\io \ot \dehop)(W)_{123} \; \cU^*_{134}
=W_{12} W_{13} \cU^*_{134} = W_{12} \Wtil_{134} \; .
\end{align*}
So it is clear that we can define a normal, unital and injective $*$-homomorphism
$$\alh: \cros \recht \Mh \ot (\cros)$$
by $\alh(z) = (\Vtil \ot 1)(1 \ot z) (\Vtil^* \ot 1)$ for all $z \in \cros$, and then $\alh$
satisfies
$\alh(\al(x)) = 1 \ot \al(x)$ for $x \in N$ and $(\io \ot \alh)(\Wtil) = W_{12} \Wtil_{134}$.
It is obvious that
$$(\io \ot \alh)\alh(\al(x)) = (\dehop \ot \io)\alh(\al(x))$$
for all $x \in N$. Further we have
\begin{align*}
(\io \ot \io \ot \alh)(\io \ot \alh)(\Wtil) &= (\io \ot \io \ot \alh)(W_{12} \Wtil_{134}) \\
&= W_{12} W_{13} \Wtil_{145} = (\io \ot \dehop)(W)_{123} \Wtil_{145} \\
&= (\io \ot \dehop \ot \io)(W_{12} \Wtil_{134}) = (\io \ot \dehop \ot \io)(\io \ot \alh)(\Wtil)
\; .
\end{align*}
Both statements together give $(\io \ot \alh)\alh(z) = (\dehop \ot \io)\alh(z)$ for all $z \in \cros$.
Hence $\alh$ is indeed an action of $(\Mh,\dehop)$ on the von Neumann algebra
$\cros$. The uniqueness statement is obvious.
\end{proof}
Let us get another formula for $\alh$.
Since $\cros \subset B(H) \ot N$, the following proposition makes sense.
\begin{proposition} \label{15}
For all $z \in \cros$ we have
$$\alh(z) = \Wtil (\io \ot \al)(z) \Wtil^* \; . $$
\end{proposition}
\begin{proof}
By definition $(\io \ot \al)\al(x) = \Wtil^*(1 \ot \al(x)) \Wtil$, so
$$\Wtil (\io \ot \al)\al(x) \Wtil^* = 1 \ot \al(x) = \alh(\al(x)) \; .$$
Next we have
\begin{align*}
\Wtil_{234} (\io \ot \io \ot \al)(\Wtil) \Wtil^*_{234} &= \Wtil_{234} W_{12} (\io \ot \io \ot
\al)(\cU^*) \Wtil^*_{234} \\
&= \Wtil_{234} W_{12} \; (\de \ot \io \ot \io)(\cU) (\io \ot \de \ot \io)(\cU^*) (1 \ot \cU^*) \;
\Wtil^*_{234} \\
&= W_{23} W_{12} W_{23}^* \cU^*_{134} = W_{12} W_{13} \cU^*_{134} \\
&= W_{12} \Wtil_{134} = (\io \ot \alh)(\Wtil) \; .
\end{align*}
Both computations together give the formula stated in the proposition.
\end{proof}
The following result will be needed in Section~5.
\begin{proposition} \label{center}
Let $x \in N$. Then $\al(x)$ commutes with the elements $\{(\om \ot \io \ot \io)(\Wtil) \mid
\om \in B(H)_* \}$ if and only if $\al(x) = 1 \ot x$.
\end{proposition}
\begin{proof}
Since $(\io \ot \al)\al(x) = \Wtil^* (1 \ot \al(x)) \Wtil$,  $\al(x)$
commutes with the elements $\{(\om \ot \io \ot \io)(\Wtil) \mid \om \in B(H)_* \}$ if and only
if $(\io \ot \al)\al(x) = 1 \ot \al(x)$. Since $\io \ot \al$ is faithful this equality is valid if and only if $\al(x) = 1 \ot x$.
\end{proof}
\subsection{Stabilization of cocycle actions}
\begin{definition} \label{16}
A cocycle action $(\al,\cU)$ of $(M,\de)$ on a von Neumann algebra $N$ is said to be stabilizable
with a unitary $X\in M \ot N$ if
$$(1 \ot X)(\io \ot \al)(X) = (\de \ot \io)(X) \cU^* \; .$$
\end{definition}
\begin{proposition} \label{17}
Let $(\al,\cU)$ be a cocycle action of $(M,\de)$ on $N$ which is stabilizable with a unitary $X\in M \ot N$. Then the formulas
$$\be:N \recht M \ot N : \be(x) = X \al(x) X^*\ and\ \Phi: z \mapsto X^* z X$$
define respectively an action of $(M,\de)$ on $N$ and a $*$-isomorphism from $M \kruisje{\be} N$ onto $\cros$
satisfying
$$\alh \Phi = (\io \ot \Phi) \hat{\be} \; .$$
\end{proposition}
\begin{proof}
Define $\be$ as above. Then one has for all $x \in N$
\begin{align*}
(\io \ot \be) \be(x) &= (1 \ot X) (\io \ot \al)(X) \; (\io \ot \al)\al(x) \; (\io \ot \al)(X^*) (1
\ot X^*) \\
&= (\de \ot \io)(X) \cU^* \; \cU (\de \ot \io)\al(x) \cU^* \; \cU (\de \ot \io)(X^*)
= (\de \ot \io)\be(x) \; .
\end{align*}
Hence $\be$ is an action of $(M,\de)$ on $N$. Considering $\Phi$ as an isomorphism of $B(H) \ot N$, we clearly have $\Phi(\be(x)) = \al(x)$ for all $x \in N$ and
\begin{align*}
(\io \ot \Phi)(W \ot 1) &= (1 \ot X^*)(W \ot 1)(1 \ot X) = (W \ot 1) (\de \ot \io)(X^*) (1 \ot
X) \\
&=(W \ot 1) \cU^* (\io \ot \al)(X^*) = \Wtil (\io \ot \al)(X^*) \; .
\end{align*}
So $\Phi(M \kruisje{\be} N) \subset \cros$; similarly one
proves the
converse inclusion. Hence $\Phi$ is an isomorphism of $M \kruisje{\be} N$ onto $\cros$.

Recalling the definition of the dual actions $\alh$ and $\hat{\be}$ and the unitary $\Vtil$
implementing them, we only have to observe that $1 \ot X$ and $\Vtil \ot 1$ commute to get the
formula $(\io \ot \Phi) \hat{\be} = \alh \Phi$.
\end{proof}
The next proposition shows that many cocycle actions
are stabilizable.
\begin{proposition} \label{18}
Let $(\al,\cU)$ be a cocycle action of $(M,\de)$ on $N$. Then $(\al \ot \io, \cU \ot 1)$ is a
cocycle action of $(M,\de)$ on $N \ot B(H)$ which is stabilizable with the unitary
$$X = V^*_{31} \cU^*_{312} \; .$$
\end{proposition}
\begin{proof}
It is easy to check that $X \in M \ot N \ot B(H)$. Next we have to prove that
$$(1 \ot X) (\io \ot \al \ot \io)(X) = (\de \ot \io \ot \io)(X) (\cU^* \ot 1) \; .$$
So we have to prove that
$$X_{341} (\io \ot \io \ot \al)(X_{231}) = (\io \ot \de \ot \io)(X_{231}) \cU^*_{234} \; .$$
Starting the computation on the left hand side we get
\begin{align*}
X_{341} (\io \ot \io \ot \al)(X_{231}) &= V^*_{13} \cU^*_{134} \; V^*_{12} (\io \ot \io \ot
\al)(\cU^*) \\
&= V^*_{13} \cU^*_{134} \; V^*_{12} \; (\de \ot \io \ot \io)(\cU) (\io \ot \de \ot \io)(\cU^*) \cU^*_{234}
\\ &= V^*_{13} \cU^*_{134} \; \cU_{134} V^*_{12} \; (\io \ot \de \ot \io)(\cU^*) \cU^*_{234} \\
&= (\io \ot \de)(V^*)_{123} (\io \ot \de \ot \io)(\cU^*) \cU^*_{234} \\
&= (\io \ot \de \ot \io)(X_{231}) \cU^*_{234} \; .
\end{align*}
This computation proves our result.
\end{proof}
Let us give now a definition of the fixed point algebra of a cocycle action.
\begin{definition} \label{19}
Let $(\al,\cU)$ be a cocycle action of $(M,\de)$ on $N$. Then $$N^\al = \{ x \in N \mid \al(x) = 1 \ot x \}$$
is called the fixed point algebra of $(\al,\cU)$. It is clear that $N^\al$ is a von
Neumann subalgebra of $N$.
\end{definition}
\begin{theorem} \label{110}
Let $(\al,\cU)$ be a cocycle action of a l.c.\ quantum group $(M,\de)$ on a von Neumann
algebra $N$ and $\alh$ the dual action of $(\Mh,\dehop)$ on the cocycle crossed product
$\cros$. Then
\begin{enumerate}
\item $(\cros)^{\alh} = \al(N) \; ,$
\item $(\cros \cup M' \ot \C)\dpr = B(H) \ot N \; ,$
\item $\cros = \bigl(\lspan\{ (\om \ot \io \ot \io)(\Wtil) \al(x) \mid x \in N, \om \in M_* \}\bigr)^{- \strong} \; .$
\end{enumerate}
\end{theorem}
\begin{proof}
First let $(\al,\cU)$ be stabilizable with a unitary $X$. Using the action
$\be$ and the isomorphism $\Phi$ of Proposition~\ref{17} we get the first two statements,
because they are valid for ordinary actions \cite{SV}, Theorems 2.6 and 2.7. Recalling that $(\io \ot \Phi)(W
\ot 1) = \Wtil (\io \ot \al)(X^*)$ and using that
$$M \kruisje{\be} N = \bigl(\lspan\{ \bigl( (\om \ot \io)(W) \ot 1 \bigr) \; \be(x) \mid x \in N, \om \in M_*
\}\bigr)^\sluit$$
we get
$$\cros = \bigl(\lspan\{ (\om \ot \io \ot \io) \bigl( \Wtil (\io \ot \al)(X^*) \bigr) \; \al(x) \mid x \in N, \om \in M_*
\}\bigr)^{- \strong} \; .$$
Now we define the following linear subspace of $\cros \; :$
$$\cO := \bigl(\lspan\{ (\om \ot \io \ot \io)(\Wtil) \al(x) \mid x \in N, \om \in M_*
\}\bigr)^{- \strong} \; .$$
Choose $\xi,\eta \in H$, $x \in N$ and an orthonormal basis $(e_i)_{i \in I}$ for $H$. Then we
have with \strong convergence
$$(\om_{\xi,\eta} \ot \io \ot \io)(\Wtil (\io \ot \al)(X^*)) \al(x) =
\sum_{i \in I} (\om_{e_i,\eta} \ot \io \ot \io)(\Wtil) \al \bigl( (\om_{\xi,e_i} \ot \io)(X^*)x
\bigr) \in \cO \; .$$
So $\cros \subset \cO$ and hence $\cO = \cros$.
This proves the result for stabilizable actions. From Proposition~\ref{18} the general result now follows.
\end{proof}
\subsection{The dual weight construction}
Theorem~\ref{110} shows
that $\al(N)$ is the fixed point algebra of the dual action $\alh$ of $(\Mh,\dehop)$ on $\cros$, so the
formula $T:=(\vfih \ot \io \ot \io) \alh$ defines a normal and faithful operator valued weight from $\cros$ to $\al(N)$ \cite{SV}, Proposition 1.3.
To define then the dual weight $\tetil$ for a given weight $\te$ on $N$ by the formula $\tetil = \te \al^{-1} T$, we have to prove first that $T$ is semifinite.
\begin{lemma} \label{111}
We still write $T:=(\vfih \ot \io \ot \io) \alh$. For every $\om \in \cI$ and $x \in N$
we have
$$(\om \ot \io \ot \io)(\Wtil) \al(x) \in \cN_T$$
and
$$T \bigl( \al(x^*) (\om \ot \io \ot \io)(\Wtil)^* (\om \ot \io \ot \io)(\Wtil) \al(x) \bigr)
= \|\xi(\om) \|^2 \al(x^*x) \; .$$
\end{lemma}
\begin{proof}
Let $N$ act on $K$.
Choose $\om \in \cI$, $x \in N$,
a vector $\xi \in H \ot K$ and an orthonormal basis $(e_i)_{i \in I}$ for $H \ot K$. Define
$$z:=\al(x^*) (\om \ot \io \ot \io)(\Wtil)^* (\om \ot \io \ot \io)(\Wtil) \al(x) \; .$$
Then we get
\begin{align*}
(\io \ot \om_{\xi,\xi}) \alh(z)
&= (\io \ot \om_{\al(x)\xi,\al(x)\xi}) \bigl( (\om \ot \io \ot \io \ot \io)(W_{12}
\Wtil_{134})^* (\om \ot \io \ot \io \ot \io)(W_{12}
\Wtil_{134}) \bigr) \\
&= \sum_{i \in I} \bigl( (\io \ot \om_{\al(x)\xi,e_i})(\Wtil) \om \ot \io \bigr) (W)^*
\bigl( (\io \ot \om_{\al(x)\xi,e_i})(\Wtil) \om \ot \io \bigr) (W) \; .
\end{align*}
Because $\vfih$ is normal we get that
\begin{align*}
\vfih (\io \ot \om_{\xi,\xi})\alh(z) &= \sum_{i \in I} \bigl\| \xi \bigl((\io \ot
\om_{\al(x)\xi,e_i})(\Wtil) \om \bigr) \bigr\|^2 \\
&= \sum_{i \in I} \| (\io \ot \om_{\al(x)\xi,e_i})(\Wtil) \xi(\om) \|^2 \\
&= \| \xi(\om) \ot \al(x) \xi \|^2 = \|\xi(\om)\|^2 \om_{\xi,\xi}(\al(x^*x)) \; .
\end{align*}
This final computation proves the lemma.
\end{proof}
Combining this lemma with Theorem~\ref{110} one gets that the operator valued weight $T$ is semifinite. So one can use the theory of operator valued weights \cite{Stra} to give the following definition.
\begin{definition} \label{112}
Let $(\al,\cU)$ be a cocycle action of $(M,\de)$ on $N$. Given a \nsf weight $\te$ on $N$ we can
define the dual \nsf weight $\tetil$ on $\cros$ by the formula
$$\tetil = \te \al^{-1} \; (\vfih \ot \io \ot \io) \alh \; .$$
\end{definition}
Let us fix a \nsf weight $\te$ on $N$ with GNS-construction $(K,\pite,\late)$. Let $\tetil$ be
the dual weight of $\te$.
Lemma~\ref{111} implies by polarization the following formula for $\tetil$ on a dense subset of $\cros$.
\begin{corollary} \label{112bis}
For all $\om,\mu \in \cI$ and $x,y \in \Nte$ the element $(\om \ot \io \ot \io)(\Wtil) \al(x)$
belongs to $\cN_{\tetil}$ and
$$\tetil \bigl( \al(y^*) (\mu \ot \io \ot \io)(\Wtil)^* (\om \ot \io \ot \io)(\Wtil) \al(x)
\bigr) = \langle \xi(\om) \ot \late(x), \xi(\mu) \ot \late(y) \rangle.$$
\end{corollary}
In order to apply the dual weight $\tetil$ further in this chapter, we need a concrete
GNS-construction for it. The corollary above suggests to define a GNS-map $\latil$ for
$\tetil$ by the formula $$\latil \bigl( (\om \ot \io \ot \io)(\Wtil) \al(x) \bigr) = \xi(\om)
\ot \late(x) \; .$$ To do this, it is not enough to know that the elements $(\om \ot \io \ot
\io)(\Wtil) \al(x)$ span a dense set of $\cros$, we also need to prove that they form a core
for a GNS-map for $\tetil$.
\begin{proposition} \label{113}
Let $(\al,\cU)$ be a cocycle action of $(M,\de)$ on $N$ and $\te$ a \nsf weight on $N$
with GNS-construction $(K,\pite,\late)$. Then there exists a unique GNS-construction $(H \ot
K,\io \ot \pite,\latil)$ for the dual weight $\tetil$ satisfying
\begin{enumerate}
\item $\lspan\{ (\om_{\eta,\la(b)} \ot \io \ot \io)(\Wtil) \al(x) \mid \eta \in H,b \in \cT_\vfi, x \in \Nte \}$ is a
\strong -- norm core for $\latil$.
\item $\latil \bigl( (\om \ot \io \ot \io)(\Wtil) \al(x) \bigr) = \xi(\om) \ot \late(x)$ for
all $\om \in \cI$ and $x \in \Nte$.
\end{enumerate}
\end{proposition}
\begin{proof}
Let $N$ act on $K$ by $\pite$, hence we may assume that $\pite=\io$.
Suppose first that $(\al,\cU)$ is stabilizable with the unitary $X$. Then we can construct the
action $\be$ and the isomorphism $\Phi$ as in Proposition~\ref{17}. Because $(\io \ot \Phi)\hat{\be}
= \alh \Phi$ we get immediately that $\tetil \Phi= \tetil_\be$, where $\tetil_\be$ is the dual
weight of $\te$ on $M \kruisje{\be} N$. Let $(H \ot K, \io, \latil_\be)$ be the
canonical GNS-construction for $\tetil_\be$ as defined in \cite{SV}, Definition 3.4. Define for $z \in \cN_{\tetil}$
$$\latil(z) = X^* \latil_\be(\Phi^{-1}(z)).$$
Then $(H \ot K, \io, \latil)$ is a GNS-construction for
$\tetil$. The isomorphism $\Phi$ and \cite{SV}, Propositions 3.10 and 7.1 give $$\lspan\{ (\om_{\eta,\la(b)} \ot \io \ot \io)(\Wtil (\io \ot \al)(X^*)) \; \al(x) \mid \eta \in H, b \in \cT_\vfi, x
\in \Nte \}$$
is a \strong -- norm core for $\latil$ and
\begin{equation} \label{eq115}
\latil \bigl( (\om \ot \io \ot \io)(\Wtil (\io \ot \al)(X^*)) \; \al(x) \bigr) =
X^* (\xi(\om) \ot \late(x))
\end{equation}
for all $\om \in \cI$ and $x \in \Nte$.
By Corollary~\ref{112bis} we can define a unique isometry $\cV: H \ot K \recht H \ot K$ such that
$$\cV (\xi(\om) \ot \late(x)) = \latil \bigl( (\om \ot \io \ot \io)(\Wtil) \al(x) \bigr)$$
for all $\om \in \cI$ and $x \in \Nte$. Later we will show that $\cV=1$ to get the second statement of the proposition.
Let
$$\cD_0 := \lspan\{ (\om_{\eta,\la(b)} \ot \io \ot \io)(\Wtil) \al(x) \mid \eta \in H, b \in \cT_\vfi, x \in \Nte \} \; $$
and $\cD$ the domain of the \strong -- norm closure of the restriction of $\latil$ to
$\cD_0$. Choose $\eta \in H$, $b \in \cT_\vfi$, $x \in \Nte$ and an orthonormal basis $(e_i)_{i \in I}$ in $H$. Define
$$z:=(\om_{\eta,\la(b)} \ot \io \ot \io)(\Wtil (\io \ot \al)(X^*)) \; \al(x)$$
and for every finite subset $I_0 \subset I$
$$z_{I_0} := \sum_{i \in I_0} (\om_{e_i,\la(b)} \ot \io \ot \io)(\Wtil) \; \al \bigl(
(\om_{\eta,e_i} \ot \io)(X^*) x \bigr) \; .$$
Then every $z_{I_0}$ belongs to $\cD_0$ and, using \cite{SV}, Proposition 7.1, we get
\begin{align*}
\latil(z_{I_0}) &= \cV \Bigl( \sum_{i \in I_0} \xi(\om_{e_i,\la(b)}) \ot (\om_{\eta,e_i} \ot
\io)(X^*)\late(x) \Bigr) \\
&= \cV \Bigl( \sum_{i \in I_0} J \si_{i/2}(b) J e_i \ot (\om_{\eta,e_i} \ot \io)(X^*) \late(x)
\Bigr) \; .
\end{align*}
So we observe that $\bigl( \latil(z_{I_0}) \bigr)_{I_0 \subset I}$ converges in norm to
$$\cV \; (J \si_{i/2}(b) J \ot 1) \; X^* \; (\eta \ot \late(x)) =
\cV \; X^* \; (J \si_{i/2}(b) J\eta \ot \late(x)) = \cV \; X^* \; (\xi(\om_{\eta,\la(b)}) \ot \late(x)) \;
.$$
Further $(z_{I_0})_{I_0 \subset I}$ converges in \strong topology to $z$. So we may conclude that $z \in
\cD$ and
$$\latil(z) = \cV \; X^* \; (\xi(\om_{\eta,\la(b)}) \ot \late(x)) \; .$$
In the beginning of the proof
we saw that the elements $z$ span a \strong -- norm core for $\latil$, so $\cD_0$ is a \strong -- norm core for $\latil$
and by Eq.~\eqref{eq115}
$$\latil(z) = X^*(\xi(\om_{\eta,\la(b)}) \ot \late(x)),$$
so $\cV=1$. Hence we have also proved the second statement of the proposition for stabilizable cocycle actions.

Consider now the general case. Choose a \nsf trace $\Tr$ on $B(H)$. Then $\te \ot \Tr$ is a
\nsf weight on $N \ot B(H)$, and $(\al \ot \io, \cU \ot 1)$ is a stabilizable cocycle action
of $(M,\de)$ on $N \ot B(H)$. Writing $\be:=\al \ot \io$, it is clear that $$M
\kruisje{\be,\cU \ot 1} \bigl( N \ot B(H) \bigr) = \bigl(\cros \bigr) \ot B(H)$$ and
$\hat{\be} = \alh \ot \io$. Then it is clear that the dual weight of $\te \ot \Tr$ equals the
tensor product weight $\tetil \ot \Tr$.

Choose now a GNS-construction $(L,\pi_\Tr,\la_\Tr)$ for $\Tr$ and let $(K \ot L,\io \ot
\pi_\Tr,\late \ot \la_\Tr)$ be the canonical GNS-construction for the tensor product weight
$\te \ot \Tr$. By the previous part of the proof we get a GNS-construction $(H \ot K \ot L,
\io \ot \io \ot \pi_\Tr, \latil_1)$ for the dual weight $\tetil \ot \Tr$, such that
\begin{enumerate}
\item $\latil_1 \bigl( (\om \ot \io \ot \io)(\Wtil) \al(x) \ot y \bigr) = \xi(\om) \ot \late(x)
\ot \la_\Tr(y)$ for all $\om \in \cI$, $x \in \Nte$ and $y \in \cN_\Tr$.
\item $\lspan\{ (\om_{\eta,\la(b)} \ot \io \ot \io)(\Wtil) \al(x) \ot y \mid \eta \in H,b \in \cT_\vfi,x
\in \Nte, y \in \cN_\Tr \}$ is a \strong -- norm core for $\latil_1$.
\end{enumerate}
To deduce the second statement above from the first part of the proof, observe that
the elements $x \ot y$ with $x \in \Nte$ and $y \in \cN_\Tr$ span a \strong -- norm core for
$\late \ot \la_\Tr$.

Choose now an arbitrary GNS-construction $(K_2,\pi,\latil_2)$ for the weight $\tetil$. The
lemma following this proposition shows that $$\lspan\{ (\om_{\eta,\la(b)} \ot \io \ot
\io)(\Wtil) \al(x) \mid \eta \in H,b \in \cT_\vfi,x \in \Nte \}$$ is a \strong -- norm core
for $\latil_2$. By Corollary~\ref{112bis} there is a unique isometry $\cV : H \ot K \recht
K_2$ such that $$\cV( \xi(\om) \ot \late(x) ) = \latil_2 \bigl( (\om \ot \io \ot \io)(\Wtil)
\al(x) \bigr)$$ for all $\om \in \cI$ and $x \in \Nte$. Because the elements $(\om \ot \io \ot
\io)(\Wtil) \al(x)$ span a core for $\latil_2$, $\cV$ is unitary and we can define a new
GNS-construction $(H \ot K, \rho, \latil)$ for $\tetil$ by $\rho(z) = \cV^* \pi(z) \cV$ and
$\latil(z) = \cV^* \latil_2(z)$. Of course, the elements $(\om_{\eta,\la(b)} \ot \io \ot
\io)(\Wtil) \al(x)$ still span a core of $\latil$ and $$\latil \bigl( (\om \ot \io \ot
\io)(\Wtil) \al(x) \bigr) = \xi(\om) \ot \late(x)$$ for all $\om \in \cI$ and $x \in \Nte$.
Let $(H \ot K \ot L, \rho \ot \pi_\Tr,\latil \ot \la_\Tr)$ be the tensor product
GNS-construction for the weight $\tetil \ot \Tr$. Then $\latil \ot \la_\Tr$ and $\latil_1$
agree on a core for both $\latil \ot \la_\Tr$ and $\latil_1$. Hence $\latil \ot
\la_\Tr=\latil_1$, from where $\rho \ot \pi_\Tr=\io \ot \io \ot \pi_\Tr$. Then finally
$\rho=\io$, which concludes the proof.
\end{proof}
Now we still have to prove the following lemma.
\begin{lemma} \label{114}
Let $N_1$ and $N_2$ be von Neumann algebras and $\vfi_i$ a \nsf weight on $N_i \; (i=1,2)$. Suppose
that $(K_1,\pi_1,\la_1)$ is a GNS-construction for $\vfi_1$ and that $(K,\pi,\la)$ is a
GNS-construction for the tensor product weight $\vfi_1 \ot \vfi_2$. Suppose further that
$\cD$ is a linear subspace of $\cN_{\vfi_1}$ such that the algebraic tensor product $\cD \odot \cN_{\vfi_2}$ is a \strong
-- norm core for $\la$.
Then $\cD$ is a \strong -- norm core for $\la_1$.
\end{lemma}
\begin{proof}
Choose a GNS-construction $(K_2,\pi_2,\la_2)$ for $\vfi_2$. Because all GNS-constructions for
$\vfi_1 \ot \vfi_2$ are isomorphic, we may assume that
$$(K,\pi,\la)=(K_1 \ot K_2,\pi_1 \ot \pi_2,\la_1 \ot \la_2) \; .$$
Choose now $a \in \cN_{\vfi_1}$. Fix an element $b \in \cN_{\vfi_2}$, $b \neq 0$ and fix
elements $c,d \in \cT_{\vfi_2}$ such that $\vfi_2(d^* b c)=1$.

Choose then a net $(x_\al)_\al$ in $\cD \odot \cN_{\vfi_2}$ such that
$$x_\al \recht a \ot b \tekst{\strong} \quad\text{and}\quad (\la_1 \ot \la_2)(x_\al)
\recht \la_1(a) \ot \la_2(b) \tekst{in norm}.$$
Then $\bigl( (\io \ot \om_{\la_2(c),\la_2(d)})(x_\al) \bigr)_\al$ will be a net in $\cD$,
converging \strong to $a$. Further it is easy to verify that
$$\la_1 \bigl( (\io \ot \om_{\la_2(c),\la_2(d)})(x_\al) \bigr) = (1 \ot \te^*_{\la_2(d
\si^2_i(c)^*)}) (\la_1 \ot \la_2)(x_\al)$$
and this converges in norm to
$$(1 \ot \te^*_{\la_2(d \si^2_i(c)^*)}) (\la_1(a) \ot \la_2(b))= \la_1(a) \; .$$
So, indeed $\cD$ is a \strong -- norm core for $\la_1$.
\end{proof}
\begin{terminology} \label{115}
If $(\al,\cU)$ is a cocycle action of $(M,\de)$ on $N$ and $\te$ is a \nsf weight on $N$ with
GNS-construction $(K,\pite,\late)$, then the GNS-construction $(H \ot K,\io \ot \pite,\latil)$
for $\tetil$ obtained in Proposition~\ref{113} is called the canonical GNS-construction for $\tetil$.
\end{terminology}
\subsection{The unitary implementation of a cocycle action and the biduality theorem}
The results of this subsection follow from Propositions~\ref{17} and~\ref{18} and from similar results for ordinary actions \cite{SV}, so we omit their proofs.
\begin{definition} \label{116}
Let $(\al,\cU)$ be a cocycle action of $(M,\de)$ on $N$ and let $\te$ be a \nsf weight on $N$
with GNS-construction $(K,\pite,\late)$. Let $(H \ot K,\io \ot \pite,\latil)$ be the canonical
GNS-construction for the dual weight $\tetil$. Denote with $\Jte$ and $\Jtil$ the modular
conjugations of $\te$ and $\tetil$ in these GNS-constructions. Then
$$U_\al = \Jtil (\Jh \ot \Jte)$$
is called the unitary implementation of $(\al,\cU)$ obtained by $\te$.
\end{definition}
To justify this terminology let us state the following
proposition.
\begin{proposition} \label{117}
With the notation of the previous definition we have
\begin{enumerate}
\item $(\io \ot \pite)\al(x) = U_\al (1 \ot \pite(x)) U_\al^* \tekst{for all} x \in N \; .$
\item $U_\al \in M \ot B(K) \; .$
\item $(\de \ot \io)(U_\al) = (\io \ot \io \ot \pite)(\cU^*) \; U_{\al \, 23} U_{\al \, 13} \; (\Jh \ot \Jh \ot
\Jte) (\io \ot \io \ot \pite)(\cU_{213}) (\Jh \ot \Jh \ot \Jte) \; .$
\item $(\Jh \ot \Jte)U_\al = U_\al^* (\Jh \ot \Jte) \; .$
\end{enumerate}
\end{proposition}
In fact, the unitary implementation $U_\al$ of $(\al,\cU)$ does not depend on the choice of
the weight $\te$ on $N$. If $\te$ and $\te'$ are both \nsf weights on $N$ with
GNS-constructions $(K,\pite,\late)$ and $(K',\pi'_\te,\la'_\te)$ and if $u$ is the unitary
from $K$ onto $K'$ intertwining $\pite$ and $\pi'_\te$ and mapping the positive cone of $K$
onto the one of $K'$, then $$U'_\al = (1 \ot u) U_\al (1 \ot u^*) \; .$$

Finally we state the biduality theorem.
\begin{proposition} \label{118}
Let $(\al,\cU)$ be a cocycle action of $(M,\de)$ on $N$. Then the formula
$$\Phi : z \mapsto \Wtil (\io \ot \al)(z) \Wtil^*$$
defines a $*$-isomorphism of $B(H) \ot N$ onto $\Mh \kruisje{\alh} (\cros)$ satisfying
\begin{align*}
\Phi(z) &= \alh(z) \tekst{for all} z \in \cros \; , \\
\Phi(x \ot 1) &= x \ot 1 \ot 1 \tekst{for all} x \in M' \; .
\end{align*}
Defining
$$\mu : B(H) \ot N \recht M \ot B(H) \ot N : \mu(z)= (\Si V^* \ot 1)\cU^* (\io \ot \al)(z) \cU (V
\Si \ot 1)$$
we have $(\cJ \ot \Phi) \mu = \alhh \Phi$, where $\cJ(x)=J \Jh x \Jh J$ for all $x \in M$.
\end{proposition}
Let us mention that up to the flip map between $B(H) \ot N$ and $N \ot B(H)$ the action $\mu$
precisely agrees with the stabilized action of $(M,\de)$ on $N \ot B(H)$ given by combining
Propositions~\ref{17} and \ref{18}. Combining the previous proposition with the formula
$$\alh(\cros) = \bigl( \Mh \kruisje{\alh} (\cros) \bigr)^{\alhhklein}$$ we get the following
corollary.
\begin{corollary} \label{119}
Let $(\al,\cU)$ be a cocycle action of $(M,\de)$ on $N$. Then we have
$$\cros = \{ z \in B(H) \ot N \mid (V^* \ot 1)\cU^* (\io \ot \al)(z) \cU (V \ot 1) = z_{13} \} \;
.$$
\end{corollary}
\subsection{Cocycle crossed products and cleft extensions of von Neumann algebras}
Here we answer the following question. Let $(M,\de)$ be a l.c.\ quantum group and let $\te: N
\recht \Mh \ot N$ be an action of $(\Mh,\dehop)$ on the von Neumann algebra $N$. When does $N$
have the structure of a cocycle bicrossed product such that $\te$ agrees with the dual action?
\begin{proposition} \label{120}
Let $N$ be a von Neumann algebra and $(M,\de)$ a l.c.\ quantum group. Suppose that $\te : N
\recht \Mh \ot N$ is an action of $(\Mh,\dehop)$ on $N$. Then the following assertions are
equivalent.
\begin{enumerate}
\item There exists a cocycle action $(\al,\cU)$ of $(M,\de)$ on $N^\te$ such that
$$(N,\te) \cong (M \kruisje{\al,\cU} N^\te, \alh) \; .$$
\item There exists a unitary $X \in M \ot N$ such that
$$(\io \ot \te)(X) = W_{12} X_{13} \; .$$
\end{enumerate}
If the second condition is satisfied the formulas
$$\al(z) = X^* (1 \ot z) X \tekst{for all} z \in N^\te \quad and \quad \cU = X^*_{23}X^*_{13} (\de \ot \io)(X)$$
define a cocycle action $(\al,\cU)$ of $(M,\de)$ on $N^\te$. The formula
$$\pi: N \recht M \kruisje{\al,\cU} N^\te : \pi(z)=X^* \te(z) X$$
defines a $*$-isomorphism satisfying
$$
(\io \ot \pi)(X) = W_{12} \cU^*, \quad\quad
\pi(x) = \al(x) \tekst{for all} x \in N^\te \quad\quad \text{and} \quad\quad \alh \pi = (\io \ot \pi)\te \; .
$$
\end{proposition}
\begin{proof}
Let us first suppose that the second statement is true and take a corresponding unitary $X$. Define
$$\al : N^\te \recht M \ot N : \al(z) = X^* (1 \ot z) X \; .$$
Then we get for all $z \in N^\te$
$$
(\io \ot \te)\al(z) = (\io \ot \te)(X^*) (1 \ot \te(z)) (\io \ot \te)(X)
= X^*_{13} W^*_{12} (1 \ot 1 \ot z) W_{12} X_{13} = \al(z)_{13} \; .
$$
So $\al(N^\te) \subset M \ot N^\te$. Next we define the unitary $\cU$ in $M \ot M
\ot N$ by
$$\cU := X^*_{23} X^*_{13} (\de \ot \io)(X) \; ,$$
from where
\begin{align*}
(\io \ot \io \ot \te)(\cU) &= (\io \ot \te)(X^*)_{234} \; (\io \ot \te)(X^*)_{134} \; (\de \ot
\te)(X)
= X^*_{24} W^*_{23} \; X^*_{14} W^*_{13} \; (\de \ot \io \ot \io)(W_{12} X_{13}) \\
&= X^*_{24} X^*_{14} (\de \ot \io)(X)_{124} = \cU_{124} \; .
\end{align*}
Hence $\cU \in M \ot M \ot N^\te$. To show that $(\al,\cU)$ is a cocycle action of $(M,\de)$
on $N^\te$, we compute, for all $z \in N^\te$: $$ (\de \ot \io)\al(z) = (\de \ot \io)(X^*) (1
\ot 1 \ot z) (\de \ot \io)(X) = \cU^* X^*_{23} X^*_{13} (1 \ot 1 \ot z) X_{13} X_{23} \cU =
\cU^* \; (\io \ot \al)\al(z) \; \cU \; . $$ Further we have
\begin{align*}
(\io \ot \io \ot \al)(\cU) \; (\de \ot \io \ot \io)(\cU) &=
X^*_{34} \; X^*_{24} X^*_{14} (\de \ot \io)(X)_{124} \; X_{34} \; \; \;
X^*_{34} (\de \ot \io)(X^*)_{124} ((\de \ot \io)\de \ot \io)(X) \\
&= X^*_{34} X^*_{24} X^*_{14} ((\de \ot \io)\de \ot \io)(X) \; .
\end{align*}
On the other hand,
\begin{align*}
(1 \ot \cU) \; (\io \ot \de \ot \io)(\cU) &= X^*_{34} X^*_{24} (\de \ot \io)(X)_{234} \;
(\de \ot \io)(X^*)_{234} X^*_{14} ((\io \ot \de) \de \ot \io)(X) \\
&= X^*_{34} X^*_{24} X^*_{14} ((\de \ot \io)\de \ot \io)(X) \; .
\end{align*}
Hence we obtain the equality $$(\io \ot \io \ot \al)(\cU) \; (\de \ot \io \ot \io)(\cU) = (1
\ot \cU) \; (\io \ot \de \ot \io)(\cU) \; .$$ So, $(\al,\cU)$ is a cocycle action of $(M,\de)$
on $N$. We claim that the map $$\pi: N \recht B(H) \ot N : \pi(z) = X^* \te(z) X \; $$ is such
that $\pi(N) \subset M \kruisje{\al,\cU} N^\te$. To show this we use Corollary~\ref{119}.
First let us check that $\pi(N) \subset B(H) \ot N^\te$. We have for all $z \in N$:
\begin{align*}
(\io \ot \te)\pi(z) &= (\io \ot \te)(X^*) \; (\io \ot \te)\te(z) \; (\io \ot \te)(X)
=X_{13}^* W^*_{12} \; (\dehop \ot \io)\te(z) \; W_{12} X_{13} \\
&= X_{13}^* \; \te(z)_{13} \; X_{13} = \pi(z)_{13} \; .
\end{align*}
Hence we get $\pi(N) \subset B(H) \ot N^\te$. To prove our claim it now suffices to observe that
\begin{align*}
(V^* \ot 1) \cU^* \; (\io \ot \al)\pi(z) \; \cU(V \ot 1) &=
(V^* \ot 1) \cU^* \; X^*_{23} X^*_{13} \; \te(z)_{13} \; X_{13} X_{23} \; \cU (V \ot 1) \\
&=(V^* \ot 1) \; (\de \ot \io)(X^*) \; \te(z)_{13} \; (\de \ot \io)(X) \; (V \ot 1) \\
&=X^*_{13} V^*_{12} \; \te(z)_{13} \; V_{12} X_{13} = X^*_{13} \te(z)_{13} X_{13} = \pi(z)_{13}
\end{align*}
for all $z \in N$, where we used the fact that $V \in \Mh' \ot M$ in the last line of the
computation.

Next, for $z \in N^\te$ we clearly have
$$\pi(z) = X^* \te(z) X = X^* (1 \ot z) X = \al(z)$$
and further
$$(\io \ot \pi)(X) = X^*_{23} (\io \ot \te)(X) X_{23} = X^*_{23} W_{12} X_{13} X_{23} = W_{12}
\cU^*$$
so we may conclude that $M \kruisje{\al,\cU} N^\te \subset \pi(N)$. Thus, we have proved that
$\pi(N) = M \kruisje{\al,\cU} N^\te$.

To finish the first part of the proof we only have to verify the formula $\alh \pi = (\io \ot
\pi) \te$, but for $z \in N$ we get
\begin{align*}
(\io \ot \pi)\te(z) &= X^*_{23} \; (\io \ot \te)\te(z) \; X_{23} = X^*_{23} \; (\dehop \ot
\io)\te(z) \; X_{23} \\
&=X^*_{23} (\tilde{V} \ot 1) \te(z)_{23} (\tilde{V}^* \ot 1) X_{23} \\
&=(\tilde{V} \ot 1) (1 \ot \pi(z)) (\tilde{V}^* \ot 1) =\alh(\pi(z)) \; .
\end{align*}
Here we used the notation $\tilde{V} = (J \ot J) \Si W \Si (J \ot J)$ and observe that
$\tilde{V} (1 \ot z) \tilde{V}^* = \dehop(z)$ for all $z \in \Mh$ and $\tilde{V} \in \Mh \ot M'$.
So we have proved that the second statement of the proposition implies the first one.

To prove the converse statement we have to show that for every cocycle
action $(\al,\cU)$ of $(M,\de)$ on $N$ there exists a unitary $X \in M \ot (\cros)$ such that
$$(\io \ot \alh)(X) = W_{12} X_{13} \; .$$
Taking $X = W_{12} \cU^*$ an easy computation proves the result.
\end{proof}
From an algebraic point of view this proposition is related to the notion of a cleft extension
\cite{Mon}, Definition 7.2.1; the definition below is given for the case of Hopf $*$-algebra
actions:

Let $(\Mh,\deh, \hat S, \hat\eps)$ be a Hopf $*$-algebra and let $\te : N \recht \Mh \ot N$ be
an action of $(\Mh,\dehop,\hat{S}^{-1},\hat\eps)$ on a unital $*$-algebra $N$, then the
extension $N^\te \subset N$ is called a cleft extension, if there exists a linear map $\ga:\Mh
\recht N$ such that
\begin{itemize}
\item $\ga$ is a comodule map, which means that
$(\io \ot \ga) \dehop = \te \ga \; .$
\item $\ga$ is convolution unitary, i.e.,
$\mult (\ga^* \ot \ga) \dehop(z) = \hat\eps(z) 1 = \mult (\ga \ot \ga^*)\dehop(z)$ for all $z
\in \Mh$, where $\mult$ is the multiplication map and $\ga^*(z) = \ga(\hat{S}(z^*))^*$.
\end{itemize}

For a similar definition in the framework of l.c.\ quantum groups we use the predual space
$M_*$ of $M$. When we turn $M_*$ into an operator space (see e.g.\ \cite{blecher}), the map
$\lambda : M_* \recht \Mh : \lambda(\om) = (\om \ot \io)(W)$ will be completely contractive.
Next, it is possible to define the comultiplication $\deh$ on the level of $M_*$ using the
extended Haagerup tensor product $\eh$ \cite{Efr-Rua}. Loosely speaking one can give a meaning
in $M_*\eh M_*$ to expressions like $$\sum_{i \in I} \om_i \ot \mu_i$$ where $I$ is an index
set, $(\om_i)$ (resp., $(\mu_i)$) is a bounded infinite row (resp., column) over $M_*$. Then
one can show that there exists a unique complete contraction from $M_*$ to $M_* \eh M_*$, also
denoted by $\dehop$, such that $$(\lambda \eh \lambda) \dehop = \dehop \lambda$$ where
$\lambda \eh \lambda$ is the obvious completely contractive embedding of $M_* \eh M_*$ into
$\Mh \ot \Mh$. We then have $$\dehop(\om_{\xi,\eta}) = \sum_{i \in I} \om_{e_i,\eta} \ot
\om_{\xi,e_i}$$ for all $\xi,\eta \in H$ and for an orthonormal basis  $(e_i)_{i \in I}$ in
$H$. The following definition looks very much like the Hopf algebraic definition above, but
has an exact operator algebraic meaning.
\begin{definition} \label{121}
Let $N$ be a von Neumann algebra, $(M,\de)$ a l.c.\ quantum group and $\te : N \recht \Mh \ot
N$ an action of $(\Mh,\dehop)$ on $N$. Then the extension $N^\te \subset N$ is called a cleft
extension  of von Neumann algebras if there exists a complete contraction $$\ga : M_* \recht
N$$ such that
\begin{itemize}
\item $\ga$ is a comodule map, which means that
$$(\lambda \eh \ga) \dehop(\om) = \te(\ga(\om)) \tekst{for all} \om \in M_*$$
where $\lambda \eh \ga$ is the natural complete contraction from $M_* \eh M_*$ to $\Mh \ot N$.
\item $\ga$ is convolution unitary, which means that
$$\mult (\ga^* \eh \ga) \dehop(\om) = \om(1) \, 1 = \mult (\ga \eh \ga^*)\dehop(\om) \tekst{for
all} \om \in M_*$$
where $\ga^*(\om) = \ga(\bar{\om})^*$ for all $\om \in M_*$, $\ga^* \eh \ga$ maps $M_* \eh M_*$
into $N \eh N$ completely contractively and where $\mult : N \eh N \recht N$ is the completely contractive multiplication map.
\end{itemize}
\end{definition}

The following result sheds some light on Proposition~\ref{120} from an algebraic point of
view.

\begin{proposition} \label{122}
Let $N$ be a von Neumann algebra, $(M,\de)$ a l.c.\ quantum group and $\te : N \recht \Mh \ot
N$ an action of $(\Mh,\dehop)$ on $N$. Then the extension $N^\te \subset N$ is a cleft
extension  of von Neumann algebras if and only if one of the equivalent statements of
Proposition~\ref{120} is satisfied.
\end{proposition}
\begin{proof}
Denoting by $\ot_{\text{max}}$ the maximal operator space tensor product, it follows from
\cite{blecher-paulsen}, Proposition~5.4 that the completely bounded maps from $M_*$ to $N$ can
be identified with $(M_* \ot_{\text{max}} N_*)^*$. Combining this with \cite{blecher2},
Theorem~2.5 it follows that the completely bounded maps $\ga$ from $M_*$ to $N$ are in
(completely isometric) one-to-one correspondence with elements $X \in M \ot N$ by the formula
$$\ga(\om) = (\om \ot \io)(X) \tekst{for all} \om \in M_* \; .$$

Let first one of the statements of Proposition~\ref{120} be satisfied. Then
 take such a unitary $X \in M \ot N$ satisfying $(\io \ot \te)(X) = W_{12} X_{13}$. Define
$\ga(\om) = (\om \ot \io)(X)$ for all $\om \in M_*$. Let $(e_i)_{i \in I}$ be an orthonormal
basis for $H$ and $\xi,\eta \in H$. Then
\begin{align*}
(\lambda \eh \ga)\dehop(\om_{\xi,\eta}) &= \sum_{i \in I} \lambda(\om_{e_i,\eta}) \ot
\ga(\om_{\xi,e_i}) \\
&= \sum_{i \in I} (\om_{e_i,\eta} \ot \io)(W) \ot (\om_{\xi,e_i} \ot \io)(X) \\
&= (\om_{\xi,\eta} \ot \io \ot \io)(W_{12} X_{13}) = \te \bigl( (\om_{\xi,\eta} \ot \io)(X) \bigr) \\
&= \te \bigl( \ga(\om_{\xi,\eta})  \bigr) \; ,
\end{align*}
so
$$(\lambda \eh \ga) \dehop(\om) = \te(\ga(\om)) \tekst{for all} \om \in M_* \; .$$
Further we have
\begin{align*}
\mult (\ga^* \eh \ga) \dehop(\om_{\xi,\eta}) &= \sum_{i \in I} \ga^*(\om_{e_i,\eta})
\ga(\om_{\xi,e_i}) \\
&= \sum_{i \in I} (\om_{e_i,\eta} \ot \io)(X^*) (\om_{\xi,e_i} \ot \io)(X) \\
&= (\om_{\xi,\eta} \ot \io)(X^* X)=\om_{\xi,\eta}(1) \, 1 \; .
\end{align*}
Thus,
$$\mult (\ga^* \eh \ga) \dehop(\om) = \om(1) \, 1 \tekst{for all} \om \in M_*,$$
which concludes the first part of the proof, because the remaining equality can be obtained in a similar way.

Vice versa, let us have a map $\ga$. Then there exists an element
$X \in M \ot N$ such that $\ga(\om) = (\om \ot \io)(X)$ for all $\om \in M_*$. The same
computations as above prove that $(\io \ot \te)(X)=W_{12} X_{13}$, and that $X$ is a unitary.
\end{proof}

\section{Cocycle bicrossed products of l.c.\ quantum groups}
We introduce cocycle matched pairs of l.c.\ quantum groups, give the cocycle bicrossed product
construction for them and show that the result is again a  l.c.\ quantum group.
\subsection{Cocycle matched pairs of l.c.\ quantum groups}
\begin{definition} \label{21}
A pair $(M_1,\de_1),(M_2,\de_2)$ is said to be a matched pair of l.c.\ quantum groups if a
triple $(\tau, \cU,\cV)$ (called a cocycle matching) satisfies the following conditions:
\begin{align*}
& \cU \in M_1 \ot M_1 \ot M_2 \tekst{and} \cV \in M_1 \ot M_2 \ot M_2 \tekst{are both
unitaries and}\\ & \tau : M_1 \ot M_2 \recht M_1 \ot M_2 \tekst{is a faithful
$*$-homomorphism;}
\end{align*}
denoting by $\si$ the flip map and defining $$\al : M_2 \recht M_1 \ot M_2 : \al(y) = \tau(1
\ot y) \quad\text{and}\quad \be : M_1 \recht M_1 \ot M_2 : \be(x) = \tau(x \ot 1)$$ we have
\begin{itemize}
\item $(\al,\cU)$ is a cocycle action of $(\Mo,\deo)$ on the von Neumann algebra $\Mt$ :
\begin{align*}
(\io \ot \al)\al(y) &= \cU (\de_1 \ot \io)\al(y) \cU^* \\
(\io \ot \io \ot \al)(\cU) (\de_1 \ot \io \ot \io)(\cU) &= (1 \ot \cU) (\io \ot
\de_1 \ot \io)(\cU) \; .
\end{align*}
\item $(\si \be, \cV_{321})$ is a cocycle action of $(\Mt,\det)$ on the von Neumann algebra
$\Mo$ :
\begin{align*}
(\be \ot \io)\be(x) &= \cV (\io \ot \deopt)\be(x) \cV^* \\
(\be \ot \io \ot \io)(\cV)
(\io \ot \io \ot \deopt)(\cV) &= (\cV \ot 1) (\io \ot \deopt \ot \io)(\cV) \; .
\end{align*}
\item $(\al,\cU)$ and $(\be,\cV)$ are matched in the following sense:
\begin{align*}
\tau_{13} (\al \ot \io) \de_2(y) &= \cV_{132} (\io \ot \de_2)\al(y) \cV_{132}^* \\
\tau_{23}\si_{23}(\be \ot \io) \de_1(x) &= \cU (\de_1 \ot \io)\be(x) \cU^* \\
(\de_1 \ot \io \ot \io)(\cV)(\io \ot \io \ot \deopt)(\cU^*) &= (\cU^* \ot 1) (\io \ot
\tau \si \ot \io) \bigl( (\be \ot \io \ot \io)(\cU^*) (\io \ot \io
\ot \al)(\cV) \bigr) (1 \ot \cV) \; .
\end{align*}
\end{itemize}
\end{definition}
If both $(\Mo,\deo)$ and $(\Mt,\det)$ are commutative, i.e., are generated by l.c.\ groups,
the above conditions split into two parts: 1) we have two ordinary actions $\al$ and $\be$
which are matched:
\begin{align*}
(\io \ot \al)\al(y) &= (\de_1 \ot \io)\al(y) \\
(\be \ot \io)\be(x) &= (\io \ot \deopt)\be(x) \\
\tau_{13} (\al \ot \io) \de_2(y) &= (\io \ot \de_2)\al(y) \\
\tau_{23}\si_{23}(\be \ot \io) \de_1(x) &=  (\de_1 \ot \io)\be(x) \\
\intertext{2) we have two cocycles $\cU$ and $\cV$ which are matched:}
(\io \ot \io \ot \al)(\cU) (\de_1 \ot \io \ot \io)(\cU) &= (1 \ot \cU) (\io \ot
\de_1 \ot \io)(\cU) \\
(\be \ot \io \ot \io)(\cV)
(\io \ot \io \ot \deopt)(\cV) &= (\cV \ot 1) (\io \ot \deopt \ot \io)(\cV) \\
(\de_1 \ot \io \ot \io)(\cV)(\io \ot \io \ot \deopt)(\cU^*) &= (\cU^* \ot 1) (\io \ot
\tau \si \ot \io) \bigl( (\be \ot \io \ot \io)(\cU^*) (\io \ot \io
\ot \al)(\cV) \bigr) (1 \ot \cV) \; .
\end{align*}
So to construct examples of cocycle matched pairs of l.c.\ groups, we first have to construct
matched pairs without cocycles, and then to solve the above cocycle equations. Sections 4 and
5 are devoted to the analysis of this special case. Let us also  stress that
Definition~\ref{21} is an operator algebraic version of the Hopf algebraic definition
\cite{Majbook}, Chapter 6.

Suppose that $(M_1,\de_1)$ and $(M_2,\de_2)$ are l.c.\ quantum groups and that $\tau : M_1 \ot
M_2 \recht M_1 \ot M_2$ is a faithful $^*$-homomorphism. Then we put $\tau' = \tau \si$. It is
easy to check that $(\tau,1,1)$ is a cocycle matching of $(M_1,\de_1)$ and $(M_2,\de_2)$ (with
trivial cocycles) if and only if $$(\tau' \ot \io) (\io \ot \tau')(\deopt \ot \io) = (\io \ot
\deopt) \tau' \qquad\text{and}\qquad (\io \ot \tau')(\tau' \ot \io)(\io \ot \de_1) = (\de_1
\ot \io) \tau' \; .$$ So, a matching $\tau$ of $(M_1,\de_1)$ and $(M_2,\de_2)$ with
trivial cocycles is a von Neumann algebraic version of the notion of an inversion of
$(M_2,\deopt)$ and $(M_1,\de_1)$, introduced in \cite{B-S1}, D\'efinition~8.1. Moreover, in
the notation of \cite{B-S1}, $(M_2,\deopt)$ is the von Neumann algebraic
counterpart of $(S_{\Wh_2},\sde_{\Wh_2})$ and $(M_1,\de_1)$ is the von Neumann algebraic
counterpart of $(\hat{S}_{W_1},\sdeh_{W_1})$, and hence the multiplicative unitary $T$ of
\cite{B-S1}, Proposition~8.7 agrees with the multiplicative unitary $\Wh$ that we
will define in Definition~\ref{22}. Observe that our framework is wider than the one of \cite{B-S1}, so we do not need to impose extra conditions on $\tau$ as it was done in \cite{B-S1}.

\subsection{Cocycle bicrossed products}\label{sec22}
Let $(\tau,\cU,\cV)$ be a cocycle matching of $(\Mo,\deo)$ and $(\Mt,\det), \al$ and $\be$ as
in Definition~\ref{21}. Let $W_1$,$(H_1,\io,\la_1)$ and $W_2$, $(H_2,\io,\la_2)$ denote the
multiplicative unitaries and GNS-constructions of $(\Mo,\deo)$ and $(\Mt,\det)$ respectively.
We write $H=H_1\ot H_2$ and denote by $\Sigma$ the flip map on $H \ot H$. As in Notation 1.2,
we also write $\Wtil = (W_1 \ot 1) \cU^* \; .$ Let us define the major ingredients of the
cocycle bicrossed product $(M,\de)$.
\begin{definition} \label{22}
Define unitaries $W$ and $\Wh$ on $H \ot H$ by
\begin{align*}
\Wh &= (\be \ot \io \ot \io) \bigl( (W_1 \ot 1) \cU^* \bigr) \; (\io \ot \io \ot \al) \bigl(
\cV (1 \ot \Wht) \bigr)\quad and\quad W = \Si \Wh^* \Si \; .
\end{align*}
Let $M = M_1 \kruisje{\al,\cU} M_2$ be the von Neumann subalgebra of $B(H_1)\ot M_2$ generated by
$$\al(M_2) \quad\text{and}\quad \{ (\om \ot \io \ot \io) \bigl( (W_1 \ot 1) \cU^* \bigr) \mid
\om \in M_{1 \, *} \}.$$
Further we define the faithful $*$-homomorphism
$$\de : M \recht B(H \ot H) : \de(z) = W^*(1 \ot z) W \; .$$
\end{definition}

In the following propositions we will gradually prove that $(M,\de)$ is a l.c.\ quantum group,
that $W$ is the associated multiplicative unitary and we will describe the dual l.c.\ quantum
group $(\Mh,\deh)$.

The proof of the following formula is obvious.
\begin{lemma}\label{23}
We have $M \subset B(H_1) \ot M_2$ and for all $z \in M$ we have
$$\deop(z) = (\be \ot \io \ot \io)(\Wtil) \; (\io \ot \io \ot \al) \bigl( \cV (\io \ot \deopt)(z)
\cV^* \bigr) \; (\be \ot \io \ot \io)(\Wtil^*) \; . $$
\end{lemma}
The definition of $M$ makes it natural to
give formulas for $\de(\al(x))$ and $(\io \ot \de)(\Wtil)$.
\begin{proposition} \label{24}
For all $x \in M_2$ we have
$$\de(\al(x)) = (\al \ot \al) \det(x) \; .$$
\end{proposition}
\begin{proof}
Choose $x \in M_2$.
We will prove that $\deop(\al(x)) = (\al \ot \al)\deopt(x)$.
Observe that
\begin{align*}
(\io \ot \io \ot \al) \bigl( \cV \; (\io \ot \deopt)\al(x) \; \cV^* \bigr) &=
(\io \ot \io \ot \al) \si_{23} \tau_{13} (\al \ot \io) \det(x) \\
&= (\tau \ot \io \ot \io) \bigl( (\io \ot \al \ot \io)(\al \ot \io) \det(x)_{1342} \bigr) \\
&= (\tau \ot \io \ot \io) \bigl( \; \bigl( (\Wtil^* \ot 1) \; (1 \ot (\al \ot \io)\det(x) ) \; (\Wtil \ot
1) \bigr)_{1342} \; \bigr) \\
&= (\be \ot \io \ot \io)(\Wtil^*) \; (\al \ot \al) \deopt(x) \; (\be \ot \io \ot \io)(\Wtil) \; .
\end{align*}
Using the previous lemma we immediately get that $\deop(\al(x)) = (\al \ot \al)\deopt(x)$.
\end{proof}
\begin{proposition} \label{25}
We have
$$(\io \ot \deop)(\Wtil) = (\Wtil \ot 1 \ot 1) \; \bigl( (\io \ot \al) \be \ot \io \ot \io \bigr)
(\Wtil) \; (\io \ot \al \ot \al)(\cV) \; .$$
\end{proposition}
\begin{proof}
Using the lemma we get
$$(\io \ot \deop)(\Wtil) = \bigl( 1 \ot (\be \ot \io \ot \io)(\Wtil) \bigr) \; (\io \ot \io \ot
\io \ot \al) \bigl( (1 \ot \cV)(\io \ot \io \ot \deopt)(\Wtil) (1 \ot \cV^*) \bigr) \;
\bigl( 1 \ot (\be \ot \io \ot \io)(\Wtil^*) \bigr) \; .$$
The definition of $\Wtil$ and the fact that $\cU$ and $\cV$ are matched imply
\begin{align*}
(1 \ot \cV) & (\io \ot \io \ot \deopt)(\Wtil) (1 \ot \cV^*) \\
&= (1 \ot \cV) (W_1 \ot 1 \ot 1) (\io \ot \io \ot \deopt)(\cU^*) (1 \ot \cV^*) \\
&= (W_1 \ot 1 \ot 1) (\deo \ot \io \ot \io)(\cV) (\io \ot \io \ot \deopt)(\cU^*) (1 \ot \cV^*)
\\
&= (W_1 \ot 1 \ot 1)(\cU^* \ot 1) \; (\io \ot \tau\si \ot \io) \bigl( (\be \ot \io \ot
\io)(\cU^*) (\io \ot \io \ot \al)(\cV) \bigr) \; .
\end{align*}
From this we may conclude that
\begin{align*}
(\io \, \ot \, & \io \ot \io \ot \al)  \bigl( (1 \ot \cV)(\io \ot \io \ot \deopt)(\Wtil)
(1 \ot \cV^*) \bigr) \\
&= (\Wtil \ot 1 \ot 1) \; \tau_{23} \si_{23} \bigl( (\be \ot \io \ot \al)(\cU^*) (\io \ot \io
\ot (\io \ot \al)\al )(\cV) \bigr) \\
&= (\Wtil \ot 1 \ot 1) \; \tau_{23} \si_{23} \bigl( (\be \ot \io \ot \al)(\cU^*) (1 \ot 1 \ot
\Wtil^*) (\io \ot \io \ot \al)(\cV)_{1245} (1 \ot 1 \ot \Wtil) \bigr) \\
&= (\Wtil \ot 1 \ot 1) \; \tau_{23} \si_{23} (\be \ot \io \ot \io \ot \io) \bigl( (\io \ot \io \ot
\al)(\cU^*)(1 \ot \Wtil^*) \bigr) \;
(\io \ot \al \ot \al)(\cV) \; (1 \ot (\be \ot
\io \ot \io)(\Wtil)) \; .
\end{align*}
Using the first formula of the proof we get
\begin{align}
(\io \ot \deop)(\Wtil) & =
(1 \ot (\be \ot \io \ot \io)(\Wtil)) \; (\Wtil \ot 1 \ot 1) \;
\tau_{23} \si_{23} (\be \ot \io \ot \io \ot \io) \bigl( (\io \ot \io \ot
\al)(\cU^*)(1 \ot \Wtil^*) \bigr) \notag
\\ & \qquad\qquad (\io \ot \al \ot \al)(\cV) \; . \label{hulp1}
\end{align}
Because for all $x \in M_1$ we have
$$ \Wtil^* (1 \ot \be(x)) \Wtil = \tau_{23}\si_{23} (\be \ot \io)\deo(x),$$
we get
$$ (1 \ot (\be \ot \io \ot \io)(\Wtil)) \; (\Wtil \ot 1 \ot 1) = (\Wtil \ot 1 \ot 1) \; \tau_{23}
\si_{23} (\be \ot \io \ot \io \ot \io) (\deo \ot \io \ot \io)(\Wtil) \; .$$
Combining this with Eq.~\eqref{hulp1} we get
$$(\io \ot \deop)(\Wtil) = (\Wtil \ot 1 \ot 1) \; \tau_{23}
\si_{23} (\be \ot \io \ot \io \ot \io) \bigl( (\deo \ot \io \ot \io)(\Wtil) (\io \ot \io \ot
\al)(\cU^*) (1 \ot \Wtil^*) \bigr) \; (\io \ot \al \ot \al)(\cV) \; .$$
Finally
\begin{align*}
(\deo \ot \io \ot \io)(\Wtil) (\io \ot \io \ot \al)(\cU^*) (1 \ot \Wtil^*) &=
W_{1,13} W_{1,23} (\deo \ot \io \ot \io)(\cU^*)(\io \ot \io \ot \al)(\cU^*) (1 \ot \Wtil^*) \\
&= W_{1,13} W_{1,23} (\io \ot \deo \ot \io)(\cU^*)(1 \ot \cU^* \Wtil^*) \\
&= W_{1,13} W_{1,23} (\io \ot \deo \ot \io)(\cU^*) W^*_{1,23} \\
&= W_{1,13} \cU^*_{134} = \Wtil_{134}
\end{align*}
and hence we get
$$(\io \ot \deop)(\Wtil) = (\Wtil \ot 1 \ot 1) \; \bigl( (\io \ot \al) \be \ot \io \ot \io \bigr)
(\Wtil) \; (\io \ot \al \ot \al)(\cV) \; .$$
\end{proof}
\begin{corollary} \label{26}
\begin{itemize}
\item $\de(M) \subset M \ot M$.
\item $(\de \ot \io)(W) = W_{13} W_{23}$.
\item $(\io \ot \de) \de = (\de \ot \io) \de$.
\item The unitaries $W$ and $\Wh$ defined on $H \ot H$ are multiplicative.
\end{itemize}
\end{corollary}
\begin{proof}
The first inclusion follows from Propositions~\ref{24} and \ref{25} and from the definition of $M$.

The second equality is equivalent to $(\io \ot
\deop)(\Wh) = \Wh_{12} \Wh_{13}$. Propositions~\ref{24} and \ref{25} imply
\begin{multline*}
(\io \ot \io \ot \deop)(\Wh) = (\be \ot \deop)(\Wtil) \; (\io \ot \io \ot \al \ot \al)(\io \ot
\io \ot \deopt) ( \cV (1 \ot \Wht)) \\ = (\be \ot \io \ot \io)(\Wtil)_{1234} \; \bigl( (\be
\ot \al)\be \ot \io \ot \io \bigr)(\Wtil) \; (\be \ot \al \ot \al)(\cV) \; (\io \ot \io \ot
\al \ot \al) \bigl( (\io \ot \io \ot \deopt)(\cV) \hat{W}_{2,23} \hat{W}_{2,24} \bigr) \; .
\end{multline*}
But we know that
\begin{align*}
(\be \ot \io \ot \io)(\cV) \; (\io \ot \io \ot \deopt)(\cV) \; \hat{W}_{2,23} \;
\hat{W}_{2,24} &= (\cV \ot 1) \; (\io \ot \deopt \ot \io)(\cV) \; \hat{W}_{2,23} \;
\hat{W}_{2,24} \\ &= (\cV \ot 1) \; \hat{W}_{2,23} \; \cV_{124} \; \hat{W}_{2,24} \; .
\end{align*}
So,
\begin{multline} \label{hulp2}
(\io \ot \io \ot \deop)(\Wh) \\ =
(\be \ot \io \ot \io)(\Wtil)_{1234} \; \bigl( (\be \ot \al)\be \ot \io \ot \io \bigr)(\Wtil) \;
(\io \ot \io \ot \al)(\cV (1 \ot \Wht))_{1234} \; (\io \ot \io \ot \al)(\cV (1 \ot \Wht))_{1256} \;
.
\end{multline}
Finally we use the equality
$$\bigl( (\be \ot \io)\be \ot \io \ot \io \bigr)(\Wtil) = (\cV \ot 1 \ot 1) \hat{W}_{2,23} (\be
\ot \io \ot \io)(\Wtil)_{1245} \hat{W}^*_{2,23} (\cV^* \ot 1 \ot 1)$$
to conclude that
$$\bigl( (\be \ot \al)\be \ot \io \ot \io \bigr)(\Wtil) =
(\io \ot \io \ot \al) (\cV (1 \ot \Wht))_{1234} (\be \ot \io \ot \io)(\Wtil)_{1256}
(\io \ot \io \ot \al) (\cV (1 \ot \Wht))_{1234}^* \; .$$
Combining this with Eq.~\eqref{hulp2} we get
$$(\io \ot \deop)(\Wh) = \Wh_{12} \Wh_{13} \; .$$

The third and fourth statement of the corollary immediately follow from the second one.
\end{proof}
So at the moment we have at our disposal a von Neumann algebra $M$, a coassociative
comultiplication $\de : M \recht M \ot M$ and a multiplicative unitary $W$. Next
we define a left invariant weight $\vfi$ on $M$ and show that
$W$ is the associated multiplicative unitary (the left regular representation).
\begin{definition} \label{27}
We denote by $\vfi$ the dual weight on the cocycle crossed product $M = M_1 \kruisje{\al,\cU}
M_2$ of the weight $\vfi_2$ on $M_2$, as defined in Definition~\ref{112} and by $(H,\io,\la)$
the canonical GNS-construction for $\vfi$ as defined in Terminology~\ref{115}, given the
GNS-construction $(H_2,\io,\la_2)$ for $\vfi_2$. This means that
\begin{itemize}
\item $\lspan\{ (\om_{\eta,\la_1(b)} \ot \io \ot \io)(\Wtil) \al(x) \mid \eta \in H_1,b \in
\cT_{\vfi_1}, x \in \cN_{\vfi_2} \}$ is a
\strong -- norm core for $\la$.
\item  $\la \bigl( (\om \ot \io \ot \io)(\Wtil) \al(x) \bigr) = \xi_1(\om) \ot \la_2(x)$ for
all $\om \in \cI_1$ and $x \in \cN_{\vfi_2}$.
\end{itemize}
\end{definition}
The following result is important.
\begin{proposition} \label{28}
The weight $\vfi$ is left invariant, which means that
$$(\io \ot \vfi)\de(z) = \vfi(z) \; 1 \tekst{for all} z \in \Mfi^+ \; .$$
Further we have, for all $z \in \Nfi$ and $\om \in M_*$:
$$\la \bigl( (\om \ot \io)\de(z)  \bigr) = (\om \ot \io)(W^*) \la(z).$$
\end{proposition}
\begin{proof}
It suffices to prove the second equality, then the first one follows as usual.
Choose $\mu \in H_1$, $b \in \cT_{\vfi_1}$ and $x \in \cN_{\vfi_2}$. Let $(e_i)_{i \in I}$ be
an orthonormal basis for $H_1$. It follows from Propositions~\ref{24} and \ref{25} that, with
\strong convergence
\begin{align*}
\deop \bigl( (\om_{\mu,\la_1(b)} \ot \io \ot \io) & (\Wtil) \al(x) \bigr)=
\sum_{i \in I} \bigl( (\om_{e_i,\la_1(b)} \ot \io \ot \io)(\Wtil) \ot 1 \ot 1 \bigr) \\
& (\al \ot \io \ot \io) \bigl( (\om_{\mu,e_i} \ot \io \ot \io \ot \io) \bigl( (\be \ot \io \ot
\io)(\Wtil) (\io \ot \io \ot \al)(\cV) \bigr) (\io \ot \al)\deopt(x) \bigr) \; .
\end{align*}
Next choose $\eta,\rho \in H$ and an orthonormal basis $(f_j)_{j \in J}$ for $H$. Then with \strong convergence we have
\begin{align*}
& \qquad\qquad (\io \ot \io \ot \om_{\eta,\rho}) \deop \bigl( (\om_{\mu,\la_1(b)}
\ot \io \ot \io)(\Wtil) \al(x) \bigr) \\
&= \sum_{(i,j) \in I \times J} (\om_{e_i,\la_1(b)} \ot \io \ot \io)(\Wtil) \;
\al \bigl( (\om_{\mu,e_i} \ot \io \ot \om_{f_j,\rho}) \bigl( (\be \ot \io \ot \io)(\Wtil) (\io
\ot \io \ot \al)(\cV) \bigr) (\io \ot \om_{\eta,f_j}\al)\deopt(x) \bigr) \; .
\end{align*}
Denote by $z_{I_0 \times J_0}$ the above expression summed over $I_0 \times J_0$ where $I_0
\subset I$ and $J_0 \subset J$ are finite subsets. Then it follows from Definition~\ref{27}
and \cite{SV}, Proposition 7.1 that all $z_{I_0 \times J_0}$ belong to $\Nfi$ and $$\la(z_{I_0
\times J_0}) = \sum_{(i,j) \in I_0 \times J_0} J_1 \si^1_{i/2}(b) J_1 e_i \ot (\om_{\mu,e_i}
\ot \io \ot \om_{f_j,\rho}) \bigl( (\be \ot \io \ot \io)(\Wtil) (\io \ot \io \ot \al)(\cV)
\bigr) \; (\io \ot \om_{\eta,f_j}\al)(\Wht) \la_2(x) \; .$$ Taking the limit over $J_0 \subset
J$ and using the fact that $\la$ is \strong -- norm closed we get that for all finite subsets
$I_0 \subset I$ the element $z_{I_0}$ defined by $$z_{I_0} := \sum_{i \in I_0}
(\om_{e_i,\la_1(b)} \ot \io \ot \io)(\Wtil) \; \al \bigl( (\om_{\mu,e_i} \ot \io \ot
\om_{\eta,\rho}) \bigl( (\be \ot \io \ot \io)(\Wtil) (\io \ot \io \ot \al)(\cV) (1 \ot (\io
\ot \al)\deopt(x)) \bigr) \bigr)$$ belongs to $\Nfi$ and $$\la(z_{I_0}) = \sum_{i \in I_0} J_1
\si^1_{i/2}(b) J_1 e_i \ot (\om_{\mu,e_i} \ot \io \ot \om_{\eta,\rho})(\Wh) \la_2(x) \; .$$
Taking the limit over $I_0 \subset I$ we get that the element $z$ defined by $$z:=(\io \ot \io
\ot \om_{\eta,\rho}) \deop \bigl( (\om_{\mu,\la_1(b)} \ot \io \ot \io)(\Wtil) \al(x) \bigr)$$
belongs to $\Nfi$ and $$ \la(z)=(J_1 \si^1_{i/2}(b) J_1 \ot 1) (\io \ot \io \ot
\om_{\eta,\rho})(\Wh)(\mu \ot \la_2(x)) \; .$$ Because $\Wh \in M_1 \ot B(H_2 \ot H)$ we get
$$\la(z) = (\io \ot \io \ot \om_{\eta,\rho})(\Wh) (\xi_1(\om_{\mu,\la_1(b)}) \ot \la_2(x)) =
(\io \ot \io \ot \om_{\eta,\rho})(\Wh) \la \bigl( (\om_{\mu,\la_1(b)} \ot \io \ot \io)(\Wtil)
\al(x) \bigr).$$ By Definition~\ref{27} the elements $(\om_{\mu,\la_1(b)} \ot \io \ot
\io)(\Wtil) \al(x)$ span a core for $\la$. So we can conclude that for all $\eta,\rho \in H$
and all $z \in \Nfi$ the element $(\io \ot \om_{\eta,\rho})\deop(z)$ belongs to $\Nfi$ and
$$\la \bigl( (\io \ot \om_{\eta,\rho})\deop(z) \bigr) = (\io \ot \om_{\eta,\rho})(\Wh) \la(z)
\; .$$ From this it follows that for all $z \in \Nfi$ and $\om \in M_*$ we have $(\om \ot
\io)\de(z) \in \Nfi$ and $$\la \bigl( (\om \ot \io) \de(z) \bigr) = (\om \ot \io)(W^*) \la(z)
\; .$$
\end{proof}
If we would already know that $(M,\de)$ is a l.c.\ quantum group there would exist a
$*$-anti-automorphism $R:z \mapsto \Jh z^* \Jh$ of $M$ (where $\Jh$ is the modular conjugation
of the left invariant weight on the dual l.c.\ quantum group $(\Mh,\deh)$)  satisfying
\begin{equation} \label{ster}
\de R = \si (R \ot R) \de \; .
\end{equation}
So, because $\vfi$ is left invariant, we get that $\vfi R$ is right invariant. Thus we will
first construct this dual l.c.\ quantum group and its modular conjugation $\Jh$ and then prove
Eq.~\eqref{ster}.

For a cocycle matching $(\tau,\cU,\cV)$ of $(\Mo,\deo)$ and $(\Mt,\det)$ define $\tautil = \si
\tau \si$, $\cUtil = \cV_{321}$ and $\cVtil = \cU_{321}$. Then one can check that
$(\tautil,\cUtil,\cVtil)$ is a cocycle matching of $(\Mt,\det)$ and $(\Mo,\deo)$. As above we
can then define a von Neumann algebra with comultiplication, a multiplicative unitary (both
acting on $H_2 \ot H_1$) and a left invariant weight with GNS-construction. Using the flip map
from $H_2 \ot H_1$ onto $H=H_1 \ot H_2$ one has the following result.
\begin{proposition} \label{29}
Define $\Mh$ as the von Neumann subalgebra of $M_1 \ot B(H_2)$ generated by $$\be(M_1)
\quad\text{and}\quad \{ (\io \ot \io \ot \om)(\cV (1 \ot \Wht)) \mid \om \in M_{2 \, *} \} \;
.$$ For all $z \in \Mh$, define $$\deh(z) = \Wh^* (1 \ot z) \Wh \; .$$ Then, $\deh$ is a
normal and faithful $*$-homomorphism of $\Mh$ into $\Mh \ot \Mh$ satisfying $$(\io \ot \deh)
\deh = (\deh \ot \io)\deh \; .$$ There exists on $\Mh$ a \nsf weight $\vfih$ with
GNS-construction $(H,\io,\lah)$ such that
\begin{itemize}
\item $\lspan\{ (\io \ot \io \ot \om_{\mu,\la_2(b)})( (1 \ot \Wht^*) \cV^*) \; \be(x) \mid \mu \in H_2, b
\in \cT_{\vfi_2}, x \in \cN_{\vfi_1} \}$ is a \strong -- norm core for $\lah$.
\item $\lah \bigl( (\io \ot \io \ot \om)( (1 \ot \Wht^*) \cV^*) \; \be(x) \bigr) = \la_1(x) \ot
\xi_2(\om)$ for all $\om \in \cI_2$ and $x \in \cN_{\vfi_1}$.
\item $\vfih$ is left invariant on $(\Mh,\deh)$ and we have,
for all $z \in \Nfih$ and $\om \in \Mh_*$:
$$\lah \bigl( (\om \ot \io)\deh(z) \bigr) = (\io \ot \om)(W) \lah(z).$$
\end{itemize}
\end{proposition}
Later on we will prove that $(\Mh,\deh)$ is the dual of $(M,\de)$. The following two lemmas
are known to be valid when $(M,\de)$ is indeed a l.c.\ quantum group with dual $(\Mh,\deh)$,
see \cite{KV2}, Corollary 2.2 and \cite{KV1}, equation (4.2).
\begin{lemma} \label{210}
Let $J$ and $\Jh$ be the modular conjugations of the weights $\vfi$ and $\vfih$
in the GNS-constructions $(H,\io,\la)$ and $(H,\io,\lah)$. Then
$$(\Jh \ot J) W = W^* (\Jh \ot J) \; .$$
\end{lemma}
\begin{proof}
Let us introduce the notations $\nab$ and $\nabh$ for the modular operators of the weights
$\vfi$ and $\vfih$. We first claim that for all $\om \in B(H)_*$
\begin{equation} \label{hulpje3}
(\om \ot \io)(W^*) J \nab^{1/2} \subset J \nab^{1/2} (\bar{\om} \ot \io)(W^*) \; .
\end{equation}
To prove this claim, choose $\om \in B(H)_*$ and $z \in \Nfi \cap \Nfi^*$. Then, using
Proposition~\ref{28}, one gets
$$
(\bar{\om} \ot \io)(W^*) \la(z) = \la ((\bar{\om} \ot \io) \de(z)) \; .
$$
Because $(\bar{\om} \ot \io) \de(z)$ belongs to $\Nfi \cap \Nfi^*$ we have that $(\bar{\om} \ot \io)(W^*)
\la(z)$ belongs to the domain of $\nab^{1/2}$ and
$$J \nab^{1/2} (\bar{\om} \ot \io)(W^*) \la(z) = \la((\om \ot \io) \de(z^*)) = (\om \ot
\io)(W^*) \la(z^*) = (\om \ot \io)(W^*) J \nab^{1/2} \la(z) \; .$$
Because $\la(\Nfi \cap \Nfi^*)$ is a core for $\nab^{1/2}$ we get the first claim.

Using Proposition~\ref{29} rather then \ref{28} one has,
for all $\om \in B(H)_*$:
\begin{equation} \label{hulpje4}
(\io \ot \om)(W) \Jh \nabh^{1/2} \subset \Jh \nabh^{1/2} (\io \ot \bar{\om})(W).
\end{equation}
Rewriting Eq.~\eqref{hulpje3} we get,
for all $\xi,\eta \in H$, $\mu \in \cD(\nab^{1/2})$ and $\rho \in \cD(\nab^{-1/2})$:
\begin{equation} \label{hulpje5}
\langle W^* (\xi \ot J \nab^{1/2} \mu), \eta \ot \rho \rangle =\langle W (\xi \ot \nab^{1/2} J \rho), \eta \ot \mu \rangle.
\end{equation}
Rewriting Eq.~\eqref{hulpje4} we get,
for all $\mu,\rho \in H$, $\xi \in \cD(\nabh^{1/2})$ and $\eta \in \cD(\nabh^{-1/2})$:
\begin{equation} \label{hulpje6}
\langle W (\Jh \nabh^{1/2} \xi \ot \mu), \eta \ot \rho \rangle =
\langle W^* (\nabh^{1/2} \Jh \eta \ot \mu), \xi \ot \rho \rangle.
\end{equation}
Next we claim that $W(\nabh \ot \nab) = (\nabh \ot \nab) W$. To prove this, choose $\xi,\eta \in
\cD(\nabh)$ and $\mu,\rho \in \cD(\nab)$. Using the two previous equations one has
\begin{align*}
\langle W(\xi \ot \mu), \nabh \eta \ot \nab \rho \rangle &=
\langle W^*(\Jh \nabh^{1/2} \eta \ot \mu), \Jh \nabh^{1/2} \xi \ot \nab \rho \rangle \\
&= \langle W ( \Jh \nabh^{1/2} \eta \ot J \nab^{1/2} \rho), \Jh \nabh^{1/2} \xi \ot J
\nab^{1/2} \mu \rangle \\
&= \langle W^* (\nabh \xi \ot J \nab^{1/2} \rho) , \eta \ot J \nab^{1/2} \mu \rangle \\
&= \langle W (\nabh \xi \ot \nab \mu), \eta \ot \rho \rangle \; .
\end{align*}
From this chain of equalities the second claim follows. Then
$W(\nabh^{1/2} \ot \nab^{1/2}) = (\nabh^{1/2} \ot \nab^{1/2}) W$. Using this and
Eq.~\eqref{hulpje5} and \eqref{hulpje6} we get \begin{align*}
\langle W(\xi \ot \mu),\eta \ot \rho \rangle &= \langle W^* (\nabh^{1/2} \Jh \eta \ot \mu),
\nabh^{-1/2} \Jh \xi \ot \rho \rangle \\
&= \langle W^* (\nabh^{1/2} \ot \nab^{1/2}) (\Jh \eta \ot \nab^{-1/2}\mu), \nabh^{-1/2} \Jh \xi \ot \rho \rangle \\
&= \langle W^* (\Jh \eta \ot \nab^{-1/2} \mu), \Jh \xi \ot \nab^{1/2} \rho \rangle \\
&= \langle W (\Jh \eta \ot J \rho), \Jh \xi \ot J \mu \rangle \\
&= \langle (\Jh \ot J) W^* (\Jh \ot J) (\xi \ot \mu), \eta \ot \rho \rangle \; .
\end{align*}
for every
$\xi \in \cD(\nabh^{1/2})$, $\eta \in \cD(\nabh^{-1/2})$,
$\mu \in \cD(\nab^{-1/2})$ and $\rho \in \cD(\nab^{1/2})$.
This means that
$$(\Jh \ot J) W^* (\Jh \ot J) = W \; .$$
\end{proof}
\begin{lemma}\label{211}
We have
$$\{ (\io \ot \om)(W) \mid \om \in B(H)_* \}^\sluit = M \; .$$
\end{lemma}
\begin{proof}
Define the following subspaces of $M$.
\begin{align*}
\cO_1 &= \bigl( \lspan\{ (\io \ot \om) \de(z) \mid \om \in M_*, z \in M \} \bigr)^\sluit \\
\cO_2 &= \{ (\io \ot \om)(W^*) \mid \om \in B(H)_* \}^\sluit \; .
\end{align*}
Using \cite{KV1}, Section 4 we get
\begin{align*}
\cO_1 &= \bigl(\lspan\{ (\io \ot \om_{\la(a),\la(b)})\de(x) \mid a,b \in \cT_\vfi, x \in \Nfi
\}\bigr)^\sluit \\
&= \bigl(\lspan \{ (\io \ot \vfi) \bigl( (1 \ot \si_i(a) b^*) \de(x) \bigr)
\mid a,b \in \cT_\vfi, x \in \Nfi \}\bigr)^\sluit \\
&= \bigl(\lspan\{ (\io \ot \om_{\la(x), \la(b \si_{-i}(a^*))})(W^*) \mid a,b \in \cT_\vfi,
x \in \Nfi \}\bigr)^\sluit = \cO_2 \; .
\end{align*}
It is clear that $\cO_1$ is self-adjoint. Because $W$ is a multiplicative unitary we get that
$\cO_2$ is an algebra. Using the previous computation we get that $\cO_1=\cO_2:=\cO$ is a
$*$-algebra, closed in \strong topology. We have $$(\io \ot \om)\de(\al(x)) = \al((\io \ot \om
\al) \det(x)) \quad\text{for all} \; x \in M_2 \tekst{and} \om \in M_*$$ and because $\al$ is
faithful we may conclude that $\al(M_2) \subset \cO$. In particular $1 \in M$ and hence $\cO$
is a von Neumann subalgebra of $M$. Then the bicommutant theorem gives, for all $z \in M_1 \ot
M$: $$(\io \ot \deop)(z) \in M_1 \ot M \ot \cO.$$
In particular, we get that $(\io \ot
\deop)(\Wtil)$ belongs to $M_1 \ot M \ot \cO$. Proposition~\ref{25} implies $$(\Wtil \ot 1 \ot
1) \; \bigl( (\io \ot \al)\be \ot \io \ot \io \bigr)(\Wtil) \; (\io \ot \al \ot \al)(\cV) \in
M_1 \ot M \ot \cO \; .$$ Because $\al(M_2) \subset \cO$ it follows that $$\bigl( (\io \ot
\al)\be \ot \io \ot \io \bigr)(\Wtil) \in M_1 \ot M \ot \cO \; .$$ Because $(\io \ot \al)\be$
is faithful we may conclude that $(\om \ot \io \ot \io)(\Wtil)$ belongs to $\cO$ for all $\om
\in M_{1 \, *}$. Together with $\al(M_2) \subset \cO$ we get that $\cO = M$. But then also
$\cO^* = M$ which concludes the proof.
\end{proof}
Then we can finally prove the needed formula.
\begin{proposition} \label{212}
The map $R(z) = \Jh z^* \Jh\ (z \in M)$ is a $*$-anti-automorphism of $M$ satisfying
$$\de R = \si (R \ot R) \de \; .$$
\end{proposition}
\begin{proof}
$R$ is a $*$-anti-homomorphism from $M$ into $B(H)$.
Lemma~\ref{210} gives
$$R \bigl( (\io \ot \om_{\xi,\eta})(W) \bigr) = (\io \ot \om_{J\eta,J\xi})(W)$$
for all $\xi,\eta \in H$. Because $R$ is continuous in
\strong topology, it follows from
Lemma~\ref{211} that $R(M) \subset M$. But clearly $R(R(M))=M$ and so $R(M)=M$. Hence $R$ is a
$*$-anti-automorphism of $M$.

Let $(e_i)_{i \in I}$ be an orthonormal basis for $H$. Because $W$ is a multiplicative unitary
we get for all $\xi,\eta \in H$, with \strong convergence
$$\de \bigl( (\io \ot \om_{\xi,\eta})(W) \bigr) = (\io \ot \io \ot \om_{\xi,\eta})(W_{13}
W_{23}) = \sum_{i \in I} (\io \ot \om_{e_i,\eta})(W) \ot (\io \ot \om_{\xi,e_i})(W) \; . $$
So we also get for all $\xi,\eta \in H$
\begin{multline*}
\si(R \ot R) \de \bigl( (\io \ot \om_{\xi,\eta})(W) \bigr) = \sum_{i \in I} (\io \ot
\om_{Je_i,J\xi})(W) \ot (\io \ot \om_{J\eta,Je_i})(W) = (\io \ot \io \ot
\om_{J\eta,J\xi})(W_{13} W_{23}) \\ = \de \bigl( (\io \ot \om_{J\eta,J\xi})(W) \bigr)
= \de R \bigl( (\io \ot \om_{\xi,\eta})(W) \bigr) \; .
\end{multline*}
Because both $\de R$ and $\si (R \ot R) \de$ are continuous in the \strong topology, it
follows from Lemma~\ref{211} that $$\de R = \si(R \ot R) \de \; .$$
\end{proof}
We have now gathered enough results to prove the main result of this section.
\begin{theorem} \label{213}
$(M,\de)$ is a l.c.\ quantum group. The weight $\vfi$ is left invariant and $(H,\io,\la)$ is a
GNS-construction for $\vfi$. The multiplicative unitary of $(M,\de)$ in this GNS-construction
is $W$. The associated dual l.c.\ quantum group is $(\Mh,\deh)$. We call $(M,\de)$ the cocycle
bicrossed product of $(\Mo,\deo)$ and $(\Mt,\det)$.
\end{theorem}
\begin{proof}
Let $R$ be as above, then $\vfi R$ is a \nsf weight on $M$ which is right invariant. So
$(M,\de)$ is a l.c.\ quantum group, and Proposition~\ref{28} gives that $W$ is its
multiplicative unitary. It follows from \cite{KV2}, Definition 1.5 that the underlying von
Neuman algebra of the dual locally compact quantum group is $$\{ (\om \ot \io)(W) \mid \om \in
M_* \}^\sluit \; .$$ Applying Lemma~\ref{211} to the cocycle matching
$(\tautil,\cUtil,\cVtil)$ of $(\Mt,\det)$ and $(\Mo,\deo)$, we get that the dual l.c.\ quantum
group is precisely $(\Mh,\deh)$.
\end{proof}
Now we have at our disposal the l.c.\ quantum group $(M,\de)$ and the GNS-construction
$(H,\io,\la)$ for the left invariant weight $\vfi$. From \cite{KV2}, Section 1 we get in a
canonical way a left invariant weight $\vfih_0$ on $(\Mh,\deh)$ with GNS-construction
$(H,\io,\lah_0)$. In the remaining part of this subsection we prove that $\vfih_0=\vfih$ and
$\lah_0=\lah$.

Whenever $K$ is a Hilbert space and $\rho \in K$, we denote by $\te_\rho$ the rank-one operator in $B(\C,K)$ given by $\te_\rho(\lambda)=\lambda \rho$ for all $\lambda \in \C$.

\begin{lemma} \label{214}
For all $\xi \in H_1 \ot H_2$, $\rho \in H_2$, $b \in \cT_{\vfi_1}$ and $x \in \cN_{\vfi_2}$
the element $z \in M$ defined by
$$z:= (\om_{\xi,U_\beta (\la_1(b) \ot \rho)} \ot \io \ot \io)(\Wh) \; \al(x),$$
where $U_\beta = \Jh (J_1 \ot \Jh_2)$ is the canonical implementation of $\be$ as in
Proposition~\ref{117}, belongs to $\Nfi$ and
$$\la(z) = \bigl( (J_1 \si^1_{i/2}(b) J_1 \ot \te_\rho^*)U_\be^* \ot 1 \bigr) \cV (1 \ot \Wh_2)
(\xi \ot \la_2(x)) \; .$$
Here $(\si^1_t)$ denotes the modular group of the weight $\vfi_1$ on $\Mo$.
\end{lemma}
\begin{proof}
Let $(e_i)_{i \in I}$ be an orthonormal basis for $H_1 \ot H_2$. Then we have, with \strong
convergence
$$z= \sum_{i \in I} ( \om_{e_i,U_\beta (\la_1(b) \ot \rho)} \be \ot \io \ot \io)(\Wtil) \; \al
\bigl( (\om_{\xi,e_i} \ot \io)(\cV (1 \ot \Wh_2)) \; x \bigr) \; .$$
Denote by $z_{I_0}$ the same expression summed over a finite subset $I_0 \subset I$. We claim
that for all $\eta \in H_1 \ot H_2$ the normal functional $\mu:=\om_{\eta,U_\beta (\la_1(b) \ot \rho)}
\be$ belongs to $\cI_1$ and
$$\xi_1(\mu)= (J_1 \si^1_{i/2}(b)J_1 \ot \te_\rho^*)U_\be^* \eta \; .$$
To prove this claim, we choose $y \in \cN_{\vfi_1}$. Then
$$\mu(y^*) = \langle \be(y^*) \eta, U_\be(\la_1(b) \ot \xi) \rangle = \langle (y^* \ot 1)
U_\be^* \eta, \la_1(b) \ot \rho \rangle = \langle (J_1 \si^1_{i/2}(b)J_1 \ot \te_\rho^*)U_\be^*
\eta, \la_1(y) \rangle \; .$$
But then for all finite subsets $I_0 \subset I$ we have $z_{I_0} \in \Nfi$ and
$$\la(z_{I_0}) = \sum_{i \in I_0} (J_1 \si^1_{i/2}(b)J_1 \ot \te_\rho^*)U_\be^* e_i \ot
(\om_{\xi,e_i} \ot \io)(\cV (1 \ot \Wh_2)) \la_2(x) \; .$$
Because $\la$ is \strong -- norm closed the claim of the lemma follows.
\end{proof}

Before we can prove the announced result, we need a general lemma on l.c.\ quantum groups.

\begin{lemma} \label{215}
Let $(M,\de)$ be a l.c.\ quantum group with left invariant weight $\vfi$ with GNS-construction
$(H,\io,\la)$, $(\Mh,\deh)$ the associated dual l.c.\ quantum group with canonical left
invariant weight $\vfih$ with GNS-construction $(H,\io,\lah)$ and $(\sih_t)$ the modular group
of the weight $\vfih$ on $\Mh$.

Suppose that $y \in \cD(\sih_{i/2})$ and that there exist elements $(p_i)_{i \in I}$ and
$(q_i)_{i \in I}$ in $\Mh$ such that $\sum p_i p_i^*$ and $\sum q_i^* q_i$ are bounded and $$1
\ot y = \sum_{i \in I} \deh(p_i) (q_i \ot 1) \; .$$ Then there exists a unique linear map $P :
M \recht M$  continuous in the \strong topology and such that for every $x \in \Nfi$ the
element $P(x)$ belongs to $\Nfi$ and $$\la(P(x)) = \Jh \sih_{i/2}(y) \Jh \la(x) \; .$$
Moreover, whenever $\om \in M_*$ and $(\om \ot \io)(W)y \in \Nfih$ we have, for all $x \in
\Nfi$: $$\langle \lah \bigl( (\om \ot \io)(W) y), \la(x) \rangle = \om \bigl( P(x)^* \bigr).$$
\end{lemma}
\begin{proof}
Denote, for all $i \in I, a_i = \Rh(q_i)$ and $b_i = \Rh(p_i)$. Then $\sum a_i a_i^*$ and
$\sum b_i^* b_i$ are bounded and $$1 \ot \Rh(y) = \sum_{i \in I} (a_i \ot 1) \dehop(b_i) \;
.$$ So we can define for every $x \in M$ an operator $P(x) \in B(H)$ as $$P(x) = \sum_{i \in
I} \Jh a_i \Jh x \Jh b_i \Jh \; .$$ Considering $A=(\Jh a_i \Jh)_{i \in I}$ as an operator in
$B(H \ot \ell^2(I),H)$ and $B = (\Jh b_i \Jh)_{i \in I}$ as an operator in $B(H,H \ot
\ell^2(I))$, we get $P(x) = A(x \ot 1)B$ for all $x \in M$, so $P$ is continuous in the
\strong topology. We claim that $$P\bigl( (\io \ot \mu)(W^*) \bigr) = (\io \ot \mu y^*)(W^*)
\;$$ for all $\mu \in \Mh_*$. Indeed, using the formula $(\Jh \ot J)W = W^* (\Jh \ot J)$,
proved in \cite{KV2}, Corollary 2.2, we have $$\sum_{i \in I} (\Jh a_i \Jh \ot 1) W^* (\Jh b_i
\Jh \ot 1) = \sum_{i \in I} (\Jh \ot J)(a_i \ot 1) \dehop(b_i) (\Jh \ot J)W^* = (\Jh \ot J)(1
\ot \Rh(y)) (\Jh \ot J) W^* = (1 \ot y^*)W^*,$$ from where the claim follows. For all $\mu \in
\hat{\cI}$ we have $\mu y^* \in \hat{\cI}$, hence $P \bigl( (\io \ot \mu)(W^*) \bigr) \in
\Nfi$ and $$\la \bigl( P \bigl( (\io \ot \mu)(W^*) \bigr) \bigr) = \la \bigl( (\io \ot \mu
y^*)(W^*) \bigr) = \Jh \sih_{i/2}(y) \Jh \la \bigl((\io \ot \mu)(W^*) \bigr) \; .$$ If now $x
\in \Nfi$, there is a net $(\mu_\al)$ in $\hat{\cI}$ such that $(\io \ot \mu_\al)(W^*) \recht
x$ in the \strong topology and $\la\bigl((\io \ot \mu_\al)(W^*) \bigr) \recht \la(x)$ in norm.
Because $P$ is continuous in \strong topology and $\la$ is \strong -- norm closed, we get
$P(x) \in \Nfi$ and $$\la \bigl( P(x) \bigr) = \Jh \sih_{i/2}(y) \Jh \la(x) \; .$$ This
concludes the first part of the proof. Now suppose that we also have $\om \in M_*$ such that
$(\om \ot \io)(W) y \in \cN_{\vfih}$. Then we get, for all $\mu \in \hat{\cI}$: $$\langle \lah
\bigl( (\om \ot \io)(W) y \bigr) , \la \bigl( (\io \ot \mu)(W^*) \bigr) \rangle = \bar{\mu}
\bigl( (\om \ot \io)(W) y \bigr) = \om \bigl( P \bigl( (\io \ot \mu)(W^*) \bigr)^*\bigr) \;
.$$ If $x \in \Nfi$ we take the same net $(\mu_\al)$ in $\hat{\cI}$ as above and obtain that
$$\langle \lah \bigl( (\om \ot \io)(W) y), \la(x) \rangle = \om \bigl( P(x)^* \bigr) \; .$$
\end{proof}
Now we are ready to prove the following nice result. \begin{proposition} \label{JJhat}
We have $\vfih_0 = \vfih$ and $\lah_0 = \lah$. Hence the left invariant
weight with GNS-construction obtained by considering $(\Mh,\deh)$ as a cocycle bicrossed
product, as in Proposition~\ref{29}, is the same as the left invariant weight with
GNS-construction obtained by considering $(\Mh,\deh)$ as the dual of $(M,\de)$, as in \cite{KV1}, Proposition 8.14.
\end{proposition}
\begin{proof}
Consider the l.c.\ quantum group $(\Mh,\deh)$ with left invariant weight $\vfih$,
GNS-construction $(H,\io,\lah)$ for $\vfih$ and multiplicative unitary $\Wh$. Applying
\cite{KV1}, Propositions~8.14 and 8.15 to the l.c.\ quantum group $(\Mh,\deh)$ we obtain a
canonical left invariant weight $\vfihh$ on the dual $(M,\de)$ of $(\Mh,\deh)$, with canonical
GNS-construction $(H,\io,\lahh)$. In view of \cite{KV1}, Proposition 8.30 it suffices to prove
that $\vfihh=\vfi$ and $\lahh=\la$.

By the uniqueness of left invariant weights, there exists a positive number $\nu >0$ such that
$\vfihh = \nu \vfi$. So $\cN_{\klvfihh} = \Nfi$. Then we can define a unitary operator $u$ on $H$
such that $u \lahh(x) = \nu^{1/2} \la(x)$ for all $x \in \Nfi$. Because both
GNS-representations involved are the identical representations, we get $u \in
M'$. From the dual version of \cite{KV1}, Proposition 8.16 we get, for all $x \in \cN_{\klvfihh}$ and $\om \in M_*$:
$$\lahh \bigl( (\om \ot \io)\de(x) \bigr) = (\om \ot \io)(W^*) \lahh(x).$$
But from Proposition 2.8 it follows, for all $x \in \Nfi$ and $\om \in M_*$:
$$\la \bigl( (\om \ot \io)\de(x) \bigr) = (\om \ot \io)(W^*) \la(x).$$
Both formulas together imply $u \in
\Mh'$. Because we already have $u \in M'$, it follows that $u \in \C 1$. Hence we can take
$\lambda \in \C \setminus \{0\}$ such that $\lahh = \lambda \la$. We have to
prove that $\lambda=1$. We can already conclude that the modular group $(\sihh_t)$ of $\vfihh$
equals the modular group $(\si_t)$ of $\vfi$ and that $\Jhh = \lambda / \bar{\lambda} \; J$.

Choose now $\xi \in H_1 \ot H_2$, $\rho \in H_2, b \in \cT_{\vfi_1}$ and $y \in
\cN_{\vfi_2} \cap \cD(\si^2_{i/2})$ such that there exist elements $(p_i)_{i \in I}$ and
$(q_i)_{i \in I}$ in $\Mt$ with $\sum p_i p_i^*$ and $\sum q_i^* q_i$ bounded and
\begin{equation} \label{eq.y}
1 \ot y = \sum_{i \in I} \det(p_i)(q_i \ot 1) \; .
\end{equation}

Remark that such elements $y$ form a dense subspace of $M_2$ in the \strong topology. Because
$\psi_2 \si^2_t = \nu_2^t \, \psi_2$ and $\psi_2 \tau^2_t = \nu_2^{-t} \psi_2$ we can take
elements $a,b \in M_2$ such that $a$ is analytic w.r.t.\ $(\si^2_t)$, $b$ is analytic w.r.t.\
$(\tau^2_t)$ and $\si^2_z(a) \in \cN_{\psi_2}^* \cN_{\vfi_2}$, $\tau^2_z(b) \in \cN_{\psi_2}$
for all $z \in \C$.  Now \cite{KV1}, Result~2.6 gives that $$y:= (\psi_2 \ot \io)(\de(a)(b \ot
1)) \in \cN_{\vfi_2}.$$ Using the relation $(\tau^2_t \ot \si^2_t) \de_2 = \de_2 \si^2_t$ we
get that $y$ is analytic w.r.t.\ $(\si^2_t)$ and $$\si^2_z(y) = \nu_2^z \, (\psi_2 \ot
\io)(\de(\si^2_z(a)) (\tau^2_z(b) \ot 1)) \; .$$ Finally observe that $y =
(\om_{\Gamma_2(b),\Gamma_2(a^*)} \ot \io)(V_2^*)$, so defining $p_i = (\om_{e_i,\Gamma_2(a^*)}
\ot \io)(V_2^*)$ and $q_i = (\om_{\Gamma_2(b),e_i} \ot \io)(V_2)$, where $(e_i)$ is an
orthonormal basis for $H$, we have \eqref{eq.y}.

Because $y \in \cN_{\vfi_2}$, it follows from Lemma~\ref{214} that the element $z \in M$ defined
by
$$z:= (\om \ot \io \ot \io)(\Wh) \; \al(y) \; , \tekst{with} \om = \om_{\xi,U_\be (\la_1(b) \ot \rho)}$$
belongs to $\Nfi$ and
$$\la(z) = \bigl( (J_1 \si^1_{i/2}(b) J_1 \ot \te_\rho^*) U_\be^* \ot 1 \bigr) \cV (1 \ot
\Wh_2) (\xi \ot \la_2(y)) \; .$$
Choose now $\mu \in \cI_2$ and $x \in \cN_{\vfi_1}$. Then
\begin{align}
\langle \la(z) &, \lah \bigl( (\io \ot \io \ot \mu)((1 \ot \Wh_2^*)\cV^* ) \; \be(x) \bigr)
\rangle \notag \\ &= \langle  \bigl( (J_1  \si^1_{i/2}(b) J_1 \ot \te_\rho^*) U_\be^* \ot 1 \bigr) \cV (1 \ot
\Wh_2) (\xi \ot \la_2(y)), \la_1(x) \ot \xi_2(\mu) \rangle \notag \\
&=\langle (\om \be(x^*) \ot \io)(\cV (1 \ot \Wh_2)) \la_2(y), \xi_2(\mu) \rangle
= \bar{\mu} \bigl( (\om \ot \io)((\be(x^*) \ot 1)\cV(1 \ot \Wh_2)) y \bigr) \notag \\
&= \om \bigl( \bigl( (\io \ot \io \ot \mu y^*)((1 \ot \Wh_2^*)\cV^*) \be(x) \bigr)^* \bigr) \; .
\label{vgltje}
\end{align}
Now apply Lemma~\ref{215} to the l.c.\ quantum group $(\Mh,\deh)$. Because $\sihh_t \al =
\si_t \al = \al \si^2_t$ for all $t \in \R$, and because $\de \al = (\al \ot \al) \det$ we get
that $\al(y)$ satisfies the conditions of this lemma. Hence there exists a \strong continuous
map $P$ from $\Mh$ to $\Mh$, such that for all $v \in \cN_{\vfih}$ we get $P(v) \in
\cN_{\vfih}$ and $$\lah(P(v)) = \Jhh \sihh_{i/2} (\al(y)) \Jhh \lah(v) \; .$$ Because $\Jhh =
\lambda / \bar{\lambda} J, \sihh_t \al = \al \si^2_t$ and $\al$ is implemented by $J(\Jh_1 \ot
J_2)$ we get, for all $v \in \cN_{\vfih}$: $$\lah(P(v)) = (1 \ot J_2 \si^2_{i/2}(y) J_2)
\lah(v).$$ Since $z = (\om \ot \io \ot \io)(\Wh) \al(y)\in \Nfi= \cN_{\klvfihh}$, it follows
from Lemma~\ref{215} that
\begin{equation}\label{vgltje2}
\langle \lahh(z),\lah(v) \rangle = \om(P(v)^*)
\quad\text{for all} \; v \in \cN_{\vfih}.
\end{equation}
If we now take $v = (\io \ot \io \ot \mu)(( 1 \ot \Wh_2^*)\cV^*) \; \be(x)$, then
$$(1 \ot J_2 \si^2_{i/2}(y) J_2) \lah(v) = \la_1(x) \ot J_2 \si^2_{i/2}(y) J_2 \xi_2(\mu) =
\la_1(x) \ot \xi_2(\mu y^*) = \lah \bigl( (\io \ot \io \ot \mu y^*)((1 \ot \Wh_2^*)\cV^*) \;
\be(x) \bigr) \; .$$
So $P(v) = (\io \ot \io \ot \mu y^*)((1 \ot \Wh_2^*)\cV^*) \; \be(x)$. Now
combine Eq.~\eqref{vgltje} and Eq.~\eqref{vgltje2} to get
$$\om(P(v)^*) = \langle \lahh(z),\lah(v) \rangle = \lambda \langle \la(z),\lah(v) \rangle =
\lambda \; \om(P(v)^*) \; .$$
Because this is valid for enough elements $v$ and $\om$ we get that $\lambda=1$ which concludes the proof.
\end{proof}
We characterize compact and discrete cocycle bicrossed products as follows.
\begin{proposition} \label{compact}
The l.c.\ quantum group $(M,\de)$ is compact if and only if $(M_2,\de_2)$ is compact and
$(M_1,\de_1)$ is discrete. $(M,\de)$ is discrete if and only if $(M_2,\de_2)$ is discrete and
$(M_1,\de_1)$ is compact.
\end{proposition}
\begin{proof}
From Definition~\ref{112} it follows that $\vfi(1) = \vfih_1(1) \; \vfi_2(1)$ from where the
first statement follows. The second statement follows from the first one and from the observation that $(\Mh,\deh)$ is a cocycle bicrossed product of $(M_2,\de_2)$ and $(M_1,\de_1)$.
\end{proof}

\section{Extensions and cocycle bicrossed products}
We define and study extensions of l.c.\ quantum groups and prove that there is a one-to-one
correspondence between cleft extensions and cocycle bicrossed products.
\subsection{Extensions and cleft extensions of l.c.\ quantum groups}
Let us start with the following result which can be deduced from \cite{Kus}, but we prefer to
give a direct proof.
\begin{proposition}\label{31}
Let $(\Mo,\deo)$ and $(M,\de)$ be l.c.\ quantum groups and $\be : M_1 \recht \Mh$ a normal
$*$-homo\-mor\-phism satisfying $\deh \be = (\be \ot \be) \deo$. Define the unitaries $Z_1$
and $Z_2$ on $H_1 \ot H$ by $$Z_1 = (\io \ot \be)(\Wh_1^*) \qquad\text{and}\qquad Z_2= (J_1
\ot \Jh) \; (\io \ot \be)(\Wh^*_1) \; (J_1 \ot \Jh)$$ and then the maps $\mu$ and $\te$ on $M$
by $$\mu(z) = Z_1(1 \ot z) Z_1^* \qquad\text{and}\qquad \te(z) = Z_2(1 \ot z) Z_2^*.$$ Then
$\mu$ is an action of $(\Mh_1,\deoh)$ on $M$ and $\te$ is an action of $(\Mh_1,\dehopo)$ on
$M$ such that
\begin{align*}
\te(z) &= (\hat{R}_1 \ot R) \mu (R(z)) \tekst{for all} z \in M \\
(\mu \ot \io)(W) &= (\io \ot \be)(\Wh_1)_{13} W_{23} \\
(\te \ot \io)(W) &= W_{23} (\io \ot \be)(\Wh_1)_{13} \\
(\io \ot \de)\mu = (\mu \ot \io) \de \quad &\text{and}\quad (\io \ot \deop)\te = (\te \ot
\io)\deop \; .
\end{align*}
Finally, for all $z \in \Nfi$ and $\om \in \Mh_{1 \, *}$ we have $(\om \ot \io)\mu(z) \in \Nfi$
and
$$\la \bigl( (\om \ot \io)\mu(z) \bigr) = \be \bigl( (\io \ot \om)(W_1) \bigr) \; \la(z) \; .$$
\end{proposition}
\begin{proof}
Observe that
\begin{align*}
(\io \ot \be)(\Wh_1^*)_{12} \; W_{23} \; (\io \ot \be)(\Wh_1)_{12} &= (\io \ot
\be)(\Wh_1^*)_{12} \; (\io \ot \dehop \be)(\Wh_1) \; W_{23} \\
&=(\io \ot \be)(\Wh_1^*)_{12} \; (\io \ot \be \ot \be)(\io \ot \deopo)(\Wh_1) \; W_{23} \\
&=(\io \ot \be)(\Wh_1)_{13} \; W_{23} \; .
\end{align*}
Then one can conclude that $\mu(M) \subset \Moh \ot M$ and
$$(\mu \ot \io)(W) = (\io \ot \be)(\Wh_1)_{13} \; W_{23} \; .$$
Now one can check that $(\deoh \ot \io)\mu = (\io \ot
\mu)\mu$, i.e., $\mu$ is an action of $(\Moh,\deoh)$ on $M$,  and $(\io \ot \de) \mu = (\mu \ot \io)\de$. We get similar formulas for $\te$ from
$$\te(z)=(\hat{R}_1 \ot R) \mu(R(z)) \tekst{for all} z \in M$$ and $\hat{R} \be = \be R_1$ (this follows from the von Neumann algebraic version of \cite{KV1}, Proposition 5.45).
So it remains to prove the final statement. Using \cite{KV1},
choose $\om \in \Mh_{1 \, *}$ and $\eta \in \hat{\cI}$. Then
\begin{equation} \label{vgl1}
(\om \ot \io)\mu \bigl( (\io \ot \eta)(W^*) \bigr) = (\io \ot \rho)(W^*) \tekst{where}
\rho(z) = (\om \ot \eta) \bigl( (1 \ot z) (\io \ot \be)(\Wh_1^*) \bigr) \tekst{for all} z \in
\Mh \; .
\end{equation}
We claim that $\rho \in \hat{\cI}$. For this we choose $z \in \cN_{\vfih}$ and compute
\begin{multline*}
\rho(z^*) = (\om \ot \eta) \bigl( (1 \ot z^*) (\io \ot \be)(\Wh_1^*) \bigr) = \eta \bigl( z^*
\be( (\io \ot \om)(W_1) ) \bigr) \\ = \langle \la\bigl( (\io \ot \eta)(W^*) \bigr) , \be \bigl(
(\io \ot \om)(W_1)^* \bigr) \lah(z) \rangle  = \langle \be \bigl( (\io \ot \om)(W_1) \bigr)
 \la\bigl( (\io \ot \eta)(W^*) \bigr) , \lah(z) \rangle \; .
\end{multline*}
From this the claim follows, so $(\io \ot \rho)(W^*) \in \Nfi$ and $$\la \bigl( (\io \ot
\rho)(W^*) \bigr) = \be \bigl( (\io \ot \om)(W_1) \bigr) \; \la\bigl( (\io \ot \eta)(W^*)
\bigr) \; .$$ Combining this with Eq.~\eqref{vgl1} and the fact that the elements $(\io \ot
\eta)(W^*)$ form a \strong -- norm core for $\la$, we get the final statement.
\end{proof}
The starting point of Proposition~\ref{31} is the von Neumann algebraic morphism $\be$. If all algebras are finite dimensional, we can dualize $\be$ to a morphism
$\tilde{\be} : M \recht \Moh$ such that
$$(\tilde{\be} \ot \io)(W) = (\io \ot \be)(\hat{W}_1) \; .$$
The link between $\tilde{\be}$ and the action $\mu$ is now given by the formula
$$\mu(z) = (\tilde{\be} \ot \io) \de(z) \tekst{for all} z \in M \; ,$$
and the action $\mu$ also makes sense in general case. The action $\mu$ satisfies $(\io \ot \de) \mu = (\mu \ot \io) \de$, and
in \cite{Kus} there is established a natural one-to-one correspondence between actions of this special type
and morphisms.

Turning now to the definition of extensions of l.c.\ quantum groups, let us note that we use
the notation $\be$ first of all for a morphism $\Mo\recht\Mh$, but also symbolically on the
arrow pointing from $(M,\de)$ to $(\Mh_1,\deoh)$ for, in a sense, the dual morphism which does
not necessarily exist on the von Neumann algebra level.
\begin{definition} \label{32}
Let $(\Mo,\deo)$, $(\Mt,\det)$ and $(M,\de)$ be l.c.\ quantum groups. We call $$(\Mt,\det)
\overset{\al}{\lrecht} (M,\de) \overset{\be}{\lrecht} (\Moh,\deoh)$$ a short exact sequence,
if $$\al: \Mt \recht M \quad\text{and}\quad \be: M_1 \recht \Mh$$ are normal, faithful
$*$-homomorphisms satisfying $$\de \al = (\al \ot \al) \det \quad\text{and}\quad \deh \be =
(\be \ot \be) \deo$$ and if $\al(M_2)  = M^\te$, where $\te$ is the action of $(\Moh,\dehopo)$
on $M$ defined in Proposition~\ref{31} starting with the morphism $\be$. In this situation we
call $(M,\de)$ an extension of $(\Mt,\det)$ by $(\Moh,\deoh)$.
\end{definition}
The faithfulness of the morphisms $\al$ and $\be$ reflects the exactness of the sequence in
the first and third place. The formula $\al(M_2)=M^\te$ reflects its exactness in the second
place. In terms of the above morphism $\tilde{\be} : M \recht \Moh$ the exactness of the
sequence in the second place is formally expressed by $$\al(M_2) = \{ z \in M \mid (\io \ot
\tilde{\be})\de(z) = z \ot 1 \} \; .$$ This formula is precisely the one used in Hopf algebra
theory, see \cite{A-D}.

If $\al$ and $\be$ give rise to a short exact sequence as in the previous definition, then the
dual sequence $$(\Mo,\deo) \overset{\be}{\lrecht} (\Mh,\deh) \overset{\al}{\lrecht}
(\Mth,\deth)$$ is also a short exact sequence. This follows from the argument in part 1.a) of
the proof of Theorem~\ref{35}.

\begin{definition} \label{33}
An extension $(M,\de)$ of l.c.\ quantum groups is said to be a cleft extension of $(\Mt,\det)$
by $(\Mh_1,\deoh)$ if the extension $\al(M_2) = M^\te \subset M$ is a cleft extension of von
Neumann algebras.
\end{definition}

From part 1.b) of the proof of Theorem~\ref{35} it follows that the dual extension of a cleft extension is again cleft.

\subsection{From cocycle bicrossed products to cleft extensions}
\begin{proposition} \label{34}
Let $(M,\de)$ be the cocycle bicrossed product of l.c.\ quantum groups $(\Mo,\deo)$ and
$(\Mt,\det)$, which is also a l.c.\ quantum group by Theorem~\ref{213}. Let $\al$ and $\be$ be
as in Definition~\ref{21}. Then $$\al: \Mt \recht M \qquad\text{and}\qquad \be:\Mo \recht
\Mh$$ and $(M,\de)$ is a cleft extension of $(\Mt,\det)$ by $(\Moh,\deoh)$.
\end{proposition}
\begin{proof}
Theorem~\ref{213} shows that $(\Mh,\deh)$, defined in Proposition~\ref{29}, is the dual l.c.\
quantum group of $(M,\de)$, so $\be$ maps $\Mo$ faithfully into $\Mh$. Definition~\ref{22}
implies that $\al$ maps $M_2$ faithfully into $M$ and Proposition~\ref{24} shows that $\de \al
= (\al \ot \al)\det$. Because the construction of $(\Mh,\deh)$ corresponds to the cocycle
matching $(\tilde{\tau},\cUtil,\cVtil)$, Proposition~\ref{24} also shows that $\deh \be = (\be
\ot\be) \deo$. Let $J$ and $\Jh$ be the modular conjugations of $\vfi$ and $\vfih$ in the
GNS-constructions $(H,\io,\la)$ and $(H,\io,\lah)$ respectively (see Definition~\ref{27} and
Proposition~\ref{29}). From Proposition~\ref{JJhat} we know that $\vfih$ is the canonical left
invariant weight on the dual $(\Mh,\deh)$ with canonical dual GNS-construction $(H,\io,\lah)$.
So, the action $\te$ of $(\Moh,\dehopo)$ on $M$ constructed in Proposition~\ref{31} satisfies
$\te(z)=Z(1 \ot z) Z^*$ for all $z \in M$ where $$Z = (J_1 \ot \Jh) \; (\io \ot \be)(\Wh_1^*)
\; (J_1 \ot \Jh) \; .$$ Up to a flip map, $\Mh$ is the cocycle crossed product $M_2
\kruisje{\si\be,\cV_{321}} M_1$ and $\hat\vfi$ is the dual weight of $\vfi_1$ with
accompanying dual weight GNS-construction $(H,\io,\lah)$. Hence it follows from
Proposition~\ref{117} that $\be$ is implemented by $\Jh (J_1 \ot \Jh_2)$. Then we get $$Z =
(J_1 \ot J_1) \Si W_1 \Si (J_1 \ot J_1) \ot 1 \; .$$

So it follows from Proposition~\ref{14} that $\te = \alh$, therefore Theorem~\ref{110} gives
$$M^\te = (M_1 \kruisje{\al,\cU} M_2)^{\alh} = \al(M_2) \;
.$$
Since $\te=\alh$, Propositions~\ref{120} and \ref{122} show that $M^\te \subset M$ is a cleft extension.
\end{proof}
\subsection{From cleft extensions to cocycle bicrossed products}
\begin{theorem} \label{35}
Let $(\Mo,\deo)$, $(\Mt,\det)$ and $(N,\Ga)$ be l.c.\ quantum groups, $\rho: \Mt \recht N$ and
$\ga: \Mo \recht \Nh$ normal and faithful $*$-homomorphisms turning $(N,\Ga)$ into a cleft
extension of $(\Mt,\det)$ by $(\Moh,\deoh)$.

Then there exists a cocycle matching $(\tau,\cU,\cV)$ of $(\Mo,\deo)$ and $(\Mt,\det)$ and an
isomorphism $$\Phi:(N,\Ga) \recht (M,\de)$$ of $(N,\Ga)$ onto the cocycle bicrossed product
$(M,\de)$, satisfying $$\Phi \rho = \al \qquad\text{and}\qquad \hat{\Phi} \ga = \be$$ where
$\al$ and $\be$ are as in Definition~\ref{21} and where $\hat{\Phi}:(\Nh,\Gah)\recht
(\Mh,\deh)$ is the isomorphism canonically associated with $\Phi$ and characterized by $$(\Phi
\ot \hat{\Phi})(\cW) = W \; ,$$ where $\cW$ is the multiplicative unitary of $(N,\Ga)$ and $W$
is the one of $(M,\de)$.
\end{theorem}
\begin{proof} {\bf 1. The definition of $\al,\cU,\be,\cV$ and $\Phi$.}

a) Represent $(\Mo,\deo)$, $(\Mt,\det)$ and $(N,\Ga)$ on the
GNS-space of a left Haar weight and denote by $K$ the Hilbert space on which $(N,\Ga)$ acts, by
$\cW$ the multiplicative unitary of $(N,\Ga)$ and by $\Jn$ and $\Jhn$ the modular
conjugations of the invariant weights on $N$ and $\Nh$.

Let $\te: N \recht \Moh \ot N$ be the action of $(\Moh,\dehopo)$ on $N$ constructed with
$\ga$, as in Proposition~\ref{31}, then $$\te(z) = Z(1 \ot z)Z^* \tekst{for all} z \in N
\qquad\text{and}\qquad Z = (J_1 \ot \Jhn)(\io \ot \ga)(\Wh^*_1) (J_1 \ot \Jhn) \; .$$ One also
has $$(\io \ot \te)(\cWh) = (\ga \ot \io)(W_1)_{12} \cWh_{13} \; .$$ Applying now
Proposition~\ref{31} to $\rho:M_2 \recht N$ we get an action $\eta : \Nh \recht \Mth \ot \Nh$
of $(\Mth,\dehopt)$ on $\Nh$ by
\begin{equation} \label{eq1}
\eta(z) = \tilde{Z} (1 \ot z) \tilde{Z}^* \tekst{for all} z \in \Nh \qquad\text{and}\qquad
\tilde{Z} = (J_2 \ot \Jn)(\io \ot \rho)(\Wh_2^*) (J_2 \ot \Jn) \; .
\end{equation}
The action $\eta$ satisfies
$$(\eta \ot \io)(\cWh) = \cWh_{23} (\io \ot \rho)(\Wh_2)_{13} \; .$$
Because $\ga$ and $\rho$ turn $(N,\Ga)$ into an extension, one has $N^\te =
\rho(M_2)$. We claim that also $\Nh^\eta = \ga(M_1)$. Indeed, \eqref{eq1} gives
$$\Nh^\eta = \Nh \cap \Jn \rho(M_2)' \Jn \; .$$
Similarly we have
$$N^\te = N \cap \Jhn \ga(M_1)' \Jhn \; .$$
So
$$\rho(M_2) = N \cap \Jhn \ga(M_1)' \Jhn \; .$$
But for all $x \in M_2$ we have $\Jhn \rho(x) \Jhn = R(\rho(x)^*) = \rho(R_2(x^*))$ and so
$\Jhn \rho(M_2) \Jhn = \rho(M_2)$. Therefore,
$$\rho(M_2) = N \cap \ga(M_1)' \; .$$
Since $\Gah \ga = (\ga \ot \ga) \deo$ we get that $\ga(M_1)$ is a two-sided coideal in
$(\Nh,\Gah)$ in the sense of \cite{Eno}, D\'efinition 3.1. Then it follows from \cite{Eno}, Th\'eor\`eme 3.3 that
$$\ga(M_1) = \Nh \cap (N \cap \ga(M_1)')' = \Nh \cap \rho(M_2)' \; .$$
Because $\Jn \ga(\Mo) \Jn = \ga(\Mo)$ one has
$$\ga(\Mo) = \Nh \cap \Jn \rho(\Mt)' \Jn = \Nh^\eta$$
which proves the claim.

Since $\ga$ and $\rho$ turn $(N,\Ga)$ into a cleft extension, there exists a unitary $X \in
M_1 \ot N$ satisfying $$(\io \ot \te)(X) = W_{1,12} X_{13} \; .$$ From Proposition~\ref{120}
it follows that $(\tilde{\al},\cUtil)$ is a cocycle action of $(\Mo,\deo)$ on $N^\te$, where
$$\tilde{\al}(z) = X^*(1 \ot z) X \tekst{for all} z \in N^\te \qquad\text{and}\qquad \cUtil =
X^*_{23}X^*_{13}(\deo \ot \io)(X).$$ Because $N^\te=\rho(\Mt)$, we can define $$\al = (\io \ot
\rho^{-1})\tilde{\al} \rho \qquad\text{and}\qquad \cU = (\io \ot \io \ot \rho^{-1})(\cU)$$ and
then $\al : \Mt \recht \Mo \ot \Mt$ is a faithful and normal $*$-homomorphism, $\cU \in \Mo
\ot \Mo \ot \Mt$ is a unitary and
\begin{align*}
(\io \ot \al)\al(x) &= \cU \; (\deo \ot \io)\al(x) \; \cU^* \tekst{for all} x \in \Mt \\
(\io \ot \io \ot \al)(\cU) \; (\deo \ot \io \ot \io)(\cU) &= (1 \ot \cU)
\; (\io \ot \deo \ot \io) (\cU) \; .
\end{align*}
Still using Proposition~\ref{120} we can define a $*$-isomorphism $$\Phi : N \recht \Mo
\kruisje{\al,\cU} \Mt$$ satisfying $$\Phi(\rho(x)) = \al(x) \tekst{for all} x \in \Mt
\qquad\text{and}\qquad (\io \ot \Phi)(X) = W_{1,12} \cU^* \; .$$ From now on, we denote $M=\Mo
\kruisje{\al,\cU} \Mt$.

b) To show that also $\Nh$ has a cocycle crossed product
structure we will find a cocycle action of $(\Mt,\det)$ on $\Mo$. Let us define a unitary $\Ytil$ in $\Nh \ot N$ by
\begin{equation} \label{eqY}
\Ytil = (\ga \ot \io)(X^*) \cWh \; .
\end{equation}
Then we have
$$(\io \ot \te)(\Ytil) = (\ga \ot \te)(X^*) (\io \ot \te)(\cWh) = (\ga \ot \io)(X^*)_{13} \;
(\ga \ot \io)(W_1^*)_{12} \; (\ga \ot \io)(W_1)_{12} \; \cWh_{13} = \Ytil_{13} \; .$$
So $\Ytil \in \Nh \ot N^\te$ and we define $Y = (\io \ot \rho^{-1})(\Ytil)\in \Nh \ot \Mt$. Observe now that
$$(\eta \ot \io)(\Ytil) = (\ga \ot \io)(X^*)_{23} \; (\eta \ot \io)(\cWh) = (\ga \ot
\io)(X^*)_{23} \; \cWh_{23}  \; (\io \ot \rho)(\Wh_2)_{13} = \Ytil_{23} (\io \ot
\rho)(\Wh_2)_{13} \; .$$
From this we may conclude that
$$(\io \ot \eta)(\Si Y^* \Si) = W_{2,12} (\Si Y^* \Si)_{13} \; .$$
Hence Proposition~\ref{120} shows that $\Nh$ has a cocycle crossed product structure and we can define
\begin{align*}
& \betil: \Nh^\eta \recht \Nh^\eta \ot \Mt : \betil(z) = Y(z \ot 1) Y^* \\
& \cVtil = Y_{12} Y_{13} (\io \ot \deopt)(Y^*) \in \Nh^\eta \ot \Mt \ot \Mt \; .
\end{align*}
Then $\be = (\ga^{-1} \ot \io) \tilde{\be} \ga : \Mo \recht \Mo \ot \Mt$ is a faithful and normal $*$-homomorphism, $\cV = (\ga^{-1} \ot \io \ot \io)(\tilde{\cV})$ is a
unitary in $\Mo \ot \Mt \ot \Mt$ and
\begin{align*}
(\be \ot \io)\be(x) &= \cV (\io \ot \deopt)\be(x) \cV^* \tekst{for all} x \in \Mo \\
(\be \ot \io \ot \io)(\cV)
(\io \ot \io \ot \deopt)(\cV) &= (\cV \ot 1) (\io \ot \deopt \ot \io)(\cV) \; .
\end{align*}
If $\Mtil$ is the von Neumann algebra acting on $H_1 \ot H_2$ generated by
$\be(M_1)$ and the elements $(\io \ot \io \ot \om)(\cV (1 \ot \Wh_2))$, it also follows from
Proposition~\ref{120} that there is a unique $*$-isomorphism $\tilde{\Phi} : \Nh \recht \Mtil$
satisfying
$$\Phitil(\ga(x)) = \be(x) \tekst{for all} x \in \Mo \qquad\text{and}\qquad (\Phitil \ot \io)(Y)
= \cV (1 \ot \Wh_2) \; . $$
{\bf 2. The definition of $\tau$.}

a) The comultiplication $\de = (\Phi \ot \Phi) \Ga \Phi^{-1}$  on $M$ makes of $(M,\de)$ a
l.c.\ quantum group and then, $\Phi : (N,\Ga) \recht (M,\de)$ is an isomorphism. We claim
that, writing $\Wtil = W_{1,12} \cU^*$ one has $$(\io \ot \deop)(\Wtil) = (\Wtil \ot 1 \ot 1)
\; ((\io \ot \al)\be \ot \io \ot \io)(\Wtil) \; (\io \ot \al \ot \al)(\cV) \; .$$ Indeed,
observe that, by definition of $Y$ and Eq.~\eqref{eqY}, we have $$\cWh = (\ga \ot \io)(X) (\io
\ot \rho)(Y) \; .$$ From this it follows that
\begin{align*}
(\io \ot \Gaop)(\cWh) &= (\ga \ot \Gaop)(X) \; (\io \ot \rho \ot \rho)(\io \ot \deopt)(Y) \\
&= (\ga \ot \Gaop)(X) \; (\io \ot \rho \ot \rho)(\cVtil^* Y_{12} Y_{13}) \\
&= (\ga \ot \Gaop)(X) \; (\ga \ot \rho \ot \rho)(\cV^*) \; (\io \ot \rho)(Y)_{12} \; (\io \ot
\rho)(Y)_{13} \; .
\end{align*}
On the other hand,
$$(\io \ot \Gaop)(\cWh) = \cWh_{12} \cWh_{13} = (\ga \ot \io)(X)_{12} (\io \ot \rho)(Y)_{12} \;
(\ga \ot \io)(X)_{13} (\io \ot \rho)(Y)_{13} \; .$$
Both computations together yield
\begin{equation}\label{eq2}
(\ga \ot \Gaop)(X) = (\ga \ot \io)(X)_{12} (\io \ot \rho)(Y)_{12} (\ga \ot \io)(X)_{13} (\io
\ot \rho)(Y^*)_{12} (\ga \ot \rho \ot \rho)(\cV) \; .
\end{equation}
But
\begin{align*}
(\io \ot \rho)(Y)_{12} (\ga \ot \io)(X)_{13} (\io \ot \rho)(Y^*)_{12} &= (\io \ot \rho \ot \io)
\bigl( Y_{12} (\ga \ot \io)(X)_{13} Y_{12}^* \bigr) \\
&= (\io \ot \rho \ot \io)(\betil \ga \ot \io)(X) = (\ga \ot \rho \ot \io)(\be \ot \io)(X) \; .
\end{align*}
Then one can conclude, applying $\ga^{-1} \ot \io \ot \io$ to Eq.~\eqref{eq2}, that
$$(\io \ot \Gaop)(X) = X_{12} \; (\io \ot \rho \ot \io)(\be \ot \io)(X) \; (\io \ot \rho \ot
\rho)(\cV) \; .$$
Apply $\io \ot \Phi \ot \Phi$ to this expression, then
\begin{equation} \label{eq3}
(\io \ot \deop)(\Wtil) = (\Wtil \ot 1 \ot 1) \; ((\io \ot \al)\be \ot \io \ot \io)(\Wtil) \;
(\io \ot \al \ot \al)(\cV)
\end{equation}
and that is precisely the claim. The formula $\Ga \rho = (\rho \ot \rho)\det$ implies
\begin{equation} \label{eq4}
\de \al = (\al \ot \al)\det \; .
\end{equation}
which together with the previous formula determines completely the comultiplication $\de$ on $M$.

b) Repeating the proof of Proposition~\ref{28} and using \eqref{eq3}, \eqref{eq4}, we prove
that the dual weight $\vfi$ on $M = \Mo \kruisje{\al,\cU} \Mt$ of the weight $\vfi_2$ on
$\Mt$, with canonical GNS-construction $(H_1 \ot H_2,\io,\la)$ in the sense of
Terminology~\ref{115}, is left invariant on the l.c.\ quantum group $(M,\de)$. Defining $$\Wh
= (\be \ot \io \ot \io)\bigl( (W_1 \ot 1) \cU^* \bigr) \; (\io \ot \io \ot \al)\bigl( \cV (1
\ot \Wh_2) \bigr) \qquad\text{and}\qquad W = \Si \Wh^* \Si,$$ where $\Si$ flips $H \ot H$ and
$H = H_1 \ot H_2$, we get that $$\la \bigl( (\om \ot \io)\de(z) \bigr) = (\om \ot \io)(W^*)
\la(z) \tekst{for all} z \in \Nfi, \om \in M_* \; .$$ So $W$ is the multiplicative unitary
associated with the l.c.\ quantum group $(M,\de)$ and the GNS-construction $(H_1 \ot
H_2,\io,\la)$ for its left invariant weight $\vfi$.

c) The  comultiplication $\detil = (\Phitil \ot \Phitil)\Gah \Phitil^{-1}$ makes
$(\Mtil,\detil)$ a l.c.\ quantum group and $\Phitil : (\Nh,\Gah) \recht (\Mtil,\detil)$ is an
isomorphism. Exactly as above for $\de$, now starting from $$\cWh = (\ga \ot \io)(X) (\io \ot
\rho)(Y) \qquad\text{and}\qquad (\Gah \ot \io)(\cWh) = \cWh_{13} \cWh_{23} \; ,$$ one can
prove that $$(\detil \ot \io)(\Wc) = (1 \ot 1 \ot \Wc) \; (\io \ot \io \ot (\be \ot
\io)\al)(\Wc) \; (\be \ot \be \ot \io)(\cU) \qquad\text{and}\qquad \detil \be = (\be \ot
\be)\deo,$$ where $\Wc=(1 \ot \Wh_2^*) \cV^*$. Analogously as above, we get that the dual
weight $\vfitil$ on $\Mtil$ of the weight $\vfi_1$ on $\Mo$, with canonical GNS-construction
$(H_1\ot H_2,\io,\latil)$, is left invariant on $(\Mtil,\detil)$, and $$\latil \bigl( (\om \ot
\io) \detil(z) \bigr) = (\io \ot \om)(W) \latil(z) \tekst{for all} z \in \cN_{\vfitil}, \om
\in \Mtil_* \; .$$ Then we may conclude that $(\Mtil,\detil)$ is precisely the dual of
$(M,\de)$.

d) Apply \cite{KV1}, Propositions 8.14 and 8.15 to the l.c.\ quantum group $(M,\de)$ to get a
canonical left invariant weight $\vfih$ on $(\Mh,\deh)$ with canonical GNS-construction $(H_1
\ot H_2,\io,\lah)$. We already have $(\Mh,\deh)=(\Mtil,\detil)$ and we have the left invariant
weight $\vfitil$ on $(\Mtil,\detil)$ with GNS-construction $(H_1 \ot H_2,\io,\latil)$. Just as
in the beginning of the proof of Proposition~\ref{JJhat}, we get $\lambda \in \C \setminus
\{0\}$ such that $\lah = \lambda \latil$. Denote by $\Jh$ the modular conjugation of $\vfih$
in the GNS-construction $(H_1 \ot H_2,\io,\lah)$. Putting $u = \bar{\lambda}/\lambda$ we get
that $u \Jh$ is the modular conjugation of $\vfitil$ in the GNS-construction $(H_1 \ot
H_2,\io,\latil)$. Combining this with Proposition~\ref{117}, we have $$\be(x) = u \Jh (J_1 \ot
\Jh_2) (x \ot 1) (J_1 \ot \Jh_2) \Jh \bar{u} = \Jh (J_1 \ot \Jh_2) (x \ot 1) (J_1 \ot \Jh_2)
\Jh$$ for all $x \in M_1$. Further we get $\Jh (J_1 \ot \Jh_2) \in B(H_1) \ot M_2$. If $J$ is
the modular conjugation of $\vfi$ in the GNS-construction $(H_1 \ot H_2,\io,\la)$, it follows
from Proposition~\ref{117} that $$\al(x) = J (\Jh_1 \ot J_2) (1 \ot x) (\Jh_1 \ot J_2) J $$
for all $x \in M_2$. Further we get $J (\Jh_1 \ot J_2) \in M_1 \ot B(H_2)$. Define a unitary
$\cR$ on $H_1 \ot H_2$ by $$\cR = J \Jh (J_1 \Jh_1 \ot \Jh_2 J_2) \; .$$ For all $x \in M_2$
one has, using that $\Jh (J_1 \ot \Jh_2) \in B(H_1) \ot M_2$,
\begin{align*}
\cR (1 \ot x) \cR^* &= J \; \Jh (J_1 \ot \Jh_2) \; (1 \ot J_2 x J_2) \; (J_1 \ot \Jh_2) \Jh \; J
\\ &= J (1 \ot J_2 x J_2) J = \al(x) \; .
\end{align*}
From \cite{KV2}, Corollary 2.12 we get the existence of a complex number $\nu$ such that $$\cR
= \nu \; \Jh J (\Jh_1 J_1 \ot J_2 \Jh_2) \; .$$ Using this and the fact that $J(\Jh_1 \ot J_2)
\in M_1 \ot B(H_2)$ we get, for all $x \in M_1$,
\begin{align*}
\cR (x \ot 1) \cR^* &= \Jh \; J(\Jh_1 \ot J_2) \; (J_1 x J_1 \ot 1) \; (\Jh_1 \ot J_2) J \; \Jh
\\ &= \Jh (J_1 x J_1 \ot 1) \Jh = \be(x)
\end{align*}
Taking both computations together we
can define a faithful $*$-homomorphism $\tau$
$$\tau : \Mo \ot \Mt \recht \Mo \ot \Mt : \tau(z) = \cR z \cR^*$$
satisfying $\al(x) = \tau(1 \ot x)$ for all $x \in \Mt$ and $\be(x) = \tau(x \ot 1)$ for all $x
\in \Mo$.

{\bf 3. The proof of the compatibility relations.}

First we prove the relation
$$\cV_{132} \; (\io \ot \det) \al(x) \; \cV^*_{132} = \tau_{13} (\al \ot \io) \det(x)$$
for all $x \in \Mt$. Recall that
$$\Wh = (\be \ot \io \ot \io)(\Wtil) \; (\io \ot \io \ot \al)(\Wc^*) \tekst{where} \Wtil = (W_1
\ot 1)\cU^* \tekst{and} \Wc = (1 \ot \Wh_2^*)\cV^*  \; .$$
Hence we obtain, for all $z \in M$,
$$\deop(z) = \Wh(z \ot 1) \Wh^* =
(\be \ot \io \ot \io)(\Wtil) \; (\io \ot \io \ot \al) \bigl( \cV (\io \ot \deopt)(z)
\cV^* \bigr) \; (\be \ot \io \ot \io)(\Wtil^*) \; .$$
Using that $\de(\al(x)) = (\al \ot \al)\det(x)$ for all $x \in \Mt$, and using the previous
formula, we get for all $x \in \Mt$
\begin{align*}
(\io \ot \io \ot \al) \bigl( \cV (\io \ot \deopt)(\al(x)) \cV^* \bigr) &= (\be \ot \io \ot
\io)(\Wtil^*) \; \deop(\al(x)) \; (\be \ot \io \ot \io)(\Wtil) \\
&= (\be \ot \io \ot \io)(\Wtil^*) \; (\al \ot \al)\deopt(x) \; (\be \ot \io \ot \io)(\Wtil) \\
&= (\tau \ot \io \ot \io) \bigl( \Wtil^*_{134} \; \bigl( (\io \ot \al)\deopt(x) \bigr)_{234} \;
\Wtil_{134} \bigr) \; .
\end{align*}
Because $(\al,\cU)$ is a cocycle action of $(\Mo,\deo)$ on $\Mt$ we have $$\Wtil^* (1 \ot
\al(y)) \Wtil = (\io \ot \al)\al(y) \tekst{for all} y \in \Mt \; .$$ Using this, we can
continue the previous computation as
\begin{align*}
(\io \ot \io \ot \al) \bigl( \cV (\io \ot \deopt)(\al(x)) \cV^* \bigr) &= (\tau \ot \io \ot \io)
\bigl( (\io \ot (\io \ot \al)\al) \deopt(x)_{2134} \bigr) \\
&=( \io \ot \io \ot \al) \tau_{12} \si_{12} (\io \ot \al)\deopt(x) \; .
\end{align*}
Because $\al$ is faithful we may conclude that
$$\cV \; (\io \ot \deopt)\al(x) \; \cV^* = \tau_{12} \si_{12} (\io \ot \al)\deopt(x)$$
for all $x \in \Mt$, from where the needed relation follows.

Using now the formula $\deh(\be(x)) = (\be \ot \be) \deo(x)$ for all $x \in \Mo$ and the fact
that $(\si\be,\cV_{321})$ is a cocycle action of $(\Mt,\det)$ on $\Mo$ we obtain similarly, for all $x \in \Mo$,
$$\cU \; (\deo \ot \io)\be(x) \; \cU^* = \tau_{23} \si_{23} (\be \ot \io)\deo(x).$$
b) Let us prove the relation
$$(\de_1 \ot \io \ot \io)(\cV)(\io \ot \io \ot \deopt)(\cU^*) = (\cU^* \ot 1) (\io \ot
\tau \si \ot \io) \bigl( (\be \ot \io \ot \io)(\cU^*) (\io \ot \io
\ot \al)(\cV) \bigr) (1 \ot \cV) \; .$$
From Eq.~\eqref{eq3} one has
\begin{align*}
(\Wtil \; \ot \; &  1 \ot 1) \; ((\io \ot \al)\be \ot \io \ot \io)(\Wtil) \;
(\io \ot \al \ot \al)(\cV) = (\io \ot \deop)(\Wtil) \\
&= \bigl(1 \ot (\be \ot \io \ot \io)(\Wtil) \bigr) \; (\io \ot \io \ot \io \ot \al)
\bigl( (1 \ot \cV) (\io \ot \io \ot \deopt)(\Wtil) (1 \ot \cV^*) \bigr) \;
\bigl(1 \ot (\be \ot \io \ot \io)(\Wtil^*) \bigr)
\end{align*}
and so
\begin{multline*}
(\io \ot \io \ot \io \ot \al) \bigl( (1 \ot \cV) (\io \ot \io \ot \deopt)(\Wtil) (1 \ot \cV^*)
\bigr) \\ = \bigl(1 \ot (\be \ot \io \ot \io)(\Wtil^*) \bigr) \;
(\Wtil \ot 1 \ot 1) \; ((\io \ot \al)\be \ot \io \ot \io)(\Wtil) \;
(\io \ot \al \ot \al)(\cV) \; \bigl(1 \ot (\be \ot \io \ot \io)(\Wtil) \bigr) \; .
\end{multline*}
Using now that for all $x \in \Mo$ we have
$$\Wtil^* (1 \ot \be(x)) \Wtil = \cU \; (\deo \ot \io)\be(x) \; \cU^* = \tau_{23}\si_{23} (\be
\ot \io)\deo(x)$$
we obtain
\begin{align*}
(\io  & \ot \io \ot \io \ot \al) \; \bigl( (1 \ot \cV) (\io \ot \io \ot \deopt)(\Wtil) (1 \ot \cV^*)
\bigr) \\ &= (\Wtil \ot 1 \ot 1) \; \tau_{23} \si_{23} \bigl( (\be \ot \io \ot \io \ot
\io)(\deo \ot \io \ot \io)(\Wtil^*) \bigr) \; ((\io \ot \al)\be \ot \io \ot \io)(\Wtil) \;
(\io \ot \al \ot \al)(\cV) \\ &\hspace{12cm} \bigl(1 \ot (\be \ot \io \ot \io)(\Wtil) \bigr) \\
&= (\Wtil \ot 1 \ot 1) \; \tau_{23}\si_{23} \bigl( (\be \ot \io \ot \io \ot
\io)(\deo \ot \io \ot \io)(\Wtil^*) \; (\be \ot \io \ot \io)(\Wtil)_{1245} \; (\io \ot \io \ot
\al)(\cV)_{1245} \; \Wtil_{345} \bigr) \; .
\end{align*}
Using now that
$$\Wtil^* (1 \ot \al(x)) \Wtil = (\io \ot \al)\al(x)$$
for all $x \in \Mo$ we get
\begin{align*}
(\io  & \ot  \io \ot \io \ot \al) \;  \bigl( (1 \ot \cV) (\io \ot \io \ot \deopt)(\Wtil) (1 \ot \cV^*)
\bigr) \\ &= (\Wtil \ot 1 \ot 1) \; \tau_{23}\si_{23} \bigl( (\be \ot \io \ot \io \ot
\io)(\deo \ot \io \ot \io)(\Wtil^*) \; (\be \ot \io \ot \io)(\Wtil)_{1245} \; \Wtil_{345} \;
(\io \ot \io \ot (\io \ot \al)\al)(\cV) \bigr) \\
&= (\Wtil \ot 1 \ot 1) \; \tau_{23}\si_{23} \bigl( (\be \ot \io \ot \io \ot \io) \bigl(
(\deo \ot \io \ot \io)(\Wtil^*) \; \Wtil_{134} \Wtil_{234} \bigr) \;
(\io \ot \io \ot (\io \ot \al)\al)(\cV) \bigr) \; .
\end{align*}
Before we continue this computation, we use the cocycle property of $\cU$ to obtain
\begin{align*}
(\deo \ot \io \ot \io)(\Wtil^*) \; \Wtil_{134} \; \Wtil_{234} &=
(\deo \ot \io \ot \io)(\cU) W^*_{1,23} W^*_{1,13} \; W_{1,13} \cU^*_{134} \; W_{1,23}
\cU^*_{234} \\
&= (\deo \ot \io \ot \io)(\cU) \; (\io \ot \deo \ot \io)(\cU^*) \; \cU^*_{234} \\
&= (\io \ot \io \ot \al)(\cU^*) \; .
\end{align*}
Then we can continue the computation above and obtain
\begin{align*}
(\io \ot \io \ot \io \ot \al) \; & \bigl( (1 \ot \cV) (\io \ot \io \ot \deopt)(\Wtil) (1 \ot \cV^*)
\bigr) \\ &= (\Wtil \ot 1 \ot 1) \; \tau_{23}\si_{23} \bigl( (\be \ot \io \ot \al)(\cU^*) \;
(\io \ot \io \ot (\io \ot \al)\al)(\cV) \bigr) \\
&= (\io \ot \io \ot \io \ot \al) \bigl( (\Wtil \ot 1) \; \tau_{23} \si_{23} \bigl( (\be \ot \io \ot
\io)(\cU^*) \; (\io \ot \io \ot \al)(\cV) \bigr) \bigr) \; .
\end{align*}
Because $\al$ is faithful it now follows that
$$(1 \ot \cV) (\io \ot \io \ot \deopt)(\Wtil) (1 \ot \cV^*) =
(\Wtil \ot 1) \; \tau_{23} \si_{23} \bigl( (\be \ot \io \ot
\io)(\cU^*) \; (\io \ot \io \ot \al)(\cV) \bigr) \; .$$
Because $\Wtil = (W_1 \ot 1) \cU^*$ this gives us
$$(1 \ot \cV) W_{1,12} \; (\io \ot \io \ot \deopt)(\cU^*) \; (1 \ot \cV^*) =
W_{1,12} (\cU^* \ot 1) \;  \tau_{23} \si_{23} \bigl( (\be \ot \io \ot
\io)(\cU^*) \; (\io \ot \io \ot \al)(\cV) \bigr)$$
so that finally
$$(\deo \ot \io \ot \io)(\cV) \; (\io \ot \io \ot \deopt)(\cU^*) \; (1 \ot \cV^*) = (\cU^* \ot 1) \; (\io \ot
\tau\si \ot \io) \bigl( (\be \ot \io \ot \io)(\cU^*) \; (\io \ot \io \ot \al)(\cV) \bigr) \; .$$
c) It follows from Eq.~\eqref{eq3} and \eqref{eq4}
that $(M,\de)$ is indeed the cocycle bicrossed product of $(\Mo,\deo)$ and $(\Mt,\det)$.
We also have the isomorphism $\Phi: (N,\Ga) \recht (M,\de)$ satisfying $\Phi \rho = \al$. Let
$\hat{\Phi}$ be the canonically associated isomorphism $\hat{\Phi} : (\Nh,\Gah) \recht
(\Mh,\deh)$ characterized by
$$(\Phi \ot \hat{\Phi})(\cW) = W \; .$$
Since $\cWh = (\ga \ot \io)(X) (\io \ot \rho)(Y)$, one has $$(\Phi \ot \Phitil)(\cW) = W \; .$$
So we obtain $\Phitil = \hat{\Phi}$ and $\hat{\Phi} \ga = \be$, which concludes the proof.
\end{proof}
\subsection{Isomorphisms of cleft extensions}
\begin{definition}\label{35bis}
Let $(\Mo,\deo)$, $(\Mt,\det)$, $(\Ma,\dea)$ and $(\Mb,\deb)$ be l.c.\ quantum groups.
Extensions $$ (\Mt,\det) \overset{\ala}{\lrecht} (\Ma,\dea) \overset{\bea}{\lrecht}
(\Moh,\deoh) \qquad\text{and}\qquad (\Mt,\det) \overset{\alb}{\lrecht} (\Mb,\deb)
\overset{\beb}{\lrecht} (\Moh,\deoh) $$ are said to be isomorphic, if there exists an
isomorphism of l.c.\ quantum groups $$\pi : (\Ma,\dea) \recht (\Mb,\deb)$$ satisfying $\pi
\ala = \alb$ and $\pih \bea = \beb$, where $\pih$ is the canonical isomorphism of
$(\Mah,\deah)$ onto $(\Mbh,\debh)$ associated with $\pi$.
\end{definition}
Then we can prove the following result.
\begin{proposition} \label{35ter}
Let $(\Mo,\deo)$ and $(\Mt,\det)$ be l.c.\ quantum groups, let $(\taua,\cUa,\cVa)$ and
$(\taub,\cUb,\cVb)$ be two cocycle matchings of $(\Mo,\deo)$ and $(\Mt,\det)$ with
corresponding actions $\ala,\alb,\bea$ and $\beb$, and let $(\Ma,\dea)$ and $(\Mb,\deb)$ be
the respective cocycle bicrossed products. If the extensions $$ (\Mt,\det)
\overset{\ala}{\lrecht} (\Ma,\dea) \overset{\bea}{\lrecht} (\Moh,\deoh) \qquad\text{and}\qquad
(\Mt,\det) \overset{\alb}{\lrecht} (\Mb,\deb) \overset{\beb}{\lrecht} (\Moh,\deoh) $$ are
isomorphic through the isomorphism $\pi$, then there exists a unitary $\cR$ in $\Mo \ot \Mt$
satisfying
\begin{align*}
\taub(z) &= \cR \; \taua(z) \; \cR^* \tekst{for all} z \in \Mo \ot \Mt \\
\cUb &= (1 \ot \cR) \; (\io \ot \ala)(\cR) \; \cUa \; (\deo \ot \io)(\cR^*) \\
\cVb &= (\cR \ot 1) \; (\bea \ot \io)(\cR) \; \cVa \; (\io \ot \deopt)(\cR^*) \\
\pi(z) &= \cR \; z \; \cR^* \tekst{for all} z \in \Ma \\
\pih(z) &= \cR \; z \; \cR^* \tekst{for all} z \in \Mah \; .
\end{align*}
If conversely $(\taua,\cUa,\cVa)$ is a cocycle matching of $(\Mo,\deo)$ and $(\Mt,\det)$, and
if $\cR$ is a unitary in $\Mo \ot \Mt$, then the formulas above define a cocycle matching
$(\taub,\cUb,\cVb)$ of $(\Mo,\deo)$ and $(\Mt,\det)$. If $(\Ma,\dea)$ and $(\Mb,\deb)$ are the
corresponding cocycle bicrossed products, one can define an isomorphism $$\pi : (\Ma,\dea)
\recht (\Mb,\deb) : \pi(z) = \cR \; z \; \cR^*$$ of l.c.\ quantum groups such that $\pi \ala =
\alb$ and $\pih \bea =\beb$.
\end{proposition}
Observe that when both $\Mo$ and $\Mt$ are commutative, the extensions given by cocycle
matchings $(\taua,\cUa,\cVa)$ and $(\taub,\cUb,\cVb)$ are isomorphic if and only if
$\taua=\taub$ and there exists a unitary $\cR$ in $\Mo \ot \Mt$ satisfying
\begin{align*}
\cUb &= (1 \ot \cR) \; (\io \ot \ala)(\cR) \; \cUa \; (\deo \ot \io)(\cR^*) \\
\cVb &= (\cR \ot 1) \; (\bea \ot \io)(\cR) \; \cVa \; (\io \ot \deopt)(\cR^*) \; .
\end{align*}
We will come back to this situation in Subsection~\ref{secgroup}.
\begin{proof}
First let $\pi$ be an isomorphism of $(\Ma,\dea)$ onto $(\Mb,\deb)$ satisfying $\pi \ala =
\alb$ and $\pih \bea= \beb$. Recall that $\pih$ satisfies $(\pi \ot \pih)(\Wa)=\Wb$ where
$\Wa$ and $\Wb$ denote the multiplicative unitaries of $(\Ma,\dea)$ and $(\Mb,\deb)$
respectively. Concerning the dual actions $\alha$ and $\alhb$ of $(\Moh,\dehopo)$ on $\Ma$ and
$\Mb$ respectively, it follows from Proposition~\ref{14} and Definition~\ref{22} that $$(\io
\ot \io \ot \alha)(\Wha) = (\bea \ot \io)(W_1)_{123} \; \Whas{1245} \qquad\text{and}\qquad
(\io \ot \io \ot \alhb)(\Whb) = (\beb \ot \io)(W_1)_{123} \; \Whbs{1245} \; . $$ Hence we get
\begin{align*}
(\io \ot \io \ot (\io \ot \pi)\alha)(\Wha) &= (\bea \ot \io)(W_1)_{123} \; (\io \ot \io \ot
\pi)(\Wha)_{1245} = (\pih^{-1} \ot \io \ot \io \ot \io) \bigl( (\beb \ot \io)(W_1)_{123} \;
\Whbs{1245} \bigr) \\ &= (\pih^{-1} \ot \alhb)(\Whb) = (\io \ot \io \ot \alhb \pi)(\Wha) \; .
\end{align*}
So $\alhb \pi = (\io \ot \pi)\alha$.

As before we define the unitaries $\Wtila$ and $\Wtilb$ as
$$\Wtila= (W_1 \ot 1) \cUa^* \qquad\text{and}\qquad \Wtilb = (W_1 \ot 1) \cUb^* \; .$$
Then we define a unitary $Y := \Wtilb^* (\io \ot
\pi)(\Wtila)\in \Mo \ot \Mb$. Proposition~\ref{14} and the formula above imply
$$(\io \ot \alhb)(Y) = \Wtilbs{134}^* W^*_{1,12} \; (\io \ot \io \ot \pi) ( W_{1,12}
\Wtilas{134} ) = Y_{134} \; .$$
Therefore, $Y \in \Mo \ot \Mb^{\alhb}$. Theorem~\ref{110} shows that there is a unitary $\cR \in \Mo \ot \Mt$ such that $Y = (\io \ot \alb)(\cR)$, hence
\begin{equation} \label{v3}
(\io \ot \pi)(\Wtila) = \Wtilb (\io \ot \alb)(\cR) \; .
\end{equation}
Let $\vfia$ and $\vfib$ be the canonical left invariant weights on $(\Ma,\dea)$ and
$(\Mb,\deb)$ with GNS-constructions $(H_1 \ot H_2,\io,\laa)$ and $(H_1 \ot H_2,\io,\lab)$. We
claim that for all $z \in \cN_{\vfia}$ we have $\pi(z) \in \cN_{\vfib}$ and $$\lab(\pi(z)) =
\cR \laa(z) \; .$$ To prove this, choose $\xi \in H_1, b \in \cT_{\vfi_1}, x \in \cN_{\vfi_2}$
and an orthonormal basis $(e_i)_{i \in I}$ in $H_1$. Define the element $$z : =
(\om_{\xi,\la_1(b)} \ot \io \ot \io)(\Wtila) \; \ala(x) \in \cN_{\vfia}\; .$$ So, with \strong
convergence $\pi(z)$ is given by $$\pi(z) = \sum_{i \in I} (\om_{e_i,\la_1(b)} \ot \io \ot
\io)(\Wtilb) \; \alb \bigl( (\om_{\xi,e_i} \ot \io)(\cR) x \bigr) \; .$$ If we denote by
$y_{I_0}$ the sum over a finite subset $I_0 \subset I$, we get that $y_{I_0} \in \cN_{\vfib}$
and $$\lab(y_{I_0}) = \sum_{i \in I_0} J_1 \si^1_{i/2}(b) J_1 e_i \ot (\om_{\xi,e_i} \ot
\io)(\cR) \la_2(x) \; .$$ Since $\lab$ is \strong -- norm closed, we get that $\pi(z) \in
\cN_{\vfib}$ and $$\lab(\pi(z)) = (J_1 \si^1_{i/2}(b) J_1 \ot 1) \cR (\xi \ot \la_2(x)) = \cR
\laa(z) \; .$$ Because such elements $z$ form a \strong -- norm core for $\laa$, the claim is
proved. Similarly $z \in \cN_{\vfia}$ if and only if $\pi(z) \in \cN_{\vfib}$. From the
formula $\lab(\pi(z)) = \cR \laa(z)$ for all $z \in \cN_{\vfia}$ one can conclude that
$$\pi(z) = \cR \; z \; \cR^* \tekst{for all} z \in \Ma \qquad\text{and}\qquad \pih(z) = \cR \;
z \; \cR^* \tekst{for all} z \in \Mah \; .$$ Then also $\alb(x) = \pi(\ala(x)) = \cR \;
\ala(x) \; \cR^*$ for all $x \in \Mt$ and $\beb(y) = \pih(\bea(x)) = \cR \; \bea(x) \; \cR^*$
for all $x \in \Mo$. Therefore, $\taub(z) = \cR \; \taua(z) \; \cR^*$ for all $z \in \Mo \ot
\Mt$. So it follows from Eq.~\eqref{v3} that $$\cUb (W_1^* \ot 1) = \Wtilb^* = (\io \ot
\alb)(\cR) \; (\io \ot \pi)(\Wtila^*) = (1 \ot \cR) \; (\io \ot \ala)(\cR) \;  \cUa \; (W_1^*
\ot 1) \; (1 \ot \cR^*)$$ and hence $$\cUb = (1 \ot \cR) \; (\io \ot \ala)(\cR) \; \cUa \;
(\deo \ot \io)(\cR^*) \; .$$ It remains to prove the formula for $\cVb$. We have that $$\Wha =
(\bea \ot \io \ot \io)(\Wtila) \; (\io \ot \io \ot \ala)(\cVa (1 \ot \Wh_2))
\qquad\text{and}\qquad \Whb = (\beb \ot \io \ot \io)(\Wtilb) \; (\io \ot \io \ot \alb)(\cVb (1
\ot \Wh_2)) \; . $$ Because $\Whb = (\pih \ot \pi)(\Wha)$ and because of Eq.~\eqref{v3}, we
get that
\begin{align*}
(\beb \ot \io \ot \io)(\Wtilb) \; (\io \ot \io \ot \alb)(\cVb (1 \ot \Wh_2)) &= (\pih \ot
\pi)(\Wha) \\ &= (\beb \ot \io \ot \io) \bigl( \Wtilb (\io \ot \alb)(\cR) \bigr) \; (\io \ot \io
\ot \alb)(\pih \ot \io)(\cVa (1 \ot \Wh_2)) \; .
\end{align*}
From the faithfulness of $\alb$ it follows that $$\cVb (1 \ot \Wh_2) = (\beb \ot \io)(\cR) \;
(\pih \ot \io)(\cVa (1 \ot \Wh_2))$$ and so $$ \cVb = (\cR \ot 1) \; (\bea \ot \io)(\cR) \;
\cVa \; (1 \ot \Wh_2) (\cR^* \ot 1) (1 \ot \Wh_2^*) = (\cR \ot 1) \; (\bea \ot \io)(\cR) \;
\cVa \; (\io \ot \deopt)(\cR^*) \; . $$ This concludes the proof of the first statement. Vice
versa, let $(\taua,\cUa,\cVa)$ be a cocycle matching of $(\Mo,\deo)$ and $(\Mt,\det)$, let
$\cR$ be a unitary in $\Mo \ot \Mt$ and define $(\taub,\cUb,\cVb)$ by the formulas mentioned
in the proposition. Then it is straightforward to check that $(\taub,\cUb,\cVb)$ is a cocycle
matching of $(\Mo,\deo)$ and $(\Mt,\det)$. Let $(\Ma,\dea)$ and $(\Mb,\deb)$ be the
corresponding cocycle bicrossed products and let $\Wha$ and $\Whb$ be the corresponding
multiplicative unitaries. Then one can check that $$(\cR \ot \cR) \Wha (\cR^* \ot \cR^*) =
\Whb \; .$$ So, defining $\pi(z) = \cR \; z \; \cR^*\ (z \in \Ma)$, we get an isomorphism of
the l.c.\ quantum groups $(\Ma,\dea)$ and $(\Mb,\deb)$ such that $\pih(z) = \cR \; z \; \cR^*\
(z \in \Mah)$. Hence we get that $\pi \ala = \alb$ and $\pih \bea = \beb$, which concludes the
proof.
\end{proof}

\section{Extensions of l.c.\ groups}
In this section we study extensions of l.c.\ groups, i.e.,  short exact sequences of the form
$$(L^\infty(G_2),\det) \overset{\al}{\lrecht} (M,\de) \overset{\be}{\lrecht}
(\cL(G_1),\deoh),$$ where $G_1$ and $G_2$ are usual l.c.\ groups.
\subsection{Extensions of discrete groups} \label{subsec42}
If both $G_1$ and $G_2$ are finite, G.I.~Kac \cite{Kac}, Theorem 3 showed that every extension has a cocycle bicrossed product structure (later it was shown in \cite{A} that it follows from H.-J.~Schneider's work \cite{Schn} that every extension of finite-dimensional Hopf algebras is cleft). We generalize the result of G.I.~Kac as follows.
\begin{theorem} \label{36}
Any extension
$$(L^\infty(G_2),\det) \overset{\al}{\lrecht} (M,\de) \overset{\be}{\lrecht} (\cL(G_1),\deoh),$$
where $G_1$ and $G_2$ are discrete groups is automatically cleft; applying then Theorem~\ref{35}, we have that $(M,\de)$ is a cocycle bicrossed product of $(L^\infty(G_1),\de_1)$ and $(L^\infty(G_2),\det)$.
\end{theorem}
The proof of this theorem will consist of several lemmas,
very much inspired by \cite{Kac}.

Define, as in the section of Preliminaries, quantum groups $(\Mo,\deo)$ and $(\Mt,\det)$,
their duals and all their ingredients related to $G_1$ and $G_2$. For $g \in G_i$ denote by
$\delta_g$ the element in $M_i\ (i=1,2)$ determined by $$\delta_g(h) = \delta_{g,h},$$ where
$\delta_{g,h}$ is the Kronecker symbol; the elements $\delta_g$ form an orthonormal basis in
$H_i\ (i=1,2)$. Next define for all $g \in G_1$ the functional $\om_g \in \Mh_{1 \, *}$ by
$$\om_g(x) = \langle x \delta_e , \delta_g \rangle \tekst{for all} x \in \Mh_1$$ and observe
that $$(\io \ot \om_g)(W_1) = \delta_g \qquad\text{and}\qquad \om_g(\lambda_h) =
\delta_{g,h}$$ for all $g,h \in G_1$. Further, $\overline{\om_g} = \om_{g^{-1}}$ and $$(\io
\ot \om_g)\deoh(x) = \om_g(x) \lambda_g \tekst{for all} x \in \Mh_1 \; .$$

Let $\al:\Mt \recht M$ and $\be : \Mo \recht \Mh$ be normal
and faithful $*$-homomorphisms satisfying $\de \al = (\al \ot \al)\det$ and $\deh \be = (\be
\ot \be)\deo$. Denote by $\mu$ the action defined in Proposition~\ref{31} such that $\mu=(\hat R_1\ot R)\te R$, $M^\mu = \al(\Mt)$ and $(\io \ot \de)\mu = (\mu \ot \io) \de$.
We fix a left invariant weight $\vfi$ on $(M,\de)$ with GNS-construction $(H,\io,\la)$.
\begin{lemma} \label{38}
For all $g \in G_1$ and $s \in G_2$ there exists a $z \in M$ such that $$(\om_g \ot \io)\mu(z)
\; \al(\delta_s) \neq 0 \; .$$
\end{lemma}
\begin{proof}
We first claim that there exist $z \in \Nfi$ and $k \in G_2$ such that $\al(\delta_k)\; (\om_g
\ot \io)\mu(z) \neq 0$. If this is not the case, then for every $z \in \Nfi$ and $k \in G_2$
we have $\al(\delta_k)\; (\om_g \ot \io)\mu(z) = 0$. Hence, using Proposition~\ref{31} $$0 =
\la\bigl( \al(\delta_k)\; (\om_g \ot \io)\mu(z) \bigr) = \al(\delta_k) \; \be((\io \ot
\om_g)(W_1)) \; \la(z) = \al(\delta_k) \be(\delta_g) \la(z) \; .$$ Summing over $k \in G_2$,
it then follows that $\be(\delta_g)=0$ which contradicts the faithfulness of $\be$. So one can
take $k \in G_2$ and $z \in \Nfi$ such that $Z:=\al(\delta_k)\; (\om_g \ot \io)\mu(z) \neq 0$.
Since $Z \in \Nfi$, we have $0 < \vfi(Z^*Z) < \infty$. Fix $g \in G_1$ and define $$I_g = \{ s
\in G_2 \mid \; \text{There exists a} \; z \in M \tekst{with} (\om_g \ot \io)\mu(z) \;
\al(\delta_s) \neq 0 \} \; .$$ Suppose $s \in G_2 \setminus I_g$. One has $(\om_g \ot \io)
\mu(y) \; \al(\delta_s) = 0$ for all $y \in M$. We claim that $\de(Z) (\al(\delta_s) \ot 1) =
0$. Indeed, $$\de(Z) = (\al \ot \al)\det(\delta_k) \; \de \bigl( (\om_g \ot \io) \mu(z) \bigr)
= (\al \ot \al)\det(\delta_k) \; ((\om_g \ot \io)\mu \ot \io)\de(z) \; .$$ For all $\eta \in
M_*$ $$(\io \ot \eta) \bigl(  ((\om_g \ot \io)\mu \ot \io)\de(z) \bigr) \; \al(\delta_s)
=(\om_g \ot \io)\mu \bigl( (\io \ot \eta)\de(z) \bigr) \; \al(\delta_s) = 0$$ and so $$((\om_g
\ot \io)\mu \ot \io)\de(z) \; (\al(\delta_s) \ot 1) = 0 \;  .$$ From this the claim follows.
But then $$(\al(\delta_s) \ot 1) \de(Z^* Z) (\al(\delta_s) \ot 1) =0 \; .$$ Applying $\io \ot
\vfi$ we obtain $$\vfi(Z^*Z) \; \al(\delta_s) = 0$$ i.e., a contradiction with the
faithfulness of $\al$. Therefore, $G_2=I_g$, which concludes the proof.
\end{proof}
\begin{lemma} \label{39}
Let $g \in G_1$. Then
$$\{ (\om_g \ot \io) \mu(z) \mid z \in M \} = \{ z \in M \mid \mu(z) = \lambda_g \ot z \} \;
.$$
\end{lemma}
\begin{proof}
Let $z \in M$. Then
$$\mu\bigl( (\om_g \ot \io)\mu(z) \bigr) = ((\om_g \ot \io)\deoh \ot \io)\mu(z) = \lambda_g \ot
(\om_g \ot \io)\mu(z) \; .$$
This proves the inclusion $\subset$. Vice versa, let $z \in M$ and suppose
$\mu(z) = \lambda_g \ot z$. Then $z=(\om_g \ot \io)\mu(z)$, which gives the converse inclusion.
\end{proof}
\begin{lemma} \label{310}
Let $g \in G_1, s,k \in G_2$ and $z \in M$. If the element $u \in M$ defined by
$$u:= \al(\delta_k) \; (\om_g \ot \io)\mu(z) \; \al(\delta_s)$$
is different from $0$, then $\frac{u}{\|u\|}$ is a partial isometry with initial projection
$\al(\delta_s)$ and final projection $\al(\delta_k)$.
\end{lemma}
\begin{proof}
Define
$$v:= (\om_g \ot \io)\mu(z)^* \; \al(\delta_k) \; (\om_g \ot \io)\mu(z)
= (\om_{g^{-1}} \ot \io)\mu(z^*) \; \al(\delta_k) \; (\om_g \ot \io)\mu(z) \; .$$
Then $u^* u = \al(\delta_s) \; v \; \al(\delta_s)$ and further, using the previous lemma,
$$\mu(v) = \bigl(\lambda_{g^{-1}} \ot (\om_{g^{-1}} \ot \io)\mu(z^*) \bigr) \; \bigl( 1 \ot
\al(\delta_k) \bigr) \; \bigl( \lambda_g \ot (\om_g \ot \io)\mu(z) \bigr) = 1 \ot v \; .$$
Hence $v \in M^\mu = \al(M_2)$ and one can take $x \in \Mt$ such that $v=\al(x)$. Then
$$u^* u = \al(\delta_s x \delta_s) = x(s) \al(\delta_s) \; .$$
From this we have $u^*u = \|u\|^2 \al(\delta_s)$ and similarly $u u^* = \|u\|^2\al(\delta_k)$, which concludes the proof.
\end{proof}
\begin{lemma} \label{311}
Let $g \in G_1$ and $s \in G_2$. Then there exists at most one $k \in G_2$ such that
$$\{ \al(\delta_k) \; (\om_g \ot \io)\mu(z) \; \al(\delta_s) \mid z \in M \} \neq \{0\} \; .$$
\end{lemma}
\begin{proof}
If there are $k,l \in G_2$, both satisfying the condition above, one can take $y,z \in M$ such that
$$u:=\al(\delta_k) \; (\om_g \ot \io)\mu(y) \; \al(\delta_s) \neq 0 \qquad\text{and}\qquad
v:=\al(\delta_l) \; (\om_g \ot \io)\mu(z) \; \al(\delta_s) \neq 0 \; .$$
We may suppose that $\|u\|=\|v\|=1$. By Lemma~\ref{310} $u$ (resp., $v$) is a partial isometry with
initial projection $\al(\delta_s)$ and final projection $\al(\delta_k)$ (resp., $\al(\delta_l)$). Hence $u v^*$ is a
partial isometry with initial projection $\al(\delta_l)$ and final projection $\al(\delta_k)$.
On the other hand, Lemma~\ref{39} gives
$$\mu(u) = \lambda_g \ot u \qquad\text{and}\qquad \mu(v) = \lambda_g \ot v,$$
so $uv^* \in M^\mu = \al(\Mt)$. Choose $x \in \Mt$ with $\al(x) = uv^*$. Then
$$\al(\delta_l) = (uv^*)^* uv^* = \al(x^*x) =\al(x x^*) = uv^* (uv^*)^*= \al(\delta_k),$$
so $l=k$.
\end{proof}
The following proposition is the main ingredient to prove Theorem~\ref{36}.
\begin{proposition} \label{312}
Let $g \in G_1$ and $s \in G_2$. Then the subspace of $M$
$$E_{g,s} = \{ (\om_g \ot \io)\mu(z) \; \al(\delta_s) \mid z \in M \}$$
is one-dimensional and spanned by a partial isometry $u_{g,s}$ with initial
projection $\al(\delta_s)$ and final projection $\al(\delta_{P_g(s)})$, where $P_g(s) \in G_2$.
For every $g \in G_1$, the map $s \mapsto P_g(s)$ is a bijection of $G_2$. Further
$$\mu(u_{g,s}) =  \lambda_g \ot u_{g,s} \; .$$
\end{proposition}
\begin{proof}
By the previous lemma there is at most one $k \in G_2$ such that $\al(\delta_k) E_{g,s} \neq
\{0\}$. On the other hand, if such a $k$  does not exist, then $\al(\delta_k) E_{g,s} = \{0\}$
for all $k \in G_2$ and summing over $k$ we get $E_{g,s} = \{0\}$, which contradicts
Lemma~\ref{38}. So, let $P_g(s)$ be the unique element in $G_2$ satisfying
$$\al(\delta_{P_g(s)}) E_{g,s} \neq \{0\} \; .$$ Then $y = \al(\delta_{P_g(s)}) y$ for all $y
\in E_{g,s}$. Suppose now that $u,v \in E_{g,s}$ and $\|u\|=\|v\|=1$. By Lemma~\ref{310}, both
$u$ and $v$ are partial isometries with initial projection $\al(\delta_s)$ and final
projection $\al(\delta_{P_g(s)})$. Lemma~\ref{39} shows that all elements $y \in E_{g,s}$
satisfy $\mu(y) = \lambda_g \ot y$, and hence $v^* u \in M^\mu = \al(\Mt)$. But $v^* u$ is a
partial isometry with initial and final projection equal to $\al(\delta_s)$. Hence $v^*u = \pm
\al(\delta_s)$ and so $u = \pm v$. Hence $E_{g,s}$ is one-dimensional and generated by a
partial isometry $u_{g,s}$ with initial projection $\al(\delta_s)$ and final projection
$\al(\delta_{P_g(s)})$. Also we get $\mu(u_{g,s}) = \lambda_g \ot u_{g,s}$.

Fix now $g \in G_1$ and consider the map $s \mapsto P_g(s)$. First of all this map is
injective. Take $s \neq t$ and suppose that $P_g(s) = P_g(t)$. Writing $P_g(s)=k$, we can take
$y,z \in M$ such that $$\al(\delta_k) \; (\om_g \ot \io)\mu(y) \; \al(\delta_s) \neq 0
\qquad\text{and}\qquad \al(\delta_k) \; (\om_g \ot \io)\mu(z) \; \al(\delta_t) \neq 0 \; .$$
Taking the adjoint we get $$\al(\delta_s) \; (\om_{g^{-1}} \ot \io)\mu(y^*) \; \al(\delta_k)
\neq 0 \qquad\text{and}\qquad \al(\delta_t) \; (\om_{g^{-1}} \ot \io)\mu(z^*) \; \al(\delta_k)
\neq 0$$ and this contradicts Lemma~\ref{311}. Therefore, the map $s \mapsto P_g(s)$ is
injective. To prove its surjectivity, choose $t \in G_2$. By the above reasoning we can take
$s\in G_2$ and $y \in N$ such that $$\al(\delta_s) \; (\om_{g^{-1}} \ot \io)\mu(y) \;
\al(\delta_t) \neq 0 \; .$$ Taking the adjoint we get $$\al(\delta_t) \; (\om_g \ot
\io)\mu(y^*) \; \al(\delta_s) \neq 0$$ and this implies that $P_g(s) =t$.
\end{proof}
\begin{proof}[Proof of Theorem~\ref{36}]
Recall that the action $\te$ introduced in Proposition~\ref{31} satisfies
$\te(z) = (\hat{R}_1 \ot R) \mu(R(z))$ for all $z \in M$.
In order to prove that $(M,\de)$ is a cleft extension, it suffices, by
Proposition~\ref{122}, to show the existence of a unitary $X \in M_1 \ot M$ satisfying $(\io \ot
\te)(X) = W_{1,12} X_{13}$. Applying $R_1 \ot \hat{R}_1 \ot R$ we have to prove the existence
of a unitary $Y$ in $\Mo \ot M$ satisfying $(\io \ot \mu)(Y) = Y_{13} W_{1,12}$.

Choose partial isometries $u_{g,s}$ as in the previous proposition. Fix $g \in G_1$. Then
$(u_{g,s})_{s \in G_2}$ is a family of partial isometries in $M$ with initial projections
$(\al(\delta_s))_{s \in G_2}$ and final projections $(\al(\delta_{P_g(s)}))_{s \in G_2}$. Because
the map $s \mapsto P_g(s)$ is a bijection of $G_2$ we can define a unitary $Y_g$ in $M$ by the sum
$$Y_g = \sum_{s \in G_2} u_{g,s}$$
converging in the \strong topology. Since $\mu(u_{g,s})=\lambda_g \ot u_{g,s}$ for all $s \in G_2$, we get $\mu(Y_g) = \lambda_g \ot Y_g$.
Identifying $\Mo \ot M$ with $\ell^\infty(G_1,M)$ we can define $Y \in \Mo \ot M$ by $Y(g) =
Y_g$. Then $Y$ is unitary and
$$(\io \ot \mu)(Y) = Y_{13} W_{1,12} \; .$$
This concludes the proof.
\end{proof}
\subsection{Extensions of l.c.\ groups}
In this subsection we discuss the particular situation of a (cocycle) matched pair of locally
compact groups. Our starting point, Definition~\ref{41}, is the same as in \cite{B-S2}. Also, but in the
absence of cocycles, S.~Baaj and G.~Skandalis \cite{B-S2} state the same formulas for the ingredients
of the bicrossed product as we do in Propositions~\ref{48} and \ref{49}.

Here we consider regular Borel measures on l.c.\ spaces. By the product measure, we mean the
regular Borel product, see e.g.\ \cite{cohn}, Section~7.6. We say that a statement is valid
almost everywhere if it is valid everywhere except on a set whose intersection with an
arbitrary compact set is a Borel set of measure zero, see e.g.\ \cite{cohn}, Section~3.3 where
they use the terminology locally almost everywhere. A function is called measurable when it is
Borel measurable on every compact set. Given a l.c.\ space $X$ with a regular Borel measure
$\mu$, we denote by $K(X)$ the space of continuous compactly supported $\C$-valued functions
on $X$ and by $\cL^\infty(X)$ the space of measurable functions from $X$ to $\C$ such that
there exist a number $M >0$ satisfying $|f(x)| \leq M$ for almost all $x \in X$. Denote by
$L^2(X)$ the usual Hilbert space of (equivalence classes of) square integrable functions. For
every $f \in \cL^\infty(X)$ one can define a multiplication operator $M_f$ on $L^2(X)$ in the
usual way, and then we have a von Neumann algebra (this follows e.g.\ from \cite{cohn},
Theorem~7.5.3) $$L^\infty(X) = \{ M_f \mid f \in \cL^\infty(X) \}.$$
\begin{definition} \label{41}
Let $G,G_1$ and $G_2$ be l.c.\ groups with fixed left invariant Haar measures and let a
homomorphism $i: G_1 \recht G$ and an anti-homomorphism $j:G_2 \recht G$ have closed images
and be homeomorphisms onto these images.  Let finally $$\te : G_1 \times G_2 \recht \Om
\subset G : (g,s) \mapsto i(g)j(s)$$ be a homeomorphism of $G_1 \times G_2$ onto an open
subset $\Om$ of $G$ having a complement of measure zero. Then we call $G_1$ and $G_2$ a
matched pair of l.c.\ groups.
\end{definition}
From this data we will construct actions $\al$
and $\be$ and a map $\tau$ to obtain a cocycle matching of $(L^\infty(G_1),\deo)$ and
$(L^\infty(G_2),\det)$, with trivial cocycles, in the sense of Definition~\ref{21}.

Denote by $\sde,\sdeo$ and $\sdet$ the modular functions of $G,G_1$ and $G_2$ and define the homeomorphism
$$\rho : G_1 \times G_2 \recht \Om^{-1} : (g,s) \mapsto j(s)i(g) \; .$$
We put
$$\cO = \te^{-1}(\Om \cap \Om^{-1}) \qquad\text{and}\qquad \cO' = \rho^{-1}(\Om \cap \Om^{-1})
\; .$$
So $\cO$ and $\cO'$ are open subsets of $G_1 \times G_2$ and $\rho^{-1} \te$ is a homeomorphism
of $\cO$ onto $\cO'$. For all $(g,s) \in \cO$ we define $\be_s(g) \in G_1$ and $\al_g(s) \in
G_2$ such that
$$\rho^{-1}(\te(g,s)) = (\be_s(g),\al_g(s)) \; .$$
Hence we get $j \bigl( \al_g(s) \bigr) \; i\bigl( \be_s(g) \bigr) = i(g)j(s)$ for all $(g,s) \in \cO$.

Let us first prove an easy lemma.
\begin{lemma}\label{42}
Let $(g,s) \in \cO$ and $h \in G_1$. Then $(hg,s) \in \cO$ if and only if $(h,\al_g(s)) \in
\cO$, and in that case
$$\al_{hg}(s) = \al_h \bigl( \al_g(s) \bigr) \qquad\text{and}\qquad \be_s(hg) =
\be_{\al_g(s)}(h) \; \be_s(g) \; .$$
Let $(g,s) \in \cO$ and $t \in G_2$. Then $(g,ts) \in \cO$ if and only if $(\be_s(g),t) \in
\cO$ and in that case
$$\be_{ts}(g) = \be_t \bigl( \be_s(g) \bigr) \qquad\text{and}\qquad \al_g(ts) =
\al_{\be_s(g)}(t) \; \al_g(s) \; .$$
Finally, for all $g \in G_1$ and $s \in G_2$ we have $(g,e) \in \cO$, $(e,s) \in \cO$, and
$$\al_g(e) = e, \quad \al_e(s) = s, \quad \be_s(e)=e \tekst{and} \be_e(g)=g \; .$$
\end{lemma}
\begin{proof}
Let $(g,s) \in \cO$ and $(h,\al_g(s)) \in \cO$. Then we know that $$j \bigl( \al_h(\al_g(s))
\bigr) \; i \bigl( \be_{\al_g(s)}(h) \bigr) = i(h) \; j \bigl( \al_g(s) \bigr) \; .$$ So we
get that $$j \bigl( \al_h(\al_g(s)) \bigr) \; i \bigl( \be_{\al_g(s)}(h) \bigr) \; i \bigl(
\be_s(g) \bigr) = i(h) \; j \bigl( \al_g(s) \bigr) \; i \bigl( \be_s(g) \bigr) = i(hg) \; j(s)
\; .$$ From this we may conclude that $(hg,s) \in \cO$ and $\be_s(hg)=\be_{\al_g(s)}(h) \;
\be_s(g)$, $\al_{hg}(s) = \al_h \bigl( \al_g(s) \bigr)$. The other statements of the lemma are
proved analogously.
\end{proof}

\begin{lemma}\label{43}
The Haar measures on $G,G_1$ and $G_2$ can be normalized in such a way that
for all positive Borel functions $f$ on $G$ we have
$$
\int f \bigl( i(g) j(s) \bigr) \; \sde(j(s)) \; dg \times ds = \int f(x) \; dx = \int f \bigl( j(s)
i(g) \bigr) \; \sde(i(g)) \; \sdeo(g^{-1}) \; \sdet(s^{-1}) \; dg \times ds \; .$$
\end{lemma}
\begin{proof}
For any $f \in K(G_1 \times G_2)$ we define $\tilde{f} \in K(G_1 \times G_2)$ by
$\tilde{f}(g,s) = f(g,s) \sde(j(s)^{-1})$ and $$I(f) = \int_\Om \tilde{f}(\te^{-1}(x)) \; dx
\; .$$ We claim that $I$ is a left invariant integral on $K(G_1 \times G_2)$. To prove this,
choose $(h,t) \in G_1 \times G_2$. Choose $f \in K(G_1 \times G_2)$ and define $F \in K(G_1
\times G_2)$ by $F(g,s)=f(hg,ts)$. We have to prove that $I(F)=I(f)$.

Choose $x \in \Om$. Let $x = i(g)j(s)$. Then
$$\tilde{F}(\te^{-1}(x)) = \tilde{F}(g,s) = f(hg,ts) \sde(j(s)^{-1}) = \tilde{f}(hg,ts)
\sde(j(t)) \; .$$
Because $\te(hg,ts)=i(h) x j(t)$ we get that $\tilde{F}(\te^{-1}(x)) = \tilde{f} \bigl( \te^{-1}(i(h)
x j(t) ) \bigr) \sde(j(t))$ for all $x \in \Om$. Hence
$$I(F) = \int_\Om \tilde{F}(\te^{-1}(x)) \; dx = \int_\Om \tilde{f} \bigl( \te^{-1}(i(h)
x j(t) ) \bigr) \sde(j(t)) \; dx = \int_\Om \tilde{f}(\te^{-1}(x)) \; dx = I(f) \; .
$$
Hence it follows that $I$ is a left invariant integral on $K(G_1 \times G_2)$.
Then we can normalize the Haar measure on $G$ in such a way that the first equality of the lemma holds for all continuous functions on $G$
with compact support in $\Om$. The same equality then follows for all positive Borel
functions as usual. The second equality follows from the first and the formula
$$\int f(x) \; dx = \int f(x^{-1}) \sde(x^{-1}) \; dx \; .$$
\end{proof}
This lemma implies in particular that a Borel set $A \subset G_1 \times G_2$ has measure zero
if and only if $\te(A)$ has measure zero in $G$. Because $\Om \cap \Om^{-1}$ has complement of
measure zero in $G$, $\cO$ and $\cO'$ have complements of measure zero in $G_1 \times G_2$.
Hence we can define a $*$-isomorphism $\tau$ by
\begin{equation}\label{tau}
\tau : L^\infty(G_1) \ot L^\infty(G_2) \recht L^\infty(G_1) \ot L^\infty(G_2) : \tau(f)(g,s)
= f(\be_s(g),\al_g(s)) \; .
\end{equation}
Lemma~\ref{42} guarantees that, defining
\begin{align*}
&\al : L^\infty(G_2) \recht L^\infty(G_1) \ot L^\infty(G_2) : \al(f) =\tau(1 \ot f)
\tekst{and}  \\  &\be : L^\infty(G_1) \recht L^\infty(G_1) \ot L^\infty(G_2) :
\be(f) =\tau(f \ot 1)
\end{align*}
we have
\begin{alignat*}{2}
(\io \ot \al)\al(f) &= (\de_1 \ot \io)\al(f) &\qquad
(\be \ot \io)\be(f) &= (\io \ot \deopt)\be(f) \\
\tau_{13} (\al \ot \io) \de_2(f) &= (\io \ot \de_2)\al(f) &\qquad
\tau_{23}\si_{23}(\be \ot \io) \de_1(f) &=  (\de_1 \ot \io)\be(f) \; .
\end{alignat*}
Therefore, $\tau$ gives a cocycle matching with trivial cocycles of $(L^\infty(G_1),\deo)$
and $(L^\infty(G_2),\det)$ in the sense of Definition~\ref{21}. Then we can
introduce cocycles by the following obvious lemma.
\begin{lemma}\label{44}
Suppose that
$$\cU  : G_1 \times G_1 \times G_2 \recht \C \qquad\text{and}\qquad
\cV  : G_1 \times G_2 \times G_2 \recht \C$$
are measurable maps with values in the unit circle $U(1)\subset \C$, satisfying
\begin{align}
\cU(g,h,\al_k(s)) \; \cU(gh, k, s) &= \cU(h,k,s) \; \cU(g, hk, s) \label{cocU}, \\
\cV(\be_s(g),t,r) \; \cV(g,s, rt) &= \cV(g,s,t) \; \cV(g, ts, r) \notag ,\\
\cV(gh, s, t) \; \bar{\cU}(g,h,ts)
& = \bar{\cU}(g,h,s) \; \bar{\cU}(\be_{\al_h(s)}(g), \be_s(h) , t) \; \cV(g, \al_h(s),
\al_{\be_s(h)}(t)) \; \cV(h,s,t) \notag
\end{align}
nearly everywhere. Then $(\tau,\cU,\cV)$ is a cocycle
matching of $(L^\infty(G_1),\deo)$ and $(L^\infty(G_2),\det)$.
\end{lemma}
Fixing cocycles $\cU$ and $\cV$ as above and using Definition~\ref{22} and Theorem~\ref{213}
we get the cocycle bicrossed product $(M,\de)$, which is a l.c.\ quantum group.

Let $\vfi_2$ be the \nsf weight on $L^\infty(G_2)$ coming from the Haar measure on $G_2$ and $\vfi$ the left invariant weight on $(M,\de)$ according to Definition~\ref{27}.
The obvious GNS-construction for $\vfi_2$ gives also the canonical
GNS-construction $(L^2(G_1 \times G_2),\io,\la)$ for $\vfi$. Let $W_1$ be the multiplicative unitary of $G_1$ and $\Wtil = (W_1 \ot 1)\cU^*$. Recall that $M$ is generated by $\al(L^\infty(G_2))$ and the elements $(\om \ot \io \ot
\io)(\Wtil)$, $\om \in L^\infty(G_1)_*$.

For every function $F \in K(G_1)$ we define the normal functional $\om_F \in L^\infty(G_1)_*$
by $\om_F(g) = \int F(x) g(x) \; dx$ for all $g \in L^\infty(G_1)$. It follows from
Definition~\ref{27} that
\begin{equation}\label{core}
\lspan\{ (\om_F \ot \io \ot \io)(\Wtil) \al(G) \mid F \in K(G_1), G \in K(G_2) \}
\end{equation}
is a \strong -- norm core for $\la$ and
$$\la \bigl((\om_F \ot \io \ot \io)(\Wtil) \al(G) \bigr) = F \ot G \; .$$
\begin{lemma} \label{45}
Let $k \in G_1$. Define
$$\cO_k = \{ s \in G_2 \mid (k,s) \in \cO \} \; .$$
Then $\cO_k$ is an open subset of $G_2$ with complement of measure zero and $\al_k$ is a
homeomorphism of $\cO_k$ onto $\cO_{k^{-1}}$ satisfying,
for all positive Borel functions $h$ on $G_2$,
$$\int_{\cO_k} h \bigl( \al_k(s) \bigr) \; ds = \int_{\cO_{k^{-1}}} h(s) F(k^{-1},s) \; ds$$
Here $F$ is the strictly positive
continuous function on $\cO$ defined by
$$F(k,s) = \sde \bigl( i(\be_s(k)^{-1}) \bigr) \; \sdeo\bigl( \be_s(k) \bigr) \; \sdet\bigl(
\al_k(s) s^{-1} \bigr) \; .$$
\end{lemma}
\begin{proof}
Lemma~\ref{43} shows that $\rho^{-1}$ preserves Borel sets of measure zero. Now
$i(k^{-1})\Om^{-1} \cap \Om^{-1}$ has complement of measure zero in $G$ and its image under
$\rho^{-1}$ is $G_1 \times \cO_k$. Hence $\cO_k$ is an open subset of $G_2$ with
complement of measure zero. Using twice Lemma~\ref{43}, we get for every continuous function
$h$ on $G$ with compact support in $\Om^{-1} \cap i(k)\Om^{-1}$
\begin{align*}
\int h \bigl( j(s) i(g) \bigr) \; & \sde(i(g)) \; \sdeo(g^{-1}) \; \sdet(s^{-1}) \; dg \times ds
= \int h(x) \; dx = \int h(i(k) x) \; dx \\
&= \int_{G_1 \times \cO_k} h \bigl( i(k) j(s) i(g) \bigr) \; \sde(i(g)) \; \sdeo(g^{-1}) \; \sdet(s^{-1}) \; dg \times
ds  \\ &= \int_{G_1 \times \cO_k} h \bigl( j(\al_k(s)) i( \be_s(k) g ) \bigr) \; \sde(i(g)) \; \sdeo(g^{-1}) \; \sdet(s^{-1}) \; dg \times
ds \\ &= \int_{G_1 \times \cO_k} h \bigl(j(\al_k(s)) i(g) \bigr) \; \sde \bigl( i(
\be_s(k)^{-1} g) \bigr) \; \sdeo \bigl( g^{-1} \be_s(k) \bigr) \;  \sdet(s^{-1}) \; dg \times ds \; ,
\end{align*}
where we have used twice the Fubini theorem and the invariance of the Haar measure on $G_2$ to
obtain the last equality. The Fubini theorem can be applied because the integrand has compact
support in $G_1 \times \cO_k$.

Then the same equality is valid for all positive Borel functions $h$ on $G$.
Taking $h$ of the form $(h_1 \ot h_2) \rho^{-1}$ we get the result.
\end{proof}
\begin{lemma}\label{46}
Suppose that $f$ is a bounded Borel function on $G_1 \times G_2$ whose support has compact closure.
Suppose that $\xi \in K(G_1 \times G_2)$. Define, for nearly all $(k,s) \in
G_1 \times G_2$
$$(\pi(f)\xi)(k,s) = \int \bar{\cU}(g,g^{-1}k,s) \; f(g, \al_{g^{-1}k}(s)) \; \xi(g^{-1}k,s) \;
dg \; .$$
Then $\pi(f)$ defines a bounded operator on $L^2(G_1 \times G_2)$ belonging to $\Nfi \subset M$ and such that $\la(\pi(f)) = f$.
\end{lemma}
Observe that $\pi(f)\xi$ is a bounded Borel function having a support with compact closure.
\begin{proof}
If $|f(g,s)| \leq M$ for all $g$ and $s$ and $f$ has its support in $L_1 \times L_2$ for
certain compact sets $L_1$ and $L_2$, then $\pi(f)\xi \in L^2(G_1 \times G_2)$ and $\|
\pi(f)\xi \| \leq M \mu_1(L_1) \|\xi\|$, where $\mu_1$ is the Haar measure on $G_1$. Hence
$\pi(f)$ defines a bounded operator. If now $f = F \ot G$ with $F \in K(G_1)$ and $G \in
K(G_2)$, then $$\pi(f) = (\om_F \ot \io \ot \io)(\Wtil) \al(G)$$ and the result follows in
this special case. When $f$ is continuous on $L_1 \times L_2$ and equals $0$ elsewhere we can
find a sequence $f_n$ of linear combinations of functions $F \ot G$ where $F$ is continuous on
$L_1$ and $G$ is continuous on $L_2$, converging uniformly to $f$. It is clear from the
inequality above that $\pi(f_n)$ converges in norm to $\pi(f)$ and that $f_n$ converges in
$L^2(G_1 \times G_2)$ to $f$. Hence $\pi(f) \in \Nfi$ and $\la(\pi(f)) = f$, because $\la$ is
closed. Finally, if $f$ is Borel, if its support is contained in $L_1 \times L_2$ and if
$|f(g,s)| \leq M$ for all $g$ and $s$, we can use the Lusin theorem to obtain a sequence $f_n$
of functions continuous on $L_1 \times L_2$ and equal to $0$ elsewhere, satisfying $|f_n(g,s)|
\leq M$ for all $g,s$ and $n$, and such that $$(\mu_1 \times \mu_2) \bigl( \{ (g,s) \mid
f_n(g,s) \neq f(g,s) \} \bigr) \recht 0,$$ where $\mu_i$ is the Haar measure on $G_i\
(i=1,2)$. Using Lemma~\ref{45} one can check that $\pi(f_n) \recht \pi(f)$ in the \strong
topology and $f_n \recht f$ in $L^2(G_1 \times G_2)$. Because $\la$ is \strong -- norm closed,
the result follows.
\end{proof}
\begin{lemma}\label{47}
There exist unique elements $\cUtil, \cVtil \in L^\infty(G_1) \ot L^\infty(G_2)$ such that,
for almost all $(g,k,s) \in G_1 \times G_1 \times G_2$, $$(\io \ot \al)(\cUtil)(g,k,s) =
\cU(g,g^{-1}k,s) \; \cU(g^{-1},k,s)$$ and such that, for almost all $(g,s,t) \in G_1 \times
G_2 \times G_2$, $$(\be \ot \io)(\cVtil)(g,s,t) = \cV(g,s,t^{-1}) \; \cV(g,t^{-1}s,t).$$ Also
we have $\cUtil(g,\al_g(s)) = \cUtil(g^{-1},s)$ and $\cVtil(\be_s(g),s) = \cVtil(g,s^{-1})$
for almost all $(g,s) \in G_1 \times G_2$.
\end{lemma}
\begin{proof}
Define a unitary $Z \in L^\infty(G_1 \times G_1 \times G_2)$ by $$Z(g,k,s) = \cU(g,g^{-1}k,s)
\; \cU(g^{-1},k,s) \quad\text{for almost all}\; (g,k,s) \in G_1 \times G_1 \times G_2 \; .$$
We claim that $Z \in L^{\infty}(G_1) \ot \al(L^\infty(G_2))$. In view of Corollary~\ref{119}
and Theorem~\ref{110} it suffices to prove that $$(1 \ot V_1^* \ot 1) \; (\io \ot \io \ot
\al)(Z) \; (1 \ot V_1 \ot 1) = Z_{124} \; ,$$ which amounts to prove, for almost all
$(g,h,k,s)$, the equality $$Z(g,h,\al_k(s)) = Z(g,hk,s).$$ But applying \eqref{cocU} to
$(g,g^{-1}h,k,s)$ and to $(g^{-1},h,k,s)$ we obtain this last formula.

Now we can take $\cUtil = (\io \ot \al^{-1})(Z)$. The second statement can
be proved similarly. Finally, to prove that $\cUtil(g,\al_g(s)) = \cUtil(g^{-1},s)$
for almost all $(g,s)$ it suffices to check that $\cUtil(g,\al_{gh}(s)) =
\cUtil(g^{-1},\al_h(s))$ for almost all $(g,h,s)$,  which is clear. Similarly we get $\cVtil(\be_s(g),s) = \cVtil(g,s^{-1})$ for almost all $(g,s)$.
\end{proof}
Then we can finally prove the main result.
\begin{proposition}\label{48}
Denote by $X=J \nab^{1/2}$ the modular operator of the left invariant weight $\vfi$ on the
l.c.\ quantum group $(M,\de)$ with GNS-construction $(L^2(G_1 \times G_2),\io,\la)$ and by
$\Xh = \Jh \nabh^{1/2}$ the modular operator of the left invariant weight $\vfih$ on the the
dual l.c.\ quantum group $(\Mh,\deh)$ with GNS-construction $(L^2(G_1 \times G_2),\io,\lah)$.
Then, for every bounded measurable function $\xi$ on $G$ whose support is contained in a
compact subset of $\cO$, we get that $\xi \in \cD(X)$, $\xi \in \cD(\Xh)$, and
\begin{align*}
(X \xi)(g,s) &= \cUtil(g^{-1},s) \; \sdeo(g^{-1}) \; \bar{\xi}(g^{-1},\al_g(s))
\\ (\Xh \xi)(g,s) &= \cVtil(g,s^{-1}) \; \sdet(s^{-1}) \; \bar{\xi}(\be_s(g),s^{-1}) \\
(J \xi)(g,s) &= \sde \bigl( i(\be_s(g)) \bigr)^{-1/2} \; \sdeo \bigl( \be_s(g) g^{-1} \bigr)^{1/2}
\; \sdet\bigl( \al_g(s) s^{-1} \bigr)^{1/2} \; \cUtil(g^{-1},s) \; \bar{\xi}(g^{-1},\al_g(s))
\\
(\Jh \xi)(g,s) &= \sde \bigl( j(\al_g(s)) \bigr)^{1/2} \; \sdeo(\be_s(g) g^{-1})^{1/2} \;
\sdet(\al_g(s) s^{-1})^{1/2} \; \cVtil(g,s^{-1}) \; \bar{\xi}(\be_s(g),s^{-1}) \; .
\end{align*}
Moreover these elements $\xi$ form a core for both $X$ and $\Xh$.

The operators $\nab$ and $\nabh$ are strictly positive multiplication operators with the
functions
\begin{align*}
\nab(g,s) &= \sde\bigl( i(\be_s(g)) \bigr)^{-1} \; \sdeo \bigl( \be_s(g) g \bigr) \; \sdet
\bigl( \al_g(s) s^{-1} \bigr)   \\ \nabh(g,s) &= \sde \bigl(
j(\al_g(s)) \bigr) \; \sdeo\bigl( \be_s(g) g^{-1} \bigr) \; \sdet \bigl( \al_g(s) s \bigr) \; .
\end{align*}

The scaling constant of $(M,\de)$ equals $1$.
\end{proposition}
\begin{proof}
Let $\cF$ be the set of bounded Borel functions on $G_1 \times G_2$ whose support is
contained in a compact subset of $\cO$ and let $f \in \cF$.
Using Lemmas~\ref{46} and \ref{47} one can check that for all $\xi \in K(G_1 \times G_2)$
\begin{align*}
(\pi(f)^* \xi)(k,s) &= \int \cU(g,k,s) \; \bar{f}(g,\al_k(s)) \; \xi(gk,s)  \; dg \\
&= \int \bar{\cU}(g,g^{-1}k,s) \; \sdeo(g^{-1}) \; \cU(g,g^{-1}k,s) \; \cU(g^{-1},k,s) \;
\bar{f}(g^{-1},\al_k(s)) \; \xi(g^{-1}k,s) \; dg \\
&=(\pi(h) \xi)(k,s)
\end{align*}
where $h(g,s) = \sdeo(g^{-1}) \; \cUtil(g^{-1},s) \; \bar{f}(g^{-1},\al_g(s))$. Because $(g,s)
\mapsto (g^{-1},\al_g(s))$ is a homeomorphism of $\cO$ we get $h \in \cF$. Then it follows that $f \in \cD(X)$ and $Xf=h$. When both $f_1$ and $f_2$ belong to $\cF$,
one can check that $\pi(f_1)\pi(f_2)=\pi(f_3)$ where
$$f_3(h,s) = \int \bar{\cU}(g,g^{-1}h,s) \; f_1(g,\al_{g^{-1}h}(s)) \; f_2(g^{-1}h,s) \; dg \;
.$$
We see that $f_3 \in \cF$. From Definition~\ref{27} and Lemma~\ref{46} it follows that $\{
\pi(f) \mid f \in \cF\}$ is a \strong -- norm core for $\la$. The computations
above show that it is also a $*$-algebra, so $\cF \subset L^2(G_1 \times G_2)$ is a core for $X$.
Let now $f$ and $h$ be in $\cF$. Then we get, using Lemma~\ref{45} and the notation
$F$ introduced there,
\begin{align*}
\langle Xf ,h \rangle &= \int \sdeo(g^{-1}) \; \cUtil(g^{-1},s) \; \bar{f}(g^{-1},\al_g(s)) \;
\bar{h}(g,s) \; dg \times ds \\
&= \int \Bigl( \int \sdeo(g^{-1}) \; \cUtil(g^{-1},s) \; \bar{f}(g^{-1},\al_g(s)) \;
\bar{h}(g,s) \; ds \Bigr) \; dg \\
&= \int \Bigl( \int F(g^{-1},s) \; \sdeo(g^{-1}) \; \cUtil(g^{-1},\al_{g^{-1}}(s)) \; \bar{f}(g^{-1},s) \;
\bar{h}(g, \al_{g^{-1}}(s)) \; ds \Bigr) \; dg \\
&= \int F(g,s) \; \cUtil(g^{-1},s) \; \bar{f}(g,s) \; \bar{h}(g^{-1},\al_g(s)) \; dg \times ds
\; .
\end{align*}
In this computation we used twice the Fubini theorem on a bounded Borel function whose support
has compact closure. Since the functions $f$ form a core for $X$, the computation above shows
that $h \in \cD(X^*)$ and $$(X^* h)(g,s)= \cUtil(g^{-1},s) \; F(g,s) \;
\bar{h}(g^{-1},\al_g(s)) \; .$$ Then for all $f \in \cF$ we get $f \in \cD(X^* X)$ and, with
$\nab(g,s)$ as in the statement of the proposition, $$(X^* X f)(g,s) = \nab(g,s) \; f(g,s).$$
Denote by $A$ the strictly positive multiplication operator with the function $(g,s) \mapsto
\nab(g,s)$. Because $X^* X = \nab$ we see that $A f = \nab f$ for all $f \in \cF$. Because the
subspace $\cF$ of $L^2(G_1 \times G_2)$ is invariant under $A^{it}$ for all $t \in \R$, it is
a core for $A$. Hence it follows that $A \subset \nab$, and so $A=\nab$ because both are
self-adjoint. From the formulas for $X$ and $\nab$ we get the formula for $J$.

Next Proposition~\ref{29} and Proposition~\ref{JJhat} give that $(\Mh,\deh)$ is also a cocycle
bicrossed product. If we define $$\itil : G_2 \recht G : \itil(s) = j(s)^{-1}
\qquad\text{and}\qquad \jtil : G_1 \recht G : \jtil(g) = i(g)^{-1}$$ then $G_2$ and $G_1$
become a matched pair of groups and the maps $\tilde{\al}$ and $\tilde{\be}$ are given by
$$\tilde{\al}_s(g) = \be_s(g) \qquad\text{and}\qquad \tilde{\be}_g(s) = \al_g(s) \; .$$ The
maps $(s,t,g) \mapsto \cV(g,t,s)$ and $(s,g,h) \mapsto \cU(h,g,s)$ satisfy \eqref{cocU}. Hence
we can make the cocycle bicrossed product, and Proposition~\ref{29} shows that this equals, up
to a flip map, $(\Mh,\deh)$. Then the formulas for $\hat{X},\Jh$ and $\nabh$ follow from the
formulas for $X$, $J$ and $\nab$.

Finally, we get immediately that $\nab$ and $\nabh$ commute strongly. Hence it follows from \cite{KV2}, Proposition~2.13 that the scaling constant equals $1$.
\end{proof}
Now one could write explicit formulas for $\tau_t(\pi(f))$ and $R(\pi(f))$ where $f$ is a
bounded Borel function whose support has compact closure. Let us give concrete formulas for
the remaining operators related to $(M,\de)$ and its dual.
\begin{proposition} \label{49}
Let $P$ be the unitary implementation of $\tau$ \cite{KV1}, Definition 6.9, $\sde_M$ the
modular element of $(M,\de)$ and $\sde_{\Mh}$ the modular element of $(\Mh,\deh)$. Then $P$ is
the multiplication operator with the function
\begin{align*}
P(g,s) &= \sde\bigl( i(g \be_s(g)^{-1}) \bigr) \; \sdeo\bigl( g^{-1} \be_s(g) \bigr) \;
\sdet\bigl( s^{-1} \al_g(s) \bigr) \; ,\\
\intertext{$\sde_M$ is the multiplication operator with the function}
\sde_M(g,s) &= \sde\bigl( j(\al_g(s)) \bigr)^{-1} \; \sdet\bigl( \al_g(s) \bigr)^{-2} \\
\intertext{and $\sde_{\Mh}$ is the multiplication operator with the function}
\sde_{\Mh}(g,s) &= \sde\bigl( i(\be_s(g)) \bigr) \; \sdeo\bigl( \be_s(g) \bigr)^{-2} \; .
\end{align*}
\end{proposition}
\begin{proof}
Using \cite{KV1}, Definition 6.9 and Proposition 8.23 we get that $\tau_t(x) = P^{it} x P^{-it}$
and $\tauh_t(y) = P^{it} y P^{-it}$ for all $t \in \R$, $x \in M$ and $y \in \Mh$. From
\cite{KV2}, Proposition 2.13 it follows that $JPJ = P^{-1}$. If $A$ is the
multiplication operator with the function $(g,s)\mapsto P(g,s)$ defined in the proposition, one can compute that the closure of $A \nab^{-1}$ is the multiplication operator with the
function
$$(g,s) \mapsto \sde(i(g)) \; \sdeo(g)^{-2}$$
which commutes with $\Mh$. Hence $A^{it} y A^{-it} = \tauh_t(y)$ for all $t \in \R$ and $y \in
\Mh$. Similarly we prove that $A^{it}$ implements $\tau_t$ for all $t \in \R$.
Because $M$ and $\Mh$ generate $B(L^2(G_1 \times G_2))$, the operators $A$ and $P$ are
proportional. Because $JAJ = A^{-1}$, it follows that $A=P$.

From \cite{KV1}, Lemma 8.8 it follows that $\nabh$ is the closure of $P J \sde_M^{-1} J$,
from where the formula for $\sde_M$ follows. The formula for $\sde_{\Mh}$ can be obtained in the same way.
\end{proof}
Lemma~\ref{45} shows that the transformation $\al_g$ of $G_2$ preserves Borel sets
of measure zero, and we get the following concrete formula for the Radon-Nikodym derivative.
$$\frac{d \al_g(s)}{ds} = \sde\bigl( i(\be_s(g)^{-1}) \bigr) \; \sde_1\bigl(\be_s(g) \bigr) \;
\sde_2\bigl( \al_g(s) s^{-1} \bigr) \; .$$
Following \cite{Maj2} we use the notation
$$\hchi(g,s) := \frac{d \al_g(s)}{ds} \qquad\text{and}\qquad \Psi(s,g) := \frac{d \be_s(g)}
{dg} \; .$$
Then we see that
\begin{equation} \label{formulachi}
\frac{\hchi(g,s)}{\hchi(g,e)} = \sde\bigl( i(g \be_s(g)^{-1}) \bigr) \; \sde_1\bigl(g^{-1} \be_s(g) \bigr) \;
\sde_2\bigl( \al_g(s) s^{-1} \bigr) \; .
\end{equation}
By symmetry, we obtain that
\begin{equation} \label{formulachi2}
\frac{\Psi(s,g)}{\Psi(s,e)} = \sde\bigl( j(s^{-1} \al_g(s)) \bigr) \; \sde_1\bigl(g^{-1} \be_s(g) \bigr) \;
\sde_2\bigl( \al_g(s) s^{-1} \bigr) \; .
\end{equation}
Because $i(g) j(s) =  j \bigl( \al_g(s) \bigr) \; i\bigl(\be_s(g)\bigr)$ we get that
$$\sde\bigl( i(g \be_s(g)^{-1}) \bigr) = \sde\bigl( j(s^{-1} \al_g(s)) \bigr) \; .$$ From this
it follows that $$\frac{\hchi(g,s)}{\hchi(g,e)} = \frac{\Psi(s,g)}{\Psi(s,e)}$$ and using the
S.~Majid's notation, we put $$\xi(g,s) := \frac{\hchi(g,s)}{\hchi(g,e)} \; .$$ The following
proposition characterizes all cocycle bicrossed products of l.c.\ groups which are Kac
algebras. Its final statement generalizes \cite{Maj2}, Theorem 2.12.
\begin{proposition} \label{charKac}
The scaling group $\tau$ of the l.c.\ quantum group $(M,\de)$ is trivial if and only if
$$\xi(g,s) = 1 \tekst{for all} (g,s) \in \cO \; .$$ The modular element $\sde_M$ is affiliated
to the centre of $M$ if and only if $$\sde\bigl( j(s^{-1} \al_g(s)) \bigr) \; \sde_2\bigl(
s^{-1} \al_g(s) \bigr)^2 = 1 \tekst{for all} (g,s) \in \cO \; .$$ The l.c.\ quantum group
$(M,\de)$ is a Kac algebra if and only if $$\xi(g,s) = 1 \quad\text{and}\quad \frac{\sde_1
\bigl( \be_s(g) \bigr)}{\sde_1(g)} = \frac{\sde_2 \bigl( \al_g(s) \bigr)}{\sde_2(s)} \; .$$
\end{proposition}
\begin{proof}
$\tau$ is trivial if and only if $P=1$ or if and only if $\xi(g,s)=1$ almost
everywhere. The continuity of $\xi$ on $\cO$ gives the first statement.

Proposition~\ref{49} shows that $\sde_M =\al(H)$ where the function $H$ on $G_2$ is defined by
$H(s) = \sde(j(s))^{-1} \; \sde_2(s)^{-2}$. Then $\sde_M$ is affiliated to the center of $M$
if and only if $\al(H) = 1 \ot H$ (Proposition \ref{center}), and the second statement
follows.

To prove the final statement, observe that $(M,\de)$ is a Kac
algebra if and only if $\tau$ is trivial and $\sde_M$ is affiliated to the center of $M$.
If $(M,\de)$ is a Kac algebra, then by the previous statements, $\xi(g,s)=1$ for all $(g,s) \in \cO$ and
$$\sde\bigl( j(s^{-1} \al_g(s)) \bigr) \; \sde_2\bigl( s^{-1} \al_g(s) \bigr)^2 = 1 \tekst{for
all} (g,s) \in \cO \; .$$
Dividing this equation by the formula for $\xi=1$ \eqref{formulachi2} we get
$$\sde_1\bigl(g \be_s(g)^{-1}\bigr) \; \sde_2\bigl( s^{-1} \al_g(s) \bigr) = 1 \tekst{for
all} (g,s) \in \cO \; .$$
Rewriting this we get the equality
$$\frac{\sde_1 \bigl( \be_s(g) \bigr)}{\sde_1(g)} = \frac{\sde_2
\bigl( \al_g(s) \bigr)}{\sde_2(s)} \; .$$
Clearly we can also go the other way around to complete the proof.
\end{proof}
\begin{corollary}\label{410}
If $\al$ or $\be$ is trivial, $(M,\de)$ and
$(\Mh,\deh)$ are Kac algebras.
\end{corollary}
\begin{proof}
Let us suppose that $\al$ is trivial. Then for all $(g,s) \in \cO$ we have $j(s) i(\be_s(g)) =
i(g)j(s)$. Hence we get that $\sde\bigl(i(g^{-1} \be_s(g)) \bigr) = 1$ for all $(g,s) \in
\cO$. Because $\al$ is trivial, we also have
$$\frac{d \al_g(s)}{d s} = 1 \tekst{for all} (g,s) \in \cO \; .$$
Applying this equation to $(g,s)$ and $(g,e)$ we get
$\xi(g,s) = 1$ for all $(g,s) \in \cO$. Combining \eqref{formulachi} with the equality $\sde\bigl(i(g^{-1} \be_s(g)) \bigr) = 1$
we get that $\sde_1\bigl(g^{-1} \be_s(g) \bigr) = 1$ for all
$(g,s) \in \cO$. Because $\al$ is trivial this gives the equation
$$\frac{\sde_1 \bigl( \be_s(g) \bigr)}{\sde_1(g)} = \frac{\sde_2
\bigl( \al_g(s) \bigr)}{\sde_2(s)} \; .$$
Using the previous proposition we obtain that $(M,\de)$ is a Kac algebra.
\end{proof}
\begin{corollary}
If both $\al$ and $\be$ preserve modular functions and Haar measures, then $(M,\de)$ and $(\Mh,\deh)$ are Kac algebras.
\end{corollary}
\begin{proof}
Because $\al$ preserves the Haar measure of $G_2$ we get $\xi(g,s) = 1$ for all $(g,s) \in
\cO$. Because $\al$ and $\be$ preserve the modular functions, we get
$$\frac{\sde_1 \bigl( \be_s(g) \bigr)}{\sde_1(g)} =1= \frac{\sde_2
\bigl( \al_g(s) \bigr)}{\sde_2(s)}$$
for all $(g,s) \in \cO$. Then we can apply Proposition~\ref{charKac}.
\end{proof}
Remark that the conditions of this corollary are fulfilled if both groups are discrete.
Indeed, any discrete group is unimodular and the Haar measure is constant at an arbitrary
point of such a group.
\begin{corollary}
If $(G_1,G_2)$ is a fixed matched pair of groups and if cocycles $\cU$ and $\cV$ satisfy
\eqref{cocU}, we get a cocycle bicrossed product $(M,\de)$. If one of these cocycle bicrossed
products is a Kac algebra, then all of them are Kac algebras.
\end{corollary}
\begin{proof}
The necessary and sufficient conditions for $(M,\de)$ to be a Kac algebra given in
Proposition~\ref{charKac} are independent of $\cU$ and $\cV$.
\end{proof}
\subsection{The group of extensions}\label{secgroup}
If $(G_1,G_2)$ is a fixed matched pair of l.c.\ groups and if $\cU$ and $\cV$ are measurable
functions satisfying \eqref{cocU}, we get by the previous subsection a cocycle matching
$(\tau,\cU,\cV)$ of $(L^\infty(G_1),\deo)$ and $(L^\infty(G_2),\det)$, and hence, by
Proposition~\ref{34}, an extension $$(L^\infty(G_2),\det) \overset{\al}{\lrecht} (M,\de)
\overset{\be}{\lrecht} (\cL(G_1),\deoh).$$ We will say that such an extension $(M,\de)$ is
associated with the matched pair $(G_1,G_2)$. Proposition~\ref{35ter} and the remark following
its statement show that the extensions given by two cocycle matchings $(\tau,\cUa,\cVa)$ and
$(\tau,\cUb,\cVb)$, are isomorphic if and only if there exists a measurable map $\cR$ from
$G_1 \times G_2$ to $U(1)$, satisfying
\begin{align*}
\cUb(g,h,s) &= \cUa(g,h,s) \; \cR(h,s) \; \cR(g,\al_h(s)) \; \bar{\cR}(gh,s) \\
\cVb(g,s,t) &= \cVa(g,s,t) \; \cR(g,s) \; \cR(\be_s(g),t) \; \bar{\cR}(g,ts)
\end{align*}
almost everywhere. If this is the case, the pairs $(\cUa,\cVa)$ and $(\cUb,\cVb)$ will be
called cohomologous.
Then the set of equivalence classes of cohomologous pairs of cocycles $(\cU,\cV)$ satisfying \eqref{cocU}, exactly corresponds to the set $\Gamma$ of classes of isomorphic extensions associated with $(G_1,G_2)$.

The set $\Gamma$ can be given the structure of an abelian group by defining
$$\pi(\cUa,\cVa) \; \cdot \; \pi(\cUb,\cVb) = \pi(\cUa \cUb,\cVa \cVb)$$
where $\pi(\cU,\cV)$ denotes the equivalence class containing the pair $(\cU,\cV)$.
\begin{terminology} \label{410bis}
The group $\Gamma$ is called the group of extensions of $(L^\infty(G_2),\det)$ by
$(\cL(G_1),\deoh)$ associated with the matched pair of l.c.\ groups $(G_1,G_2)$. The unit of
this group corresponds to the equivalence class of cocycles cohomologous to the trivial
cocycles. The corresponding extension is called the split extension; all other extensions are
called non-trivial extensions.
\end{terminology}
\subsection{Cocycles in the sense of S. Baaj and G. Skandalis}
Consider the split extension $(M_0,\de_0)$ of $(L^\infty(G_2),\det)$ by $(\cL(G_1),\deoh)$
associated with the matched pair of l.c.\ groups $(G_1,G_2)$. The associated (dual)
multiplicative unitary, acting on $L^2(G_1 \times G_2 \times G_1 \times G_2)$, is given by
$$\Wh_0 = ((\be \ot \io)(W_1) \ot 1) \; (1 \ot (\io \ot \al)(\Wh_2)) \; .$$ Similarly to
\cite{B-S1}, D\'efinition 8.24 we call a function $\Theta : G_1 \times G_2 \times G_1 \times
G_2 \recht U(1)$ a $\Wh_0$-cocycle, if $\Theta$ is measurable and the unitary $\Theta \Wh_0$
is multiplicative. Here we consider $\Theta$ as a unitary in $L^{\infty}(G_1 \times G_2 \times
G_1 \times G_2)$.

Let us discuss the link between $\Wh_0$-cocycles $\Theta$ and cocycles $\cU$ and $\cV$.
Let first $\cU$ and $\cV$ be cocycles in the sense of Lemma~\ref{44}. Then $(\tau,\cU,\cV)$ is a cocycle matching of $(L^\infty(G_1),\deo)$ and
$(L^\infty(G_2),\det)$ and hence we obtain the multiplicative unitary
$$\Wh = (\al \ot \io \ot \io)((W_1 \ot 1) \cU^*) \; (\io \ot \io \ot \be)(\cV (1 \ot \Wh_2)) \;.$$
One can check that for $\xi \in L^2(G_1 \times G_2 \times G_1 \times G_2)$ we have
\begin{align*}
(\Wh_0 \xi)(g,s,h,t) &=\xi(g,\al_{\be_s(g)^{-1}h}(t)s,\be_s(g)^{-1}h, t) \\
(\Wh \xi)(g,s,h,t) &= \bar{\cU}(\be_s(g),\be_s(g)^{-1} h, t) \;
\cV(g,s,\al_{\be_s(g)^{-1}h}(t)) \;
\xi(g,\al_{\be_s(g)^{-1}h}(t)s,\be_s(g)^{-1}h, t)
\end{align*}
nearly everywhere. Hence, defining
\begin{equation}\label{Theta}
\Th(g,s,h,t)= \bar{\cU}(\be_s(g),\be_s(g)^{-1} h, t) \;
\cV (g,s,\al_{\be_s(g)^{-1}h}(t))
\end{equation}
we get that $\Th$ is a $\Wh_0$-cocycle.

Vice versa, if $\Th$ is a $\Wh_0$-cocycle, the
multiplicativity of $\Th \Wh_0$ amounts, almost everywhere, to the equation
\begin{multline} \label{mu}
\Th(g,s,h,t) \; \Th(g, \al_{g^{-1}}(\al_h(t) \al_g(s) ),k,r) \; \Th( \be_s(g)^{-1}h,t,
\be_{\al_h(t) \al_g(s)}(g^{-1})k,r) \\
= \Th(h,t,k,r) \; \Th(g,s,h,\al_{\be_t(h)^{-1}k}(r)t).
\end{multline}

If we know that the previous equality holds whenever its arguments make sense (by this we mean
whenever $(h,t)$, $(g,s)$, $(g^{-1},\al_h(t) \al_g(s)), \ldots \in \cO$), then we can define
$$\cU(g,h,t) = \bar{\Th}(g,e,gh,t) \; \Th(g,e,g,e) \qquad\text{and}\qquad \cV(g,s,t) =
\Th(g,s,\be_s(g),t) \; \bar{\Th}(e,e,e,t) \; .$$ Some computations reveal that $\cU$ and
$\cV$ are cocycles, i.e.\ they satisfy \eqref{cocU}, and $\Theta$ is given by
Eq.~\eqref{Theta}. We cannot prove the same result when we only know that \eqref{mu}
holds almost everywhere. Nevertheless the formulas for the cocycles $\cU$ and $\cV$ given
above are sufficient to find these cocycles starting from a $\Wh_0$-cocycle $\Theta$. If both
$G_1$ and $G_2$ are discrete, the extra assumption is redundant because then 'nearly
everywhere' means everywhere.

\section{Examples of cocycle bicrossed products}
\subsection{A brief survey of known examples}
{\bf 1. The Kac algebra case: applications of Corollary \ref{410}.}

Let l.c.\ quantum groups $G_1$, $G_2$ form a matched pair and let one of the actions, say
$\beta$, be trivial. Then $G_2$ is a normal subgroup of $G$ (see the proof of  Corollary
\ref{410}) and $G_1$ acts on $G_2$ by (inner in $G$) automorphisms (this is seen from Lemma
\ref{42}). Now, as it was explained in \cite{Maj2}, the corresponding split extension is
nothing but a semidirect product of the commutative Kac algebra $L^{\infty}(G_2)$ with the
group $G_1$ acting on it (the construction of a semidirect product of a Kac algebra with a
l.c. group is due to J.~De Canniere \cite{Can}). In the above situation, there is also the
induced action of $G_1$ on $\cL(G_2)$, and the dual Kac algebra of the above split extension
can be again computed following J.~De Canniere (for a detailed explanation see \cite{Maj2}).

Turning now to cocycles, the second equation \eqref{cocU} shows that $\cV(g,t,s)$ can be viewed as a function $\eta^g(t,s)$ on $G_1$ with values in the
group $Z_2(G_2;U(1))$ of $U(1)$-valued 2-cocycles on $G_2$.
Let us try to find a non-trivial solution for
\eqref{cocU} under the hypothesis that $\cU(\cdot)=1$. Then
the third equation \eqref{cocU} gives
$$
\eta^{gh}(t,s)=h\cdot\eta^g(t,s) \; \; \eta^h(t,s),
$$
where $h\cdot\eta^g(t,s):=\eta^g(\alpha_h(t),\alpha_h(s))$,
which means, by definition, that $g\mapsto \eta^g$ is an equivariant map from $G_1^{\text{op}}$ to $Z_2(G_2;U(1))$.

Such an equivariant map can be found at least in some special cases when both $G_1$ and $G_2$
are either finite abelian groups or finite-dimensional real vector spaces. Indeed, in
\cite{B-S1}, 8.26 the Kac-Paljutkin example \cite{KP1}, where $G_1=\Z/2\Z$ acts on
$G_2=\Z/2\Z\times \Z/2\Z$ by permutations, has been treated this way. It turns out that this
way we find the unique non-trivial extension, so the group of extensions here is $\Z/2\Z$
(\cite{Kac},\cite{Mas1}). More generally, if $G_1=\Z/2\Z$ acts on $G_2=\Z/m\Z\times \Z/m\Z$
($m\geq 2$ is natural) by permutations, the group of extensions is $\Z/m\Z$ (\cite{Kac,Mas1}),
concrete examples of non-trivial extensions can be found in \cite{Mas1,Sek,Vai}. Let us
mention some other matched pairs of finite groups that fit into this framework and for which
the group of extensions is non-trivial: a) $G_1=\Z/m\Z$ ($m\geq 2$ is natural) acts on
$G_2=\Z/m\Z\times \Z/m\Z$ as follows: $x\cdot a^ib^j=a^{i+j}b^j$ ($x,a,b$ are the
corresponding generators, $i,j=0,1,...,m-1$) (see \cite{Mas1}); b) $G_1=\Z/2\Z$ acts on
$G_2=\Z/2m\Z$ ($m\geq 2$ is natural) by inversion, then $G$ is the dihedral group (see
\cite{Nik,Vai}).

In \cite{B-S1}, 8.26, concerning the other Kac-Paljutkin example  \cite{KP2} with $G_1=\R$
acting on $\R^2$ by $\alpha_g(x)=exp(gK)(x)$ (where $K$ is a real $2\times 2$ matrix), the
above equivariant map $g\mapsto \eta^g$ was computed as well as for the Rieffel's example
\cite{Rief}, where $G_1=\R^2$ acts on $G_2=\R^3$ by $(a,b)\cdot(x,y,z)=(x+ay+bz,y,z)$. So, in
both cases the group of extensions is non-trivial. The same is true for some other matched
pairs where both $G_1$ and $G_2$ are finite-dimensional real vector spaces: a) $G_1=\R^n$
acting on $G_2=\R^{n+1}$ by $\alpha_a(b,t)=(b,t+ \langle a,b \rangle)$, where $a,b\in \R^n,
t\in \R,n\geq 2$ is natural (see \cite{E-V},\cite{Vai}); b) $G_1=\R^n$ acting on $G_2=\R^{m}$
by certain family of linear operators \cite{VanDaele}.

Another situation when non-trivial solutions of \eqref{cocU} can be easily found, is the one
with both actions trivial. Then $\cU(\cdot)$ can be choosen as a multiplicative function
$\zeta:G_2\to Z_2(G_1;U(1))$ (i.e., $\zeta^{ts}=\zeta^t \zeta^s$) and $\cV(\cdot)$ as a
multiplicative function $\eta:G_1\to Z_2(G_2;U(1))$ (i.e., $\eta^{gh}=\eta^g \eta^h$). If both
$G_1$ and $G_2$ are finite, the whole group of extensions has been computed in \cite{Kac},
3.5.

{\bf 2. The Kac algebra case: application of Corollary 4.18.}

a) Let $G=S_n$ be the symmetric group of order $n\geq 4, G_1= S_{n-1}$ its subgroup which fixes the point $n$ and $G_2=\Z/n\Z$ the subgroup generated by the cyclic permutation. Then it is shown in \cite{B-S1}, 8.23b that $G_1,G_2$ form a
matched pair with two non-trivial actions (for instance, the action of $s\in G_1$ on $a^i\in G_2$ is given by $\al_s (a^i)=a^{s(i)}, i=0,1,...,n-1$). The corresponding group of extensions is trivial \cite{Mas1}, 4.1, and the split extension is clearly a Kac algebra.

b) Following \cite {Majbook}, 6.2.16, consider the pair $G_1=T_+(n,\R), G_2=T_-(n,\R)$, formed
by the groups of upper (resp., lower) triangular matrices over $\R$ with 1 on the diagonal and
equipped with the actions $$ \alpha_g(s)=1+(s-1) g^{\text{T}},\
\beta_s(g)=(1+\theta(s)(g^{-1}-1))^{-1}, $$ where $\theta$ is the composition of inversion and
matrix transposition and $^{\text{T}}$ is matrix transposition. These formulas slightly differ
from \cite{Majbook}, 6.2.16 because there another convention on a matched pair is adopted.
Being connected nilpotent groups, both $G_1,G_2$ are unimodular \cite{Helg}, I.1.6 (morover,
one can check that their Haar measures are in fact the product measures of the Lebesgue
measures of all non-diagonal entries $g_{ij}$), so the actions preserve modular functions.
Then one can check that for $i>j$ we have $$ (\al_g(s))_{ij}=s_{ij}+\sum_{j<k<i} s_{ik} \;
\te(g)_{kj}, $$ from where it is clear that the transformation $\al_g$ preserves the Haar
measure of $G_2$. So the corresponding split extension is a Kac algebra because of Corollary
4.18. We do not know if other extensions exist in this situation, but if they exist, they are
all Kac algebras (see Corollary 4.19).

{\bf 3. The S. Majid's example \cite{Maj1}}

Using the Iwasawa decomposition of a complex semisimple Lie group $G=KAN$ \cite{Helg}, I.5.1, one can take $G_1=K,G_2=AN$, where $G_1$
and $G_2$ are respectively the maximal compact and the maximal solvable subgroups of $G$. For instance, for $SL(2,\C)$ viewed as a real Lie group we have $G_1=SU(2,\C)$ ($g\in G_1$ is defined by a pair $x,y\in \C$ such that $|x|^2+|y|^2=1$) and $G_2=\{s\in \R^3 | s_3>-1\}$ with the operation
$st=s+(s_3+1)t$ (see \cite{Majbook}, 8.3.13). Here and in what follows $s_i$ and $e_i$ denote the $i-th$ coordinate of the vector $s$ and $i$-th coordinate vector in $\R^3, i=1,2,3$.
Then the mutual actions of $G_1,G_2$ are
\begin{align*}
\alpha_g(s) &=\text{Rot}_g(s-\frac{|s|^2}{2(s_3+1)} e_3)+
\frac{|s|^2}{2(s_3+1)}e_3,
\\
\beta_s(x,y) &=\frac{((s_3+1)x+(s_1+is_2) \bar{y},y)}{\sqrt{|(s_3+1)x+(s_1+is_2)
\bar{y}|^2+|y|^2}}
\end{align*}
(see \cite{Majbook}, 8.3.18), where $\text{Rot}$ is the usual 3-vector rotation associated to
an element of $SU(2,\C)$. Now we are able to construct the corresponding split extension, and
to decide if it is a Kac algebra, using Proposition \ref{charKac}. One can compute (see
\cite{Maj1}) that the function $\xi(g,s)$ from this Proposition differs from 1, so the split
extension is not a Kac algebra. Again, we do not know if there exist other extensions, but if
they exist, they are all not Kac algebras (see Corollary 4.19).

\subsection{On infinitesimal objects of l.c. quantum groups}

In order to give concrete examples of cocycle bicrossed products we have to take concrete matched pairs of groups $(G_1,G_2)$ and to solve Eq. \eqref{cocU} to find cocycles $\cU$ and $\cV$.
Then one can perform the cocycle bicrossed product construction to get a l.c. quantum group.

If $G_1$ and $G_2$ are Lie groups, we also want to describe, at least formally, a Hopf algebra
consisting of generators of the l.c. quantum group coming from the matched pair $(G_1,G_2)$
with cocycles $\cU$ and $\cV$. This infinitesimal Hopf algebra is, in fact, an algebraic
cocycle bicrossed product Hopf algebra.

We assume that the cocycles $\cU$ and $\cV$ satisfy \eqref{cocU} and we require additionally
$\cU(g,e,s),\cU(e,g,s)$ and $\cV(e,s,r)$ to be well defined measurable functions in $s\in G_2$ for all $g\in G_1$ (resp., in $s,r \in G_2$). These extra requirements will always be fulfilled in the concrete examples that we study below. Then one can choose in the class of cohomologous cocycles such representatives that $\cU(e,e,s)=1$ and, taking in the first equation \eqref{cocU} subsequently $h=k=e$ and $g=h=e$, we get $\cU(g,e,s) = \cU(e,g,s) = 1.$
Now taking in the third equation \eqref{cocU} $g=h=e$, we get $\cV(e,s,r) =1.$

Let $\de_1$ be the coproduct on $L^\infty(G_1)$, $\de_2$ the coproduct on $L^\infty(G_2)$ and
$(M,\de)$ the cocycle bicrossed product. Let $\gone$ be the Lie algebra of $G_1$ and $\cH$ its
universal enveloping algebra which is a Hopf algebra with the usual symmetric coproduct
$$\de_\cH(X) = X \ot 1 + 1 \ot X \hskip 1cm (\forall X \in \gone).$$ We consider $\cH$ as
acting on smooth functions on $G_1$ by left invariant differential operators (see \cite{Helg},
II.4), i.e., $$H_g[A] = H_e[A(g \cdot)]$$ for all $g \in G_1$, $H \in \cH$ and all smooth
functions $A$ on $G_1$. Using Sweedler notation, it is clear that, for $H,G \in \cH$ and $A$ a
smooth function on $G_1$, we have $$\sum H_{(1) \, e}[A] \; H_{(2) \, e}[B] = H_e[AB]$$ and
$$(HG)_e [A] = H_e[ g \mapsto G_g[A]] = H_e[g \mapsto G_e[ \de_1(A)(g,\cdot)]] = (H_e \ot
G_e)[\de_1(A)] \; .$$ To make the link with the von Neumann algebra setting, we define
formally, for $H \in \cH$, the unbounded functional $\om_H$ on $L^\infty(G_1)$ by $$\om_H(A) =
H_e[A] \; .$$ By the previous computation we get, for all $H,G \in \cH$ and $a,b \in
L^\infty(G_1)$, $$(\om_H \ot \om_G) \de_1 = \om_{HG} \quad\text{and}\quad \om_H(ab) = \sum
\om_{H_{(1)}}(a) \; \om_{H_{(2)}}(b).$$

Next suppose that we have some natural algebra $\cA$ of functions on $G_2$ such that
$(\cA,\de_2)$ is a Hopf algebra. We compute how to match the Hopf algebras $\cH$ and $\cA$, in
the sense of \cite[Section 6.3]{Majbook}. Then we will construct their algebraic cocycle
bicrossed product in the sense of S.~Majid.

An arbitrary element $H\ot A$ of the algebraic tensor product  $\cH \ot \cA$ can be considered formally as an unbounded operator affiliated to the von Neumann algebra $M$ by the identification
$$H \ot A \quad \longleftrightarrow \quad (\om_H \ot \io \ot \io)(\Wtil) \; \al(A)$$
where, as before, $\Wtil = (W_1 \ot 1) \cU^*$.

In order to compute the algebra structure that $M$ brings into the vector space $\cH \ot \cA$, let us observe that
\begin{equation} \label{eqWtil}
(\de_1 \ot \io \ot \io)(\Wtil) \; (\io \ot \io \ot \al)(\cU^*) = \Wtil_{134} \; \Wtil_{234}
\; .
\end{equation}
This can be verified by the following computation.
\begin{align*}
(\de_1 \ot \io \ot \io)(\Wtil) \; (\io \ot \io \ot \al)(\cU^*) &= (\de_1 \ot \io \ot
\io)(W_{1,12} \cU^*) \; (\io \ot \io \ot \al)(\cU^*) \\
&=W_{1,13} W_{1,23} \; (\io \ot \de_1 \ot \io)(\cU^*) \; \cU^*_{234} \\
&=W_{1,13} \cU^*_{134} \; W_{1,23} \cU^*_{234} = \Wtil_{134} \; \Wtil_{234}
\end{align*}
where we used the cocycle equation for $\cU$.
Taking $H,G \in \cH$ and $A,B \in \cA$ we get
\begin{align*}
(\om_H \ot \io \ot \io) & (\Wtil) \; \al(A) \; (\om_G \ot \io \ot \io)(\Wtil) \; \al(B) =
(\om_H \ot \om_G \ot \io \ot \io)(\Wtil_{134} \; \al(A)_{34} \; \Wtil_{234}) \; \al(B) \\
&= (\om_H \ot \om_G \ot \io \ot \io)(\Wtil_{134} \; \Wtil_{234} \; (\io \ot \al)\al(A)_{234})
\; \al(B) \\
&= (\om_H \ot \om_G \ot \io \ot \io) \bigl( (\de_1 \ot \io \ot \io)(\Wtil) \; (\io \ot \io \ot
\al)( \cU^* \al(A)_{23}) \bigr) \; \al(B) \\
&= \sum ((\om_{H_{(1)}} \ot \om_{G_{(1)}})\de_1 \ot \io \ot \io)(\Wtil) \; \al\bigl(
(\om_{H_{(2)}} \ot \om_{G_{(2)}} \ot \io)(\cU^*) \; (\om_{G_{(3)}} \ot \io)\al(A) \; B \bigr)
\end{align*}
where we used Eq.~\eqref{eqWtil} and that $\al$ is an action. Identifying again with $\cH \ot
\cA$, we get
\begin{equation} \label{verm}
(H \ot A)(G \ot B) = \sum H_{(1)} G_{(1)} \ot \hchi(H_{(2)} \ot G_{(2)}) \; (G_{(3)} \acts A)
\: B
\end{equation}
where
\begin{align*}
& \acts : \cH \ot \cA \recht \cA : H \acts A = (\om_H \ot \io) \al(A) \\
& \hchi : \cH \ot \cH \recht \cA : \hchi(H \ot G) = (\om_H \ot \om_G \ot \io)(\cU^*) \; .
\end{align*}
More specifically, we get that $(H \acts A)(r) = H_e[A(\al_{\ \! \pmb{\cdot} \ \!}(r))]$ and $\hchi(H \ot G)(r) =
(H_e \ot G_e)[\cU^*(\cdot,\cdot,r)]$ for all $r \in G_2$.
We see that, up to some left/right conventions, the algebra $\cA$ becomes this way a cocycle $\cH$-module algebra, and that $\cH \ot \cA$ is the cocycle cross product algebra (compare \eqref{verm} with the formulas appearing in \cite{Majbook}, Proposition 6.3.7).
We can go further and describe generators and relations for the algebra $\cH \ot \cA$.

One can check that for $X,Y \in \gone$ we have
$$(X \ot 1)(Y \ot 1) = XY \ot 1 + 1\ot \hchi(X \ot Y)$$
and for $X \in \gone$ and $A \in \cA$ we have
$$(X \ot 1)(1 \ot A) = X \ot A \quad\text{and}\quad (1 \ot A)(X \ot 1) = X \ot A + 1 \ot (X
\acts A) \; .$$
For this we use the relations $\hchi(X \ot 1) = \hchi(1 \ot X) = 0$ for all $X \in \gone$ and $\hchi(1 \ot 1) = 1$ which follow from $\cU(g,e,s)=\cU(e,g,s)=1$ for all $g \in G_1$ and $s\in G_2$.
So the algebra $\cH \ot \cA$ is the algebra generated by
$$\{ A \mid A \in \cA \} \quad\text{and}\quad \{X \mid X \in \gone\}$$
with commutation relations
\begin{align} \label{comrel}
& [A,B] = 0 \quad\text{when}\; A,B \in \cA \\ & [A,X] = X \acts A \quad\text{when}\; X \in
\gone, A \in \cA \notag \\ & [X,Y] = [X,Y]_{\gone} + \hchi(X \ot Y - Y \ot X) \notag
\quad\text{when}\; X,Y \in \gone
\end{align}
where $[X,Y]_{\gone}$ denotes the Lie bracket of the Lie algebra $\gone$.

Next we want to compute the coproduct on the algebra  $\cH \ot \cA$ inherited from the
quantum group $(M,\de)$. Using Propositions~\ref{24} and \ref{25} we get, for all $H \in \cH$
and $A \in \cA$,
\begin{align*}
\deop & \bigl( (\om_H \ot \io \ot \io)(\Wtil) \; \al(A) \bigr) \\ &= (\om_H \ot \io \ot \io \ot \io
\ot \io) \bigl( (\Wtil \ot 1 \ot 1) \; ((\io \ot \al)\be \ot \io \ot \io)(\Wtil) \; (\io \ot
\al \ot \al)(\cV) \bigr) \; (\al \ot \al)\deopt(A) \\
&= \sum \bigl( (\om_{H_{(1)}} \ot \io \ot \io)(\Wtil) \ot 1 \ot 1) \; (\al \ot \io \ot
\io) ((\om_{H_{(2)}} \ot \io)\be \ot \io \ot \io)(\Wtil) \\ & \qquad\qquad\qquad \; (\al \ot \al)
(\om_{H_{(3)}} \ot \io \ot \io)(\cV) \; \bigl(\al(A_{(2)}) \ot \al(A_{(1)}) \bigr) \; .
\end{align*}
Using the chain rule one sees that it is possible to define a linear map $$\beh: \cH \recht
\cH \ot \cA : \beh(H) = \sum H^{(\bar{1})} \ot H^{(\bar{2})}$$ such that for all smooth
functions $B$ on $G_1$, $r \in G_2$ and $H \in \cH$ we have $$\bigl( (\om_H \ot \io)\be(B)
\bigr)(r) = H_e[B(\be_r(\cdot))] = \sum H^{(\bar{1})}_e[B] \; H^{(\bar{2})}(r) \; .$$ Fix some
coordinates in the neighbourhood of the unit element $e$ of $G_1$ and let $X^i$ be the vector
field in this  neighbourhood which takes the derivative to the $i$-th coordinate. Let $Y^i \in
\gone$ be such that $Y^i$ and $X^i$ agree at $e$. If now $X \in \gone$ then $$\beh(X) = \sum_i
Y^i \ot \beh_i(X) \quad\text{where}\; \beh_i(X)(r) = X_e[\be^i_r(\cdot)]$$ where $\be^i_r(g)$
denotes the $i$-th coordinate of $\be_r(g)$. Also define $$\Psi : \cH \recht \cA \ot \cA :
\Psi(H) = \si(\om_H \ot \io \ot \io)(\cV)$$ where $\si$ is the flip. Observe that
$\Psi(H)(r,s) = H_e[\cV(\cdot,s,r)]$ for all $r,s \in G_2$. Using the notation $$\Psi(H) =
\sum \Psi(H)^{(1)} \ot \Psi(H)^{(2)}$$ and using Sweedler notation for $\de_2(A)$ we get
$$\de(H \ot A) = \sum \bigl( {H_{(2)}}^{(\bar{1})} \ot \Psi(H_{(3)})^{(1)} A_{(1)} \bigr) \ot
\bigl( H_{(1)} \ot {H_{(2)}}^{(\bar{2})} \Psi(H_{(3)})^{(2)} A_{(2)} \bigr) \; .$$ Comparing
this formula with the formula in \cite{Majbook}, Proposition 6.3.8, we get that this coproduct
turns $\cH \ot \cA$ into a cocycle cross coproduct coalgebra (again, up to some left/right
conventions).

In terms of the generators $A \in \cA$ and $X \in \gone$ for the algebra $\cH \ot A$ we get
more specifically
\begin{align} \label{comult}
& \de(A) = \de_2(A) \quad\text{for all}\; A \in \cA \\ & \de(X) = 1 \ot X + \hat{\be}(X) +
\Psi(X) \quad\text{for all}\; X \in \gone \; . \notag
\end{align}
To obtain these formulas we use that $\Psi(1) = 1 \ot 1$ which follows from $\cV(e,r,s)=1$ for all $r,s \in G_2$.

\begin{remark} \label{Liebialg} When $G_1, G_2$ is
a matched pair of Lie groups, it is helpful to look also at the corresponding matched pair of
Lie algebras $\fg_1,\fg_2$ (two Lie algebras form a matched pair if the direct sum of their
linear spaces is again a Lie algebra - see \cite{Majbook}, 8.3). On the other hand, a Lie
algebra (resp., dual vector space to a Lie algebra) can be viewed as a Lie bialgebra with zero
cobracket (resp., bracket). Recall that the notion of a Lie bialgebra is due to Drinfel'd
\cite{Drin}. Now for any pair of Lie bialgebras of the form $\fg_1,\fg^*_2$ (where
$\fg_1,\fg_2$ is a matched pair of Lie algebras) all the Lie bialgebra extensions, i.e., exact
sequences $0\to \fg^*_2\to \fh\to \fg_1\to 0$ of Lie bialgebras, are described in \cite{Mas2},
1.12 in terms of compatible actions of $\fg_1$ and $\fg_2$ on each other, as on vector spaces,
and a pair of compatible 2-cocycles for these actions. Here 2-cocycles are defined as linear
maps $\fg_1\bigwedge \fg_1\to \fg^*_2$ (resp., $\fg_2\bigwedge \fg_2\to \fg^*_1$) verifying
certain cocycle identities \cite{Majbook}, 2.3.

We consider this theory of Lie bialgebra extensions as an infinitesimal version of our theory
and, even if there is no any strict claim in this direction, we believe it helps to understand
better some concrete situations. For instance, if one of the Lie groups forming the matched
pair is $\R$, the corresponding cocycle on the Lie bialgebra level must be cohomologous to
zero, according to the above definition. From \cite{Mas2}, Theorem~4.11 it follows that the
same picture holds for Hopf algebras. This suggests that the situation is similar for l.c.\
quantum groups. If both $G_1=G_2=\R$, there should only exist the split extension.
\end{remark}
\subsection{The S.~Baaj and G.~Skandalis' example}
In this subsection we discuss an example which was already presented by S.~Baaj and
G.~Skandalis in a talk in Oberwolfach \cite{Skand}. We nevertheless include this example
explicitly, because we want to compute the associated infinitesimal Hopf algebra and show how
this example is a deformation of the $ax+b$-group.

Consider the group $$G=\{ (a,b) \mid a \in \R \setminus \{0\} ,b \in \R \} \quad\text{with}\;
(a,b)(c,d)=(ac, d+cb) \; .$$ Define $G_1 = G_2 = \R \setminus \{0\}$, with multiplication as
group operation. Then define $$i:G_1 \recht G : i(g) = (g,g-1) \quad\quad j: G_2 \recht G :
j(s) = (s,0) \; .$$ This way $(G_1,G_2)$ is a matched pair of groups in the sense of
Definition~\ref{41}. A direct computation then gives, for $g,s \in \R \setminus \{0\}$ and
$s(g-1)+1 \neq 0$, $$\al_g(s) = \frac{gs}{s(g-1)+1} \quad\quad \be_s(g) = s(g-1)+1.$$ Taking
$\cU=\cV=1$, we can construct the bicrossed product l.c. quantum group $(M,\de)$ having the
following nice property.
\begin{proposition}
The l.c. quantum group $(M,\de)$ is self-dual, i.e., $(\Mh,\deh)\cong (M,\de)$.
\end{proposition}
\begin{proof}
Proposition~\ref{29} and Theorem~\ref{213} show that the dual $(\Mh,\deh)$ is again a
bicrossed product, obtained by interchanging $\al$ and $\be$. Now defining the isomorphism $u
: G_1 \recht G_2 : u(g) = g^{-1}$, one verifies that $$u(\be_s(g)) = \al_{u^{-1}(s)}(u(g))$$
for all $g,s \in \R \setminus \{0\}$. Hence, interchanging $\al$ and $\be$, we get an
isomorphic matched pair and so an isomorphic l.c. quantum group.
\end{proof}
\begin{proposition} \label{nontriv1}
The l.c. quantum group $(M,\de)$ is not a Kac algebra, it is non-compact,
non-discrete and non-unimodular. The scaling group is non-trivial and the left and right
Haar weights are not traces.
\end{proposition}
\begin{proof}
Proposition~\ref{compact} shows that $(M,\de)$ is non-compact and
non-discrete. The modular functions of $G_1$ and $G_2$ are trivial and
$\sde(a,b)=|a|$. Using the notation introduced in Propositions~\ref{48} and \ref{49} we get
$$P(g,s) = \Bigl|\frac{g}{s(g-1)+1}\Bigr| \quad\text{and}\quad \nab(g,s) = \Bigl|\frac{1}{s(g-1)+1}\Bigr| \; .$$
Because $P$ is non-trivial, the scaling group of $(M,\de)$ is non-trivial, so $(M,\de)$ is not a Kac algebra. Because $\nab$ is non-trivial, the modular
automorphism group of the left Haar weight $\vfi$ is non-trivial. Hence $\vfi$ is not a trace.
Because the right Haar weight $\psi$ on $(M,\de)$ is given by $\psi = \vfi R$, also $\psi$ is
not a trace. Finally, Proposition~\ref{49} gives
$$\sde_M(g,s) = \Bigl|\frac{s(g-1)+1}{gs}\Bigr|$$
and hence $(M,\de)$ is non-unimodular.
\end{proof}

Let us now describe, following the scheme of the previous subsection, the infinitesimal Hopf
algebra of the l.c. quantum group $(M,\de)$. Let $\cA$ be the algebra of polynomials on $G_2$
in the variable $s$, generated by $A(s)=s$ for all $s \in G_2$. Then $\de_2(A)= A \ot A$. Let
$X$ denote the generator of the Lie algebra $\gone$ of $G_1$ given by $$X_x[B] = \frac{d}{dy}
B(xy) \big|_{y=1} \; .$$ From \eqref{comrel} and \eqref{comult}, we get that this
infinitesimal Hopf algebra is generated by elements $A$ and $X$, where $A$ is invertible and
$$ [A,X] = X \acts A \; , \quad \de(A) = A \ot A \quad\text{and}\quad \de(X) = 1 \ot X +
\beh(X) \; . $$ Now we get that $$ (X \acts A)(r) = \frac{d}{dx} \al_x(r) \big|_{x=1} = r(1-r)
\quad\text{and}\quad \beh(X) = X \ot \bigl( r \mapsto \frac{d}{dx} \be_r(x) \big|_{x=1} \bigr)
= X \ot A \; . $$ Hence the infinitesimal Hopf algebra is generated by elements $A$ and $X$,
where $A$ is invertible and $$ [A,X] = A(1-A) \; , \quad \de(A) = A \ot A \quad\text{and}\quad
\de(X) = X \ot A + 1 \ot X \; . $$
\begin{remark}\label{exttriv}
a) From Remark \ref{Liebialg} it follows that, on the level
of Lie bialgebras, in the above situation there is no non-trivial extension. The same is true on the level of Hopf algebras (see \cite{Mas2}, 4.11).

b) The above example is closely related to \cite{Majbook}, 6.2.17,
where $G_1=G_2=(\R,+)$ and the mutual actions are only defined in the neighbourhood of the origin $(0,0)$. This suffices to compute the corresponding infinitesimal Hopf algebra \cite{Majbook}, 6.2.18; the last one coincides with our infinitesimal Hopf algebra.

c) Let $q > 0$ and consider the automorphism of $G_1$ given by $\mu_q(g) = g^q$. Here we define $g^q = -(-g)^q$ whenever $g < 0$. Putting $i_q = i \; \mu_q$ and keeping $j$ as above, we get an isomorphic matched pair of groups, with mutual actions $\al^q$ and $\be^q$ given by
$$\al^q_g(s) = \frac{g^q s}{s(g^q - 1) +1} \quad\text{and}\quad \be^q_s(g) = \bigl( s (g^q - 1) +1 \bigr)^{\frac{1}{q}} \; .$$
If now $q \recht 0$ we observe that for $g > 0$, $\al^q_g(s) \recht s$ and $\be^q_s(g) \recht g^s$. For every $q > 0$ we can construct the bicrossed product l.c.\ quantum group with associated multiplicative unitary $W^q$ on $L^2(G_1) \ot L^2(G_2)$, and they are all isomorphic. Denoting by $P$ the orthogonal projection of $L^2(G_1)$ onto the space of vectors $\xi \in L^2(G_1)$ satisfying $\xi(g) =0$ whenever $g<0$, we get that $(P \ot 1 \ot P \ot 1)W^q(P \ot 1 \ot P \ot 1)$ converges in \strong topology to a multiplicative unitary which is isomorphic (using Fourier transformation) with the multiplicative unitary of the $ax+b$-group. In this sense our example is a deformation of the $ax+b$-group.

This can also be seen on the Hopf algebraic level. Writing $\Xtil = qX$ we see that the infinitesimal Hopf algebra of our example is the Hopf algebra generated by elements $A$ and $\Xtil$, where $A$ is invertible and
$$
[A,\Xtil] = qA(1-A) \; , \quad
\de(A) = A \ot A \quad\text{and}\quad \de(\Xtil) = \Xtil \ot A + 1 \ot \Xtil \; .
$$
If we now formally take the limit $q \recht 0$ we get the Hopf algebra of polynomials on the $ax+b$-group.

d) In \cite{Wor9} S.L.~Woronowicz constructs, on the operator algebra level, another
deformation of the $ax+b$-group. Let $q \in U(1)$. The underlying Hopf algebra of Woronowicz'
example is generated by elements $A$ and $B$, where $A$ is invertible and $$ AB = q BA \; ,
\quad \de(A) = A \ot A \quad\text{and}\quad \de(B) = B \ot A + 1 \ot B \; . $$ Let now $r \in
\R$ and define $\Btil = B + r(1-A)$. Then the Hopf algebra of Woronowicz is generated by $A$
and $\Btil$, where $A$ is invertible and $$ A\Btil = q \Btil A + r(1-q) A(1-A) \; , \quad
\de(A) = A \ot A \quad\text{and}\quad \de(\Btil) = \Btil \ot A + 1 \ot \Btil \; . $$ If we now
take $q = \exp(i \lambda)$ and $r= \frac{1}{\lambda}$ and let $\lambda \recht 0$ we formally
obtain the Hopf algebra generated by $A$ and $\Btil$, where $A$ is invertible and $$ A\Btil =
\Btil A - i A(1-A) \; , \quad \de(A) = A \ot A \quad\text{and}\quad \de(\Btil) = \Btil \ot A +
1 \ot \Btil \; . $$ Identifying $X$ with $i \Btil$ this Hopf algebra is the infinitesimal Hopf
algebra of $(M,\de)$.

It should be remarked that Woronowicz' example cannot be obtained as the cocycle bicrossed product of two l.c.\ groups, because it follows from the work of A.~Van Daele \cite{VD5} that the scaling constant of Woronowicz' example is different from $1$, contradicting Proposition~\ref{48}.

\end{remark}

\subsection{Example: split extension}

Define $H=SL_2(\R)$. Consider $\Z/2\Z$ as a normal subgroup $K$ of $H$ by identifying it with
$\{1,-1\}$ and define $G= H/K$. Then we define
\begin{align*}
G_1 &= \{(a,b) \mid a > 0, b \in \R \} \quad\text{with}\; (a,b)(c,d) = (ac, ad + \frac{b}{c})
\\ G_2 &= (\R,+)
\end{align*}
and
$$i(a,b) = \begin{pmatrix} a & b \\ 0 & \frac{1}{a} \end{pmatrix} \operatorname{mod} K
\quad\quad j(x) = \begin{pmatrix} 1 & 0 \\ x & 1 \end{pmatrix} \operatorname{mod} K \; .$$
This way $(G_1,G_2)$ is a matched pair of groups. A direct computation gives that, whenever
$a+bx \neq 0$
$$\al_{(a,b)}(x) = \frac{x}{a(a+bx)} \quad\quad \be_x(a,b) = \begin{cases} (a+bx,b) \tekst{if}
a+bx > 0 \\ (-a-bx,-b) \tekst{if} a+bx < 0 \end{cases} \; . $$
Taking $\cU=\cV=1$, we can construct the bicrossed product l.c. quantum group $(M,\de)$ and its dual $(\Mh,\deh)$. The following result is similar to
Proposition~\ref{nontriv1}.
\begin{proposition} \label{nontriv2}
The l.c. quantum groups $(M,\de)$ and $(\Mh,\deh)$ are not Kac algebras, they are non-compact,
and non-discrete. The scaling groups are non-trivial. The left and right Haar weights of
$(M,\de)$ and $(\Mh,\deh)$ are not traces. Finally $(M,\de)$ is unimodular, but $(\Mh,\deh)$ is
non-unimodular.
\end{proposition}
\begin{proof}
The proof is completely similar to the one of
Proposition~\ref{nontriv1} and uses again Propositions~\ref{48} and
\ref{49}. Let us mention only that the groups $G$ and $G_2$ are unimodular, $\sde_1(a,b)=\frac{1}{a^2}$ and we get
that
$$\nabh(a,b,s)=P(a,b,s)=\frac{a^2}{(a+bs)^2} \; , \quad \nab(a,b,s) = \frac{1}{a^2(a+bs)^2} \; , \quad
\sde_M=1 \quad\text{and}\quad \sde_{\Mh}(a,b,s) = (a+bs)^4 \; .$$
\end{proof}

Let us now describe the infinitesimal Hopf algebras of $(M,\de)$ and $(\Mh,\deh)$, following
again the strategy of Subsection 5.2. Let $\cA$ be the algebra of polynomials on $G_2$ in the
variable $x$ generated by $A(x)=x$. Take generators $X,Y$ for the Lie algebra $\gone$ of $G_1$
given by $$X_{(a,b)}[B] = \frac{d}{dx} B((a,b)(x,0)) \big|_{x=1} \qquad\qquad Y_{(a,b)}[B] =
\frac{d}{dx} B((a,b)(1,x)) \big|_{x=0}$$ whenever $(a,b) \in G_1$ and when $B$ is a smooth
function on $G_1$. So we observe that $[X,Y]_{\gone} = 2Y$. To determine the infinitesimal
Hopf algebra of $(M,\de)$ we compute $$(X \acts A)(x) = \frac{d}{da} \al_{(a,0)}(x)
\big|_{a=1} = -2x \qquad\qquad (Y \acts A)(x) = \frac{d}{db} \al_{(1,b)}(x) \big|_{b=0} =
-x^2$$ and taking the obvious coordinates on $G_1$ we see that
\begin{align*}
& \beh_1(X)(x) = \frac{d}{da} \be^1_x(a,0) \big|_{a=1} = 1 \qquad\qquad
\beh_2(X)(x) = \frac{d}{da} \be^2_x(a,0) \big|_{a=1} = 0 \\
& \beh_1(Y)(x) = \frac{d}{db} \be^1_x(1,b) \big|_{b=0} = x \qquad\qquad
\beh_2(Y)(x) = \frac{d}{db} \be^2_x(1,b) \big|_{b=0} = 1 \; .
\end{align*}
So we get that
$$X \acts A = -2A \; , \quad Y \acts A = -A^2 \; , \quad \beh(X) = X \ot 1 \quad\text{and} \quad
\beh(Y) = X \ot A + Y \ot 1 \; .$$
So this Hopf algebra has generators $A,X$ and $Y$ with
relations
\begin{align*}
& [X,A] = 2A \; , \quad [Y,A] = A^2 \quad\text{and} \quad [X,Y] = 2Y \\
& \de(A) = A \ot 1 + 1 \ot A \\
& \de(X) = X \ot 1 + 1 \ot X \\
& \de(Y) = Y \ot 1 + X \ot A + 1 \ot Y \; .
\end{align*}
To obtain the infinitesimal Hopf algebra of the dual l.c. quantum group $(\Mh,\deh)$,
interchange $\al$ and $\be$. We take as $\cA$ the algebra of polynomials on $G_1$ in the
variables $a$ and $b$ generated by $A(a,b)=a$ and $B(a,b)=b$. We denote by $X$ the generator
of the Lie algebra $\gtwo$ of $G_2$ given by $$X_x[B] = \frac{d}{dy} B(x + y) \big|_{y=0}$$
whenever $x \in G_1$ and when $B$ is a smooth function on $G_1$. Similar computations as above
yield that $$X \acts A = B \; , \quad X \acts B = 0 \quad\text{and} \quad \beh(X) = X \ot
A^{-2} \; .$$ So this Hopf algebra is generated by elements $X,A$ and $B$, with $A$ invertible
and
\begin{align*}
& [A,X] = B \; , \quad [A,B] = 0 \quad\text{and} \quad [B,X] = 0 \\
& \de(A) = A \ot A \\
& \de(B) = A \ot B + B \ot A^{-1} \\
& \de(X) = X \ot A^{-2} + 1 \ot X \; .
\end{align*}

\begin{remark}
In the same spirit as in Remark~\ref{exttriv}.c) we show that $(M,\de)$ and $(\Mh,\deh)$ are
deformations of usual groups. Let again $q > 0$, and define the l.c.\ group $G_1^q$ on the
same space as $G_1$, but with multiplication given by $(a,b) \; \cdot_q \; (c,d) = (ac, a^q d
+ \frac{b}{c^q})$. Then we have a natural isomorphism $\mu_q : G_1^q \recht G_1$ given by
$\mu_q(a,b) = (a^q,qb)$. Defining $i_q = i \; \mu_q$ and putting $j$ as above we get a matched
pair $G_1^q, G_2$ of l.c.\ groups which is isomorphic with our original matched pair. The
mutual actions are given by $$\al^q_{(a,b)}(x) = \frac{x}{a^q(a^q+qbx)} \quad\quad
\be^q_x(a,b) = \begin{cases}\bigl((a^q+qbx)^{\frac{1}{q}},b\bigr) \tekst{if} a^q+qbx > 0 \\
\bigl((-a^q-qbx)^{\frac{1}{q}},-b \bigr) \tekst{if} a^q+qbx < 0 \end{cases} \; . $$ So we
observe that for $q \recht 0$ we get $\al^q_{(a,b)}(x) \recht x$ and $\be^q_x(a,b) \recht
\bigl( a \exp(bx),b \bigr)$. Having at hand the matched pair $G_1^q, G_2$ we get a bicrossed
product l.c.\ quantum group, isomorphic with $(M,\de)$ and, using the obvious unitary from
$L^2(G_1^q)$ onto $L^2(G_1)$, the associated multiplicative unitaries $W_q$ all can act on the
Hilbert space $L^2(G_1) \ot L^2(G_2)$. Then we see that for $q \recht 0$, $W_q$ converges in
\strong topology to a multiplicative unitary which is isomorphic (using Fourier
transformation) with the multiplicative unitary of the Heisenberg group. In this sense
$(M,\de)$ is a deformation of the Heisenberg group.

Again this can be seen on the Hopf algebra level by putting $\Xtil = qX$ and $\tilde{Y}=qY$,
in the description of the infinitesimal Hopf algebra of $(M,\de)$. Taking formally the limit
$q \recht 0$, we get the Hopf algebra of polynomials on the Heisenberg group.

But there is more. We can also consider the automorphism $\nu_q$ of $G_2$ given by $\nu_q(x)=qx$. Keeping $i$ as above and putting $j_q = j \; \nu_q$, we can proceed in exactly the same way as in the first part of this remark. We get again a family of isomorphic l.c.\ quantum groups, such that the associated dual multiplicative unitaries $\Wh_q$ converge in \strong topology to a multiplicative unitary which is isomorphic with the multiplicative unitary of the following group $H$. The space of $H$ is $\{(a,b,x) \mid a > 0, b,x \in \R \}$, and $(a,b,x) \cdot (c,d,y) = (ac, ad + \frac{b}{c}, \frac{x}{c^2} + y)$.
Hopf algebraically we obtain the same result by putting $\Xtil=qX$ in the description of the infinitesimal Hopf algebra for $(\Mh,\deh)$ and then taking $q \recht 0$. Hence $(\Mh,\deh)$ is a deformation of some generalized $ax+b$-group.
\end{remark}

\subsection{Example: non-trivial extensions}

Let $(G_1,G_2)$ and $\al,\be$ be as in the previous
subsection. We look for a cocycle $\cU$ and Remark \ref{Liebialg} suggests to take $\cV=1$. Referring to Lemma~\ref{44} we look for a measurable function $\cU$ on $G_1 \times G_1 \times G_2$ such that
\begin{align*}
\cU(g,h,\al_k(s)) \; \cU(gh,k,s) &= \cU(h,k,s) \; \cU(g,hk,s), \\
\cU(g,h,s) \; \cU(\be_{\al_h(s)}(g),\be_s(h),t) &= \cU(g,h,t+s)
\end{align*}
almost everywhere. Define a function $A(\cdot)$ by
$$\cU(g,h,s) = \exp(i A(g,h,s)) \; .$$
So $A(\cdot)$ should satisfy
\begin{align}
A(g,h,\al_k(s)) + A(gh,k,s) &= A(h,k,s) + A(g,hk,s) \quad(\text{mod} 2 \pi), \label{ster1} \\
A(g,h,s) + A(\be_{\al_h(s)}(g),\be_s(h),t) &= A(g,h,t+s) \quad(\text{mod} 2 \pi) \label{ster2}
\end{align}
almost everywhere. For $s \in G_2$ and $g,h \in G_1$, we define $$\phi_s(g,h) =
(\be_{\al_h(s)}(g),\be_s(h)) \; .$$ Then $\phi_{t+s} = \phi_t \phi_s$ almost everywhere, and
Eq.~\eqref{ster2} becomes
\begin{equation} \label{ster3}
A(g,h,s)+A(\phi_s(g,h), t) = A(g,h,t+s) \quad(\text{mod} 2 \pi) \; .
\end{equation}
It is natural to look for a solution of this equation of the form
\begin{equation} \label{ster4}
A(g,h,s) = P \int_0^s f(\phi_r(g,h)) \; dr \; ,
\end{equation}
with a function $f$ on $G_1 \times G_1$ such that for
almost all $g,h \in G_1$ the function
$$r \mapsto f(\phi_r(g,h))$$
has a principal value integral over any interval in $\R$. So in
order to fulfill Eq.~\eqref{ster1}, our function $f$ should satisfy, for almost all $g,h,k
\in G_1$ and $s \in G_2$
\begin{equation} \label{vglf}
P \int_0^{\al_k(s)} f(\phi_r(g,h)) \; dr + P \int_0^s f(\phi_r(gh,k)) \; dr =
P \int_0^s f(\phi_r(h,k)) \; dr + P \int_0^s f(\phi_r(g,hk)) \; dr \quad(\text{mod} 2 \pi) \; .
\end{equation}
Differentiating to $s$ we should look for a function $f$ satisfying
$$ \frac{d}{ds} \al_k(s) \; f(\tilde g,\tilde h)+f(\tilde g\tilde h,
\tilde k)=f(\tilde h, \tilde k) + f(\tilde g, \tilde h \tilde k) $$
where $\tilde g= \beta_{\alpha_{hk}(s)}(g),\tilde h=
\beta_{\alpha_k(s)}(h),\tilde k=\beta_{s}(k)$. Observing that
$$\al_{\be_s(k)}(t) + \al_k(s) = \al_k(t+s)$$
and differentiating to $t$ at $t=0$ we get that
$$\frac{d}{ds} \al_k(s) = \frac{d}{dt} \al_{\tilde k}(t) \big|_{t=0} \; .$$
Hence, $f$ should satisfy
$$\frac{d}{dt} \al_k(t) \big|_{t=0} \; f(g,h)+f(gh,k)=f(h,k) + f(g,hk)$$
for almost all $g,h,k \in G_1$. Putting $g=(a,b)$, $h=(c,d)$ and $k=(l,m)$ we see that $f$ should satisfy
\begin{equation} \label{star1}
\frac{1}{l^2} f(a,b,c,d) + f(ac,ad+\frac{b}{c},l,m) = f(c,d,l,m) + f(a,b,cl,cm+\frac{d}{l})
\end{equation}
where $f$ is defined on $\R_0^+ \times \R \times \R_0^+ \times \R$.
\begin{proposition}\label{funceq}
A general sufficiently smooth solution $f(a,b,c,d)$ of \eqref{star1} is given by
$$f(a,b,c,d)=\lambda \frac{b \log c}{a c^2} + \frac{1}{c^2} B(a,b) + B(c,d) -B(ac,ad+\frac{b}{c})$$
where $\lambda \in \R$ and $B$ is a smooth function on $G_1$.
\end{proposition}
\begin{proof}
First, one can compute directly that
\begin{equation} \label{star2}
f(a,b,c,d) = \frac{1}{c^2} B(a,b) + B(c,d) - B(ac,ad + \frac{b}{c})
\end{equation}
where $B$ is an arbitrary function on $G_1$, is a solution of Eq.~\eqref{star1}. We will call
such a solution trivial. Because the solutions of \eqref{star1} form a linear space, we will
only determine a general solution of \eqref{star1} modulo a trivial solution.

In the course of the proof we use the subscript notations $f_2$, $f_{13}$, $\ldots$ to denote
the partial derivative of $f$ to the first variable and to the first and third variable
respectively. Take the derivative to $b$ of a sufficiently smooth solution $f$ of
Eq.~\eqref{star1} and evaluate it in $a=1$ and $b=0$. Then we get $$f_2(c,d,l,m) = c \,
h(cl,cm+\frac{d}{l}) - \frac{c}{l^2} \, h(c,d)$$ where $h(c,d) = f_2(1,0,c,d)$. If $H(\cdot)$
is a smooth function such that $H_2(c,d) = h(c,d)$, this gives $$f(c,d,l,m) = cl \,
H(cl,cm+\frac{d}{l}) - \frac{c}{l^2} \, H(c,d) + k(c,l,m)$$ where $k$ is some smooth function.
Hence, modulo a trivial solution, $$f(c,d,l,m) = l \, H(l,m) + k(c,l,m),$$ i.e., $f$ does not
depend on its second variable: $$f(a,b,c,d) = g(a,c,d) \; .$$ Now Eq.~\eqref{star1} gives
$$\frac{1}{l^2} g(a,c,d) + g(ac,l,m) = g(c,l,m) + g(a,cl,cm+\frac{d}{l}) \; .$$ Taking the
derivative to $a$ and evaluating it in $a=1$ we get
\begin{equation} \label{star4}
g_1(c,l,m) = \frac{1}{c} \, g_1(1,cl,cm + \frac{d}{l}) - \frac{1}{cl^2} \, g_1(1,c,d) \; .
\end{equation}
Taking the derivative to $d$ it follows that
\begin{equation} \label{star3}
\frac{1}{cl} \, g_{13}(1,cl,cm+\frac{d}{l}) = \frac{1}{cl^2} \, g_{13}(1,c,d) \; .
\end{equation}
With $l=1$ we get that $g_{13}(1,c,cm+d) = g_{13}(1,c,d)$ and hence $g_{13}(1,c,d) = q(c)$ for
some smooth function $q(\cdot)$. Plugging this into Eq.~\eqref{star3} we get $$\frac{1}{cl}
q(cl) = \frac{1}{cl^2} q(c)$$ and so there exists a number $\lambda \in \R$ such that $q(c) =
\frac{\lambda}{c}$. Then it follows that $$g_1(1,c,d) = \frac{\lambda d}{c} + r(c)$$ for some
smooth function $r(\cdot)$. Using Eq.~\eqref{star4} we have $$g_1(c,l,m) = \frac{\lambda
m}{cl} + \frac{1}{c} \, r(cl) - \frac{1}{cl^2} \, r(c) \; .$$ If $R(\cdot)$ is a primitive for
the function $\frac{1}{x} r(x)$ we get $$g(c,l,m) = \frac{\lambda m}{l} \log c + R(cl) -
\frac{1}{l^2} R(c) + t(l,m)$$ for some smooth function $t(\cdot)$. Since $f(a,b,c,d) =
g(a,c,d)$, we have, modulo a trivial solution, $$f(c,d,l,m) = \frac{\lambda m}{l} \log c +
t(l,m) + R(l)$$ or $$f(c,d,l,m) = \frac{\lambda m}{l} \log c + v(l,m)$$ for some smooth
function $v(\cdot)$. Plugging this into Eq.~\eqref{star1} we get $$\frac{1}{l^2} v(c,d) =
v(cl,cm + \frac{d}{l}) \; .$$ With $l=1$ it follows that $v$ does not depend on its second
variable. Then there is a number $\nu \in \R$ such that $v(c,d) = \frac{\nu}{c^2}$. But also
$f(a,b,c,d) = \frac{\nu}{c^2}$ is a trivial solution with $B(a,b) = \nu$. Hence we may
conclude that $$f(a,b,c,d) = \lambda \frac{d}{c} \log a$$ modulo a trivial solution, where
$\lambda \in \R$. One can check directly that this is indeed a solution.

It will be more convenient, taking $B_0(a,b) =\lambda \frac{b}{a} \log a$, to get a general
sufficiently smooth solution of Eq.~\eqref{star1} of the form $$f(a,b,c,d) = \lambda \frac{b
\log c}{a c^2} + \frac{1}{c^2} B(a,b) + B(c,d) - B(ac,ad+\frac{b}{c})$$ where $\lambda \in \R$
and $B(\cdot)$ is an arbitrary smooth function.
\end{proof}
Since two solutions differing with a trivial solution, give rise to cohomologous cocycles, we
can work on with the solution $$f_\lambda(a,b,c,d) = \lambda \frac{b \log c}{a c^2}$$ where
$\lambda \in \R$ is some parameter. Making an easy computation one observes that
$$f_\lambda(\phi_r(a,b,c,d)) = \frac{\lambda b}{(c+dr)(ac + (ad + \frac{b}{c})r)} \log |c+dr|
\; .$$ As a function of $r$ we indeed have a principal value integral over any interval in
$\R$, for almost all $a,b,c,d$. One verifies that $$P \int_{-\infty}^{+\infty}
f_\lambda(\phi_r(a,b,c,d)) \; dr = \begin{cases} \frac{\lambda}{2} \pi^2 \quad\text{if}\;
\frac{d}{b}(ad + \frac{b}{c}) > 0 \\ -\frac{\lambda}{2} \pi^2 \quad\text{if}\; \frac{d}{b}(ad
+ \frac{b}{c}) < 0 \end{cases} \; .$$ Define $A_\lambda(\cdot)$ with $f_\lambda(\cdot)$ by
\eqref{ster4} and check Eq.~\eqref{vglf}. When $s=0$ this equation is trivially true. When $s$
grows, the equation remains true, because the differentiated equation is fulfilled. Also,
passing a pole is not a problem because the substitution rule can be applied to this principal
value integrals. The only problem arises when $s$ passes the pole where $\al_k(s)$ grows to
infinity. Then the principal value integral $$P \int_0^{\al_k(s)} f_\lambda(\phi_r(g,h)) \;
dr$$ changes to
\begin{align*}
P \int_{0 \recht +\infty, -\infty \recht \al_k(s)}f_\lambda(\phi_r(g,h)) \; dr &=
P \int_0^{\al_k(s)} f_\lambda(\phi_r(g,h)) \; dr + P \int_{-\infty}^{+\infty} f_\lambda(\phi_r(g,h)) \; dr \\ &=
P \int_0^{\al_k(s)} f_\lambda(\phi_r(g,h)) \; dr \pm \frac{\lambda}{2} \pi^2 \; .
\end{align*}
Hence, Eq.~\eqref{vglf} remains true iff $\lambda$ has the form
$$\lambda=\frac{4n}{\pi} \tekst{with} n \in \Z \; .$$
Thus, for any $n \in \Z$, the functions
$$\cU_n(g,h,s) = \exp(i A_{\frac{4n}{\pi}}(g,h,s))$$
and $\cV=1$ satisfy \eqref{cocU}, so
they give the necessary data to perform the cocycle bicrossed product construction. Let us describe the corresponding l.c. quantum group $(M_n,\de_n)$ and its dual $(\Mh_n,\deh_n)$.

Clearly $(M_0,\de_0)$ and $(\Mh_0,\deh_0)$ are precisely the examples
studied in the previous subsection, because $\cU_0=1$. Thanks to Propositions~\ref{48} and
\ref{49}, Proposition~\ref{nontriv2} remains valid for $(M_n,\de_n)$
and $(\Mh_n,\deh_n)$, also when $n \neq 0$. Let us mention the relations between $(M_n,\de_n)$ and $(M_0,\de_0)$, and between $(\Mh_n,\deh_n)$ and
$(\Mh_0,\deh_0)$. Represent $\Mh_n$ on $L^2(G_1 \times G_2)$ as in Proposition~\ref{29}.
Because $\cV=1$, we get $\Mh_n = \Mh_0$ for all $n \in \Z$. Further we have
$$\Wh_n^* = \cC_n \Wh_0^* \quad\text{where}\quad \cC_n = (1 \ot (\io \ot \al)(\Wh_2^*)) \;
(\be \ot \io \ot \io)(\cU_n) \; (1 \ot (\io \ot \al)(\Wh_2)) \; .$$
An easy computation shows that $\cC_n$ is the multiplication operator with the function
$$\cC_n(g,s,h,t) = \cU_n(\be_{\al_h(t)s}(g),h,t) \; .$$
It is also clear that
$$\deh_n(z) = \cC_n \deh_0(z) \cC_n^* \tekst{for all} z \in \Mh_n=\Mh_0 \; .$$
From Proposition~\ref{29} it follows that the left Haar weight $\vfih_n$ equals $\vfih_0$.
Propositions~\ref{48} and \ref{49} imply that $\Jh_n =\Jh_0$ for all $n \in \Z$ and they also
give us that the operators $\nab_n,\nabh_n,P_n,\sde_{M_n}$ and $\sde_{\Mh_n}$ do not depend on
$n$. In particular it follows that $\tauh_n = \tauh_0$, but we see that $J_n \neq J_0$ and
then also $\Rh_n \neq \Rh_0$ for $n \neq 0$.

As in the previous examples we want to compute the infinitesimal Hopf algebras  of
$(M_n,\de_n)$ and $(\Mh_n,\deh_n)$. We start with $(M_n,\de_n)$ and we follow the strategy of
Subsection 5.2. In the previous subsection we already defined the algebra $\cA$ of polynomials
on $G_2$, with generator $A$, and we defined the generators $X$ and $Y$ for the Lie algebra
$\gone$ of $G_1$. If $\cH$ denotes the universal enveloping algebra of $\gone$, we have to
look at the map
\begin{equation} \label{eq.chi}
\hchi_n: \cH \ot \cH \recht A : \hchi_n(H \ot G)(r) = (H_e \ot G_e)[\cU_n^*(\cdot,\cdot,r)] \;
.
\end{equation}
An easy computation yields that, for all $\lambda \in \R$ and $r \in G_2$, we have $$(X_e \ot
Y_e)[f_\lambda(\phi_r(\cdot,\cdot))] = 0 \quad\text{and}\quad (Y_e \ot X_e)
[f_\lambda(\phi_r(\cdot,\cdot))]= \lambda$$ and then it follows that $$\hchi_n(X \ot Y - Y \ot
X) = -i \frac{4n}{\pi} A \; .$$ So, the infinitesimal Hopf algebra of $(M_n,\de_n)$ has
generators $X,Y$ and $A$ and relations
\begin{align*}
& [X,A] = 2A \; , \quad [Y,A] = A^2 \quad \text{and} \quad [X,Y] = 2Y - i \frac{4n}{\pi} A \\
& \de(A) = A \ot 1 + 1 \ot A \\
& \de(X) = X \ot 1 + 1 \ot X \\
& \de(Y) = Y \ot 1 + X \ot A + 1 \ot Y \; .
\end{align*}

To obtain the infinitesimal Hopf algebra of $(\Mh_n,\deh_n)$ we have to interchange $\al$ and
$\be$ and we have to use the cocycle $\cV_n(s,g,h)= \cU_n(h,g,s)$ and take $\cU=1$. Now we
take for $\cA$ the algebra of polynomials on $G_1$ in the variables $a,\log a$ and $b$
generated by $A(a,b)=a$, $(\log A)(a,b)=\log a$ and $B(a,b)=b$. As before we denote by $X$ the
generator of the Lie algebra $\gtwo$ of $G_2$. If $\cH$ denotes the universal enveloping
algebra of $\gtwo$ we have to look at the map
\begin{equation} \label{eq.psi}
\Psi_n : \cH \recht \cA \ot \cA : \Psi_n(H)(g,h) = H_e[\cV_n(\cdot,h,g)] \; .
\end{equation}
So we immediately get that $$\Psi_n(X)(g,h) = \frac{d}{ds} \cU_n(g,h,s) \big|_{s=0} = i
f_{\frac{4n}{\pi}}(g,h)$$ and hence $$\Psi_n(X) = i \frac{4n}{\pi} \; A^{-1} B \ot A^{-2} \log
A \; .$$ The infinitesimal Hopf algebra of $(\Mh_n,\deh_n)$ has generators $A,\log A,B$ and
$X$ with $A$ invertible and with relations
\begin{align*}
& B \tekst{is central}, \quad [A,\log A] = 0 \; , \quad [A,X]=B \quad\text{and}\quad
[\log A,X] = A^{-1} B \\
& \de(A) = A \ot A \; , \quad \de(\log A) = \log A \ot 1 + 1 \ot \log A \\
& \de(B) = A \ot B + B \ot A^{-1} \\
& \de(X) = X \ot A^{-2} + 1 \ot X + i \frac{4n}{\pi} \; A^{-1}B \ot A^{-2} \log A \; .
\end{align*}
\begin{remark}
One can show, using the Lie bialgebra version of the G.I.~Kac' exact sequence \cite{Kac},
(3.14) obtained by Masuoka (see \cite{Mas2}, 2.10), that the group of extensions for the above
example, on the level of Lie bialgebras or Hopf algebras, is exactly $\R$. In Eq.
\eqref{eq.chi} and \eqref{eq.psi} above we see that the maps $\hchi$ and $\Psi$ can be defined
for an arbitrary real parameter $\lambda$, and also the infinitesimal Hopf algebras make sense
for an arbitrary real parameter $\lambda$ instead of $\frac{4n}{\pi}$. But, as we have seen,
on the level of l.c. quantum groups only for the extensions with $\lambda=\frac{4n}{\pi} \ (n
\in \Z)$, the cocycles make sense on the operator algebra level.
\end{remark}

\end{document}